\documentclass[11pt]{article}
\usepackage{mathrsfs}
\usepackage{latexsym,bm}
\usepackage{color}
\usepackage{amssymb, graphicx, amsmath}  
\usepackage{CJK}
\usepackage{indentfirst}

\def \[{\begin{equation}}
\def \]{\end{equation}}

\newtheorem{theorem}{Theorem}[section]

\newtheorem{assumption}{Assumption}[section]

\newtheorem{lemma}{Lemma}[section]

\newtheorem{corollary}{Corollary }[section]

\begin{document}

\begin{CJK*}{GBK}{SONG}
\fontsize{11}{11}

\begin{center}
{\huge \bf The global existence of the smoothing solution for the
Navier-Stokes equations }
  \\[0.5cm]

\medskip

  {\bf  Wang Jianfeng}\\
   Department of Mathematics, Nanjing University, Nanjing, 210093,
   \\ P.
R. China.

\end{center}

\bigskip
{\narrower \noindent  {\bf Abstract.} This paper discussed the
global existence of the smoothing solution for the Navier-Stokes
equations. At first, we construct the theory of the linear equations
which is about the unknown four variables functions with constant
coefficients. Secondly, we use this theory to convert the
Navier-Stokes equations into the simultaneous of the first order
linear partial differential equations with constant coefficients and
the quadratic equations. Thirdly, we use the Fourier transformation
to convert the first order linear partial differential equations
with constant coefficients into the linear equations, and we get the
explicit general solution of it. At last, we convert the quadratic
equations into the integral equations or the question to find the
fixed-point of a continuous mapping. We use the theories about the
Poisson's equation, the heat-conduct equation, the Schauder
fixed-point theorem to prove that the fixed-point is exist, hence
the smoothing solution for the Navier-Stokes equations is globally
exist.

    \par }
\bigskip \hrule
 \section{Introduction}We consider the dynamical equations for an viscous and
incompressible fluid just as follows. {
\begin{eqnarray}
 &&div u=0 \\&&\cfrac{\partial u}{\partial
 t}-\mu \Delta u+\sum\limits_{k=1}^{3}u_{k}\cfrac{\partial u}{\partial
 x_{k}}+\cfrac{1}{\rho}\
 grad p=F\end{eqnarray}}
where\ $u=(u_{1},\ u_{2},\ u_{3})$\ is the velocity vector, it is
the macro motion velocity of the material particle.\ $\mu$\ is the
dynamic viscosity coefficient, we assume it is a constant.\ $\rho$\
is the density of the material particle, it is the mass of per
volume fluid. Because the fluid is incompressible, we can assume\
$\rho$\ is a constant, too.\ $p$\ is intensity of the pressure, it
is the pressure on per area fluid, the direction is perpendicular to
such area.\ $F=(F_{1},\ F_{2},\ F_{3})$\ is the density of the body
force, it is the external force of per unit mass.\ $u=(u_{1},\
u_{2},\ u_{3})$\ and\ $p,\ F$\ are all functions on the variables of
the time\ $t$\ and position\ $x=(x_{1},\ x_{2},\ x_{3})$\ . We
assume\ $\tau=\cfrac{1}{\rho}\ ,$\ $\tau$\ is called the specific
volume of the fluid, it is just the volume of per unit mass. These
equations are the second order partial differential equations about
four unknown functions:\ $u=(u_{1},\ u_{2},\ u_{3})$\ ,\ $p$\ , and
they are completed. These equations are called the Navier-Stokes
equations. \\In order to discuss it more conveniently, we often
rewrite the Navier-Stokes equations as follows.
 {
\begin{eqnarray}&&\cfrac{\partial u_{1}}{\partial
 x_{1}}+\cfrac{\partial u_{2}}{\partial
 x_{2}}+\cfrac{\partial u_{3}}{\partial
 x_{3}}=0 \\&&\cfrac{\partial u_{j}}{\partial
 t}-\mu(\cfrac{\partial^{2} u_{j}}{\partial
 x_{1}^{2}}+\cfrac{\partial^{2} u_{j}}{\partial
 x_{2}^{2}}+\cfrac{\partial^{2} u_{j}}{\partial
 x_{3}^{2}}\
 )+\sum\limits_{k=1}^{3}u_{k}\cfrac{\partial u_{j}}{\partial
 x_{k}}+\tau \cfrac{\partial p}{\partial
 x_{j}}=F_{j}\ ,\  j=1\
 ,\ 2\
 ,\ 3\
 .
\end{eqnarray}}And we assume as follows.\begin{assumption}
\label{Assumption}(1)We only discuss the Navier-Stokes equations on
the region as follows,\\$t\in [0,\ T],\ T$\ is given,\ $(x_{1},\
x_{2},\ x_{3})^{T}\in K_{1}$\ ,\ $K_{1}$\ is a bounded and closed
set in\ $R^{3}$\ , and when\ $t$\ is not in\ $[0,\ T]$, or\
$(x_{1},\ x_{2},\ x_{3})^{T}$\ is not in\ $K_{1}$\ ,\ $u_{j}\equiv0\
,\ p\equiv0\ ,\ F_{j}\equiv0\ ,\ j=1\
 ,\ 2\
 ,\ 3\
 .$\\(2)In the region\ $K_{1}^{\prime}=[0 ,\ T]\times K_{1}$\ ,\ $u_{j}\in C^{2}\ ,\ p\in
C^{1}\ ,\ F_{j}\in C^{1}\ ,\ j=1\
 ,\ 2\
 ,\ 3\
 .$\\(3)The boundary of\ $K_{1}$\ ,\ $\partial K_{1}$\ satisfies the exterior ball condition,
\ $\forall\ (x_{1},\ x_{2},\ x_{3})^{T}\in \partial K_{1}$\ , there
exists a ball\ $B_{\rho}(y)\subset R^{3}\setminus K_{1}^{0}$\ , such
that\ $B_{\rho}(y)\cap K_{1}=(x_{1},\ x_{2},\ x_{3})^{T}$\ , where\
$y=(y_{1},\ y_{2},\ y_{3})^{T}\ ,\ B_{\rho}(y)=\{z\in R^{3}|\
|z-y|\leq \rho\ ,\ \rho> 0\}$\ ,\ $K_{1}^{0}$\ is the interior of\
$K_{1}$\ .\\(4)\ $K_{1}$\ and\ $F_{j}\ ,\ j=1\
 ,\ 2\
 ,\ 3\
 ,$\ satisfy the following,$$\max_{1\leq j\leq 9}\{\
|F^{-1}(\alpha_{j1}^{T}Y_{1})|\ ,|F^{-1}(\alpha_{j2}^{T}Y_{1})|\
,|F^{-1}(\alpha_{j1}^{T}Y_{1})F^{-1}(\alpha_{j2}^{T}Y_{1})|\}\leq\cfrac{
\theta(1-\theta)}{M_{T,\ 3}M(K_{1})(2\theta+M_{T,\ 3}M(K_{1}))}\
,$$where\ $0<\theta<1$\ , and\ $F^{-1}(\alpha_{j1}^{T}Y_{1})\ ,\
F^{-1}(\alpha_{j2}^{T}Y_{1})\ ,\ M_{T,\ 3}\ ,\ M(K_{1})$\ will be
defined in the section 3.
 \end{assumption}After that we will show that the smoothing solution for the Navier-Stokes
equations is globally exist in the region\ $K_{1}^{\prime}$\ .
\section{ Main conclusion }\setcounter{equation}{0}
\begin{theorem} \label{Theorem2-1}We consider the linear
equations as follows.
$$AX=\beta\
,$$ where\ $A=(a_{ij})_{m\times s}\in R^{m\times s}$\ is a constant
matrix,$$X=\left(X_{1}(x_{1},\ x_{2},\ x_{3},\ t),\ X_{2}(x_{1},\
x_{2},\ x_{3},\ t),\ \cdots,\ X_{s}(x_{1},\ x_{2},\ x_{3},\
t)\right)^{T}$$ is the unknown\
 $s$\
dimensional four variables functional vector,
$$\beta=\left(b_{1}(x_{1},\ x_{2},\ x_{3},\ t),\ b_{2}(x_{1},\ x_{2},\ x_{3},\ t),\ \cdots,\ b_{m}(x_{1},\ x_{2},\ x_{3},\ t)\right)^{T}$$
 is the known
\ $m$\ dimensional four variables functional
vector.\\
We can get the following conclusions.\\(1)A necessary and sufficient
condition for the existence of the solution of this equations
is$$\forall\ (x_{1},\ x_{2},\ x_{3},\ t)^{T}\in R^{4}\ ,\ rank(A,\
\beta(x_{1},\ x_{2},\ x_{3},\ t))=rank(A)\ ,\ $$ where\
$\beta(x_{1},\ x_{2},\ x_{3},\ t)$\ is the value of\ $\beta$\ , when
\ $(x_{1},\ x_{2},\
x_{3},\ t)^{T}$\ is given.\\
(2)If the solution of this equations is exist,\ $rank(A)=r,\ r<s$\
,\ $X_{0}$\ is a particular solution of this equations, and the
constant linearly independent solutions of the equations\ $AX=0$\
are\\ $\eta_{1}\ ,\ \eta_{2}\ ,\cdots,\ \eta_{s-r}$\ , we denote\
$A(\eta)=[\eta_{1}\ ,\ \eta_{2}\ ,\cdots,\ \eta_{s-r}]$\ , then its
general solution can be expressed as$$X=X_{0}+A(\eta)Z_{1}\ ,$$
where\ $Z_{1}$\ is
  the arbitrary\ $s-r$\ dimensional four variables functional
  vector.\end{theorem}
\begin{theorem} \label{Theorem2-2}If\ $u_{1},\ u_{2},\ u_{3}\in C^{2}\ ,\ p\in
C^{1}\ ,\ F_{1},\ F_{2},\ F_{3}\in C^{1}$\ , the Navier-Stokes
equations can be converted into the simultaneous of the first order
linear partial differential equations with constant coefficients and
the quadratic equations, just as the following.\[\label
{A-Problem}\begin{cases}\cfrac{\partial (HZ+h)}{\partial
t}=H_{0}Z+h_{0}\\\cfrac{\partial (HZ+h)}{\partial x_{i}}=H_{i}Z\ ,\
i=1,\ 2,\ 3.\end{cases}
\]
\[\label {A-Problem}\begin{cases}e_{47}^{T}Z=(e_{3}^{T}Z)(-\alpha_{1}^{T}Z),\
e_{48}^{T}Z=(e_{7}^{T}Z)(e_{1}^{T}Z),\
e_{49}^{T}Z=(e_{11}^{T}Z)(e_{2}^{T}Z),\\
e_{50}^{T}Z=(e_{3}^{T}Z)(e_{4}^{T}Z),\
e_{51}^{T}Z=(e_{7}^{T}Z)(e_{5}^{T}Z),\
e_{52}^{T}Z=(e_{11}^{T}Z)(e_{6}^{T}Z),\\
e_{53}^{T}Z=(e_{3}^{T}Z)(e_{8}^{T}Z),\
e_{54}^{T}Z=(e_{7}^{T}Z)(e_{9}^{T}Z),\
e_{55}^{T}Z=(e_{11}^{T}Z)(e_{10}^{T}Z)\ ,
\end{cases}
\]where\begin{eqnarray*}&&H=(-\alpha_{1},\ -\alpha_{2},\ -\alpha_{3},\ -\alpha_{4},\
e_{1},\ e_{2},\ e_{3},\ e_{4},\ e_{5},\ e_{6},\ e_{7},\ e_{8},\
e_{9},\ e_{10},\ e_{11},\ e_{16})^{T}\ ,\\&&H_{0}=(e_{21},\ e_{17},\
e_{27},\ e_{37},\ e_{22},\ e_{23},\ -\alpha_{2},\ e_{31},\ e_{32},\
e_{33},\ -\alpha_{3},\ e_{41},\ e_{42},\ e_{43},\ -\alpha_{4},\
e_{12})^{T}\ ,\\&& H_{1}=(e_{18},\ e_{21},\ e_{31},\ e_{41},\
e_{24},\ e_{25},\ -\alpha_{1},\ e_{28},\ e_{34},\ e_{35},\ e_{4},\
e_{38},\ e_{44},\ e_{45},\ e_{8},\ e_{13})^{T}\ ,\\&&
H_{2}=(e_{24},\ e_{22},\ e_{32},\ e_{42},\ e_{19},\ e_{26},\ e_{1},\
e_{34},\ e_{29},\ e_{36},\ e_{5},\ e_{44},\ e_{39},\ e_{46},\
e_{9},\ e_{14})^{T}\ ,\\&& H_{3}=(e_{25},\ e_{23},\ e_{33},\
e_{43},\ e_{26},\ e_{20},\ e_{2},\ e_{35},\ e_{36},\ e_{30},\
e_{6},\ e_{45},\ e_{46},\ e_{40},\ e_{10},\ e_{15})^{T}\
,\\&&\alpha_{1}=(0_{4},\ 1,\ 0_{4},\  1,\  0,\  0,\  0,\  0,\ 0,\
0,\ 0_{39})^{T},
                                    \\&& \alpha_{2}=(0_{12},\  \tau,\  0_{4},\  -\mu,\ -\mu,\  -\mu,\
0_{26},\  1,\  1,\  1,\  0_{6} )^{T},\\&& \alpha_{3}=(0_{13},\
\tau,\   0_{13},\   -\mu,\   -\mu,\   -\mu,\   0_{19},\   1,\   1,\
1,\  0_{3} )^{T},\\&&  \alpha_{4}=(0_{14},\   \tau,\   0_{22},\
-\mu,\   -\mu,\   -\mu,\   0_{11},\   0,\   1,\ 1,\   1)^{T}\
,\end{eqnarray*} and\ $e_{i}$\ is the\ $ith$\ $55$\ dimensional unit
coordinate vector,\ $1 \leq i \leq 55$\ ,\ $h_{0}=(0_{6},\ F_{1},\
0_{3},\ F_{2},\ 0_{3},\ F_{3},\ 0)^{T}\ ,$\ $h=(0,\ F_{1},\ F_{2},\
F_{3},\ 0_{12})^{T}\ ,$\ $Z$\ is an unknown\ $55$\ dimensional
functional vector. Moreover$$(e_{3},\ e_{7},\ e_{11},\
e_{16})^{T}Z=(u_1,\ u_2,\ u_3,\ p)^{T}$$ is the solution of the
Navier-Stokes equations.\end{theorem}
\begin{theorem} \label{Theorem2-3}We consider the first order linear partial differential
equations with constant coefficients and the quadratic equations in
the theorem 2.2, moreover when\ $t$\ is not in\ $[0,\ T]$\ , or\
$(x_{1},\ x_{2},\ x_{3})^{T}$\ is not in\ $K_{1}$\ ,\ $Z\equiv0$\ ,
then we can get the following conclusions.\\(1)We can use the
Fourier transformation to convert\ (2.1)\ into the linear equations
as follows.\[\label {A-Problem}\left(
             \begin{array}{c}
               i\xi_{0}H-H_{0}\\
               i\xi_{1}H-H_{1} \\
               i\xi_{2}H-H_{2}\\
               i\xi_{3}H-H_{3} \\
             \end{array}
           \right) F(Z)=\left(
                         \begin{array}{c}
                           F(h_{0})-i\xi_{0}F(h) \\
                           -i\xi_{1}F(h) \\
                           -i\xi_{2}F(h) \\
                           -i\xi_{3}F(h) \\
                         \end{array}
                       \right)\ ,\]
                       where\begin{eqnarray*}&&F(Z)=\int_{R^{4}}Ze^{-i\xi_{0}t-i\sum_{j=1}^{3}\xi_{j}x_{j}}dtdx_{1}dx_{2}dx_{3}\
                       ,\\&& F(h)=\int_{R^{4}}he^{-i\xi_{0}t-i\sum_{j=1}^{3}\xi_{j}x_{j}}dtdx_{1}dx_{2}dx_{3}\
                       ,\\
                       &&F(h_{0})=\int_{R^{4}}h_{0}e^{-i\xi_{0}t-i\sum_{j=1}^{3}\xi_{j}x_{j}}dtdx_{1}dx_{2}dx_{3}\
                       .\end{eqnarray*}We can get the explicit
                       general solution of it,
                       $$F(Z)=Y_{1}+A_{1}(\eta)Z_{1}\ , \mbox{or}\
                       Z=F^{-1}(Y_{1}+A_{1}(\eta)Z_{1})\ ,$$
                       where\begin{eqnarray*}Y_{1}&=&(i\xi_{2}y_{3},\ i\xi_{3}y_{3},\
y_{3},\ i\xi_{1}y_{7},\ i\xi_{2}y_{7},\ i\xi_{3}y_{7},\ y_{7},\
i\xi_{1}y_{11},\ i\xi_{2}y_{11},\ i\xi_{3}y_{11},\ y_{11},\\&&
i\xi_{0}y_{16},\ i\xi_{1}y_{16},\ i\xi_{2}y_{16},\ i\xi_{3}y_{16},\
y_{16},\ (i\xi_{0})^{2}y_{3},\ (i\xi_{1})^{2}y_{3},\
(i\xi_{2})^{2}y_{3},\ (i\xi_{3})^{2}y_{3},\
i\xi_{0}i\xi_{1}y_{3},\\&& i\xi_{0}i\xi_{2}y_{3},\
i\xi_{0}i\xi_{3}y_{3},\ i\xi_{1}i\xi_{2}y_{3},\
i\xi_{1}i\xi_{3}y_{3},\ i\xi_{2}i\xi_{3}y_{3},\
(i\xi_{0})^{2}y_{7},\ (i\xi_{1})^{2}y_{7},\ (i\xi_{2})^{2}y_{7},\
(i\xi_{3})^{2}y_{7},\\&& i\xi_{0}i\xi_{1}y_{7},\
i\xi_{0}i\xi_{2}y_{7},\ i\xi_{0}i\xi_{3}y_{7},\
i\xi_{1}i\xi_{2}y_{7},\ i\xi_{1}i\xi_{3}y_{7},\
i\xi_{2}i\xi_{3}y_{7},\ (i\xi_{0})^{2}y_{11},\
(i\xi_{1})^{2}y_{11},\ (i\xi_{2})^{2}y_{11},\\&&
(i\xi_{3})^{2}y_{11},\ i\xi_{0}i\xi_{1}y_{11},\
i\xi_{0}i\xi_{2}y_{11},\ i\xi_{0}i\xi_{3}y_{11},\
i\xi_{1}i\xi_{2}y_{11},\ i\xi_{1}i\xi_{3}y_{11},\
i\xi_{2}i\xi_{3}y_{11},\ 0_{9\times1}^{T})^{T}\
,\end{eqnarray*}and\begin{eqnarray*}&&
y_{3}=\cfrac{i\xi_{1}(i\xi_{1}F(F_{1})+i\xi_{2}F(F_{2})+i\xi_{3}F(F_{3}))}{a^{2}b\tau}-\cfrac{F(F_{1})}{a\tau}\
,\\&&y_{7}=\cfrac{i\xi_{2}(i\xi_{1}F(F_{1})+i\xi_{2}F(F_{2})+i\xi_{3}F(F_{3}))}{a^{2}b\tau}-\cfrac{F(F_{2})}{a\tau}\
,
\\&&y_{11}=\cfrac{i\xi_{3}(i\xi_{1}F(F_{1})+i\xi_{2}F(F_{2})+i\xi_{3}F(F_{3}))}{a^{2}b\tau}-\cfrac{F(F_{3})}{a\tau}\
,\\&&y_{16}=\cfrac{i\xi_{1}F(F_{1})+i\xi_{2}F(F_{2})+i\xi_{3}F(F_{3})}{ab\tau}\
,\end{eqnarray*}where$$a=\cfrac{\mu[(i\xi_{1})^{2}+(i\xi_{2})^{2}+(i\xi_{3})^{2}]-i\xi_{0}}{\tau}\neq
                       0\ ,\
b=\cfrac{(i\xi_{1})^{2}+(i\xi_{2})^{2}+(i\xi_{3})^{2}}{a}\neq 0\ ,$$
when\ $(\xi_{1},\ \xi_{2},\ \xi_{3})\neq 0$\ .\
$A_{1}(\eta)=(\eta_{1} ,\ \eta_{2} ,\ \eta_{3} ,\ \eta_{4} ,\
\eta_{5} ,\ \eta_{6} ,\ \eta_{7} ,\ \eta_{8},\
                       \eta_{9})$\ , and\begin{eqnarray*}\eta_{j}&=&(i\xi_{2}y_{3},\ i\xi_{3}y_{3},\
y_{3},\ i\xi_{1}y_{7},\ i\xi_{2}y_{7},\ i\xi_{3}y_{7},\ y_{7},\
i\xi_{1}y_{11},\ i\xi_{2}y_{11},\ i\xi_{3}y_{11},\ y_{11},\\&&
i\xi_{0}y_{16},\ i\xi_{1}y_{16},\ i\xi_{2}y_{16},\ i\xi_{3}y_{16},\
y_{16},\ (i\xi_{0})^{2}y_{3},\ (i\xi_{1})^{2}y_{3},\
(i\xi_{2})^{2}y_{3},\ (i\xi_{3})^{2}y_{3},\
i\xi_{0}i\xi_{1}y_{3},\\&& i\xi_{0}i\xi_{2}y_{3},\
i\xi_{0}i\xi_{3}y_{3},\ i\xi_{1}i\xi_{2}y_{3},\
i\xi_{1}i\xi_{3}y_{3},\ i\xi_{2}i\xi_{3}y_{3},\
(i\xi_{0})^{2}y_{7},\ (i\xi_{1})^{2}y_{7},\ (i\xi_{2})^{2}y_{7},\
(i\xi_{3})^{2}y_{7},\\&& i\xi_{0}i\xi_{1}y_{7},\
i\xi_{0}i\xi_{2}y_{7},\ i\xi_{0}i\xi_{3}y_{7},\
i\xi_{1}i\xi_{2}y_{7},\ i\xi_{1}i\xi_{3}y_{7},\
i\xi_{2}i\xi_{3}y_{7},\ (i\xi_{0})^{2}y_{11},\
(i\xi_{1})^{2}y_{11},\ (i\xi_{2})^{2}y_{11},\\&&
(i\xi_{3})^{2}y_{11},\ i\xi_{0}i\xi_{1}y_{11},\
i\xi_{0}i\xi_{2}y_{11},\ i\xi_{0}i\xi_{3}y_{11},\
i\xi_{1}i\xi_{2}y_{11},\ i\xi_{1}i\xi_{3}y_{11},\
i\xi_{2}i\xi_{3}y_{11},\ e_{j}^{T})^{T}\ ,\end{eqnarray*}here\
$e_{j}$\ is the\ $jth$\ $9$\ dimensional unit coordinate vector,\ $1
\leq j \leq 9$\ , moreover when\ $j=1,\ 2,\ 3$\
,$$y_{3}=\cfrac{\xi_{1}^{2}}{a^{2}b\tau}+\cfrac{1}{a\tau}\ ,\
y_{7}=\cfrac{-i\xi_{1}i\xi_{2}}{a^{2}b\tau}\ ,\
y_{11}=\cfrac{-i\xi_{1}i\xi_{3}}{a^{2}b\tau}\ ,\
y_{16}=\cfrac{-i\xi_{1}}{ab\tau}\ ,$$when\ $j=4,\ 5,\ 6$\
,$$y_{7}=\cfrac{\xi_{2}^{2}}{a^{2}b\tau}+\cfrac{1}{a\tau}\ ,\
y_{3}=\cfrac{-i\xi_{1}i\xi_{2}}{a^{2}b\tau}\ ,\
y_{11}=\cfrac{-i\xi_{2}i\xi_{3}}{a^{2}b\tau}\ ,\
y_{16}=\cfrac{-i\xi_{2}}{ab\tau}\ ,$$when\ $j=7,\ 8,\ 9$\
,$$y_{11}=\cfrac{\xi_{3}^{2}}{a^{2}b\tau}+\cfrac{1}{a\tau}\ ,\
y_{3}=\cfrac{-i\xi_{1}i\xi_{3}}{a^{2}b\tau}\ ,\
y_{7}=\cfrac{-i\xi_{2}i\xi_{3}}{a^{2}b\tau}\ ,\
y_{16}=\cfrac{-i\xi_{3}}{ab\tau}\ .$$
            $Z_{1}=(Z_{1j})_{9\times1}$\ , moreover we assume$$F^{-1}(Y_{1}+A_{1}(\eta)Z_{1}I_{\{(\xi_{1},\ \xi_{2},\
\xi_{3})\neq 0\}})=F^{-1}(Y_{1}+A_{1}(\eta)Z_{1})\ ,$$ and\
$Z_{1}\in\ \Omega_{1}$\ , where
$$\Omega_{1}=\{Z_{1}|H[F^{-1}(Y_{1}+A_{1}(\eta)Z_{1})]=H[F^{-1}(Y_{1}+A_{1}(\eta)Z_{1})]I_{K_{1}^{\prime}}\in
C^{1}(K_{1}^{\prime}).\}\ .$$ (2)We can convert\ (2.2)\ into the
                       integral equations as follows,
\[\label {A-Problem}\begin{cases}Z_{11}=F[F^{-1}(e_{3}^{T}(Y_{1}+A_{1}(\eta)Z_{1}))F^{-1}(-\alpha_{1}^{T}(Y_{1}+A_{1}(\eta)Z_{1}))]=f_{1}(Z_{1})\ ,
\\Z_{12}=F[F^{-1}(e_{7}^{T}(Y_{1}+A_{1}(\eta)Z_{1}))F^{-1}(e_{1}^{T}(Y_{1}+A_{1}(\eta)Z_{1}))]=f_{2}(Z_{1})\ ,\\
Z_{13}=F[F^{-1}(e_{11}^{T}(Y_{1}+A_{1}(\eta)Z_{1}))F^{-1}(e_{2}^{T}(Y_{1}+A_{1}(\eta)Z_{1}))]=f_{3}(Z_{1})\
,\\
Z_{14}=F[F^{-1}(e_{3}^{T}(Y_{1}+A_{1}(\eta)Z_{1}))F^{-1}(e_{4}^{T}(Y_{1}+A_{1}(\eta)Z_{1}))]=f_{4}(Z_{1})\
,\\
Z_{15}=F[F^{-1}(e_{7}^{T}(Y_{1}+A_{1}(\eta)Z_{1}))F^{-1}(e_{5}^{T}(Y_{1}+A_{1}(\eta)Z_{1}))]=f_{5}(Z_{1})\
,\\
Z_{16}=F[F^{-1}(e_{11}^{T}(Y_{1}+A_{1}(\eta)Z_{1}))F^{-1}(e_{6}^{T}(Y_{1}+A_{1}(\eta)Z_{1}))]=f_{6}(Z_{1})\
,\\
Z_{17}=F[F^{-1}(e_{3}^{T}(Y_{1}+A_{1}(\eta)Z_{1}))F^{-1}(e_{8}^{T}(Y_{1}+A_{1}(\eta)Z_{1}))]=f_{7}(Z_{1})\
,\\
Z_{18}=F[F^{-1}(e_{7}^{T}(Y_{1}+A_{1}(\eta)Z_{1}))F^{-1}(e_{9}^{T}(Y_{1}+A_{1}(\eta)Z_{1}))]=f_{8}(Z_{1})\
,\\
Z_{19}=F[F^{-1}(e_{11}^{T}(Y_{1}+A_{1}(\eta)Z_{1}))F^{-1}(e_{10}^{T}(Y_{1}+A_{1}(\eta)Z_{1}))]=f_{9}(Z_{1})\
.\end{cases}
\]This is also the question to find the fixed-points
of\ $f(Z_{1})$\ , where\ $f(Z_{1})=(f_{j}(Z_{1}))_{9\times1}$\ . We
can use the theories about the Poisson's equation, the heat-conduct
equation, the Schauder fixed-point theorem to prove that the
fixed-point is exist in the region\ $K_{1}^{\prime}$\ , and the
smoothing solution for the Navier-Stokes equations is globally exist
in the region\ $K_{1}^{\prime}$\ .\end{theorem}

\section{Proof}\setcounter{equation}{0}
{\it Proof of theorem2-1}. (1)Necessity. Because the solution of\
$AX=\beta$\ is exist, we assume the particular solution is\ $X_{0}$\
, then\ $\forall\ (x_{1},\ x_{2},\ x_{3},\ t)^{T} \in R^{4} $\ , we
know that$$AX_{0}(x_{1},\ x_{2},\ x_{3},\ t)=\beta(x_{1},\ x_{2},\
x_{3},\ t)$$satisfies,
$$X_{0}(x_{1},\ x_{2},\ x_{3},\ t)\
  \mbox{and} \
   \beta(x_{1},\ x_{2},\ x_{3},\ t)$$ are values of\ $X_{0}$\ and\
$\beta$\ , when\
   $ (x_{1},\ x_{2},\ x_{3},\ t)^{T}$\ is given. Hence
   $$rank(A,\ \beta(x_{1},\ x_{2},\ x_{3},\ t))=rank(A),\ \forall\
(x_{1},\ x_{2},\ x_{3},\ t)^{T} \in R^{4}\ . $$ Sufficiency.
Because\
 $\forall\
 (x_{1},\ x_{2},\ x_{3},\ t)^{T} \in R^{4}$\ , we can get
 $$rank(A,\ \beta(x_{1},\ x_{2},\ x_{3},\ t))=rank(A)\
 ,\ $$then the solution of$$AX=\beta(x_{1},\ x_{2},\ x_{3},\ t)$$is exist. We
assume the particular solution is\ $X_{0}(x_{1},\ x_{2},\ x_{3},\
t)$\ , then\ $\forall\ (x_{1},\ x_{2},\ x_{3},\ t)^{T}
  \in R^{4}$\ , \\we let
$$X=\left(X_{1}(x_{1},\ x_{2},\ x_{3},\ t),\ X_{2}(x_{1},\ x_{2},\
x_{3},\ t),\ \cdots,\ X_{s}(x_{1},\ x_{2},\ x_{3},\ t)\right)^{T}
                  =X_{0}(x_{1},\ x_{2},\ x_{3},\ t)\
                  ,\ $$and we can get\
$X$\ satisfies\ $AX=\beta$\ .\\ (2)When the solution of\ $AX=\beta$\
is exist, we assume\ $X_{0}$\ is a particular solution,\ $X$\ is the
arbitrary solution of this equations, then\ $\forall\ (x_{1},\
x_{2},\ x_{3},\ t)^{T} \in R^{4}$\ , we can get
$$A(X(x_{1},\ x_{2},\ x_{3},\ t)-X_{0}(x_{1},\ x_{2},\ x_{3},\ t))=0\ ,\
$$ according to the linear equations
theory, there exists a unique team of real numbers
$$C_{1},\ C_{2},\ \cdots,\ C_{s-r}$$such that
$$X(x_{1},\ x_{2},\ x_{3},\ t)-X_{0}(x_{1},\ x_{2},\ x_{3},\ t)=C_{1}\eta_{1}+C_{2}\eta_{2}+\cdots+C_{s-r}\eta_{s-r}\
,\ $$we let
$$Z_{1}(x_{1},\ x_{2},\ x_{3},\ t)=\left(C_{1},\ C_{2},\ \cdots,\ C_{s-r}\right)^{T}\
,\ $$then
$$X(x_{1},\ x_{2},\ x_{3},\ t)=X_{0}(x_{1},\ x_{2},\ x_{3},\ t)+A(\eta)Z_{1}(x_{1},\ x_{2},\ x_{3},\ t)\
.$$ Because\
 $(x_{1},\ x_{2},\ x_{3},\ t)^{T}$\ is arbitrary, we can learn that$$X=X_{0}+A(\eta)Z_{1}\ .$$
On the other hand$$A[X_{0}+A(\eta)Z_{1}]=\beta+0Z_{1}=\beta\ .$$
         This is to
say $$\forall Z_{1}\ ,\ X_{0}+A(\eta)Z_{1}$$is the solution of the
equation\ $AX=\beta$\ . Hence
$$X_{0}+A(\eta)Z_{1}$$is the general solution of the
equation\ $AX=\beta$\ .
\\{\it Proof of theorem2-2}.\ We can rewrite the Navier-Stokes
equations as following.
$$AX=\beta\ ,$$ where$$A=\left(
\begin{array}{ccccccccccccccc}
1,\  & 0,\ & 0,\  & 0,\  & 0_{4},\  & 1,\  & 0_{4},\ & 1,\  & 0,\  & 0,\  & 0,\  & 0,\  & 0,\  & 0,\  & 0_{39} \\
0,\  & 1,\  & 0,\  & 0,\  & 0_{12},\  & \tau,\  & 0_{4},\  & -\mu,\ &-\mu,\   & -\mu,\  & 0_{26},\ & 1,\  & 1,\  & 1,\  & 0_{6} \\
0,\ & 0,\  & 1,\  & 0,\  & 0_{13},\  & \tau,\  & 0_{13},\  & -\mu,\  & -\mu,\  & -\mu,\  & 0_{19},\  & 1,\  & 1,\  & 1,\  & 0_{3} \\
0,\  & 0,\  & 0,\  & 1,\  & 0_{14},\  & \tau,\  & 0_{22},\  & -\mu,\  & -\mu,\  & -\mu,\  & 0_{11},\  & 0,\  & 1,\ & 1,\  & 1 \\
\end{array}
\right)_{4\times 59}
$$
$0_{n}$\ is the row vector which is made up of\ $n$\ zeroes,\ $X$\
includes\ $u_{1},\ u_{2},\ u_{3},\ p$\ and all their first order
partial derivative, and all the second order partial derivative of\
$u_{1},\ u_{2},\ u_{3}$\ , and all the products which they are in
the Navier-Stokes equations,
\begin{eqnarray*}&X=&(\cfrac{\partial u_{1}}{\partial x_{1}}\
,\ \cfrac{\partial u_{1}}{\partial t}\ ,\ \cfrac{\partial u_{2}
}{\partial t}\ ,\ \cfrac{\partial u_{3}}{\partial t}\ ,\
\cfrac{\partial u_{1}}{\partial x_{2}}\ ,\ \cfrac{\partial
u_{1}}{\partial x_{3}}\ ,\ u_{1}\ ,\ \cfrac{\partial u_{2}
}{\partial x_{1}}\ ,\ \cfrac{\partial u_{2}}{\partial x_{2}}\ ,\
\cfrac{\partial u_{2}}{\partial x_{3}}\ ,\ u_{2}\ ,\ \cfrac{\partial
u_{3} }{\partial x_{1}}\ ,\ \cfrac{\partial u_{3}}{\partial x_{2}}\
,\ \cfrac{\partial u_{3}}{\partial x_{3}}\ ,\ u_{3}\ ,\
\\&
&\cfrac{\partial p }{\partial t}\ ,\ \cfrac{\partial p }{\partial
x_{1}}\ ,\ \cfrac{\partial p}{\partial x_{2}}\ ,\ \cfrac{\partial
p}{\partial x_{3}}\ ,\ p\ ,\ \cfrac{\partial^{2} u_{1}}{\partial
 t^{2}}\
 ,\ \cfrac{\partial^{2} u_{1}}{\partial
 x_{1}^{2}}\
 ,\ \cfrac{\partial^{2} u_{1}}{\partial
 x_{2}^{2}}\
 ,\ \cfrac{\partial^{2} u_{1}}{\partial
  x_{3}^{2}}\
 ,\ \cfrac{\partial^{2} u_{1}}{\partial
 t\partial x_{1} }\
 ,\ \cfrac{\partial^{2} u_{1}}{\partial
 t\partial x_{2}}\
 ,\ \cfrac{\partial^{2} u_{1}}{\partial
 t\partial x_{3}}\
 ,\end{eqnarray*}\begin{eqnarray*}&&\cfrac{\partial^{2} u_{1}}{\partial
  x_{1}\partial x_{2}}\
 ,\ \cfrac{\partial^{2} u_{1}}{\partial
  x_{1}\partial x_{3}}\
 ,\ \cfrac{\partial^{2} u_{1}}{\partial
  x_{2}\partial x_{3}}\
 ,\ \cfrac{\partial^{2} u_{2}}{\partial
 t^{2}}\
 ,\ \cfrac{\partial^{2} u_{2}}{\partial
 x_{1}^{2}}\
 ,\ \cfrac{\partial^{2} u_{2}}{\partial
 x_{2}^{2}}\
 ,\ \cfrac{\partial^{2} u_{2}}{\partial
  x_{3}^{2}}\
 ,\ \cfrac{\partial^{2} u_{2}}{\partial
 t\partial x_{1} }\
 ,\ \cfrac{\partial^{2} u_{2}}{\partial
 t\partial x_{2}}\
 ,\ \cfrac{\partial^{2} u_{2}}{\partial
 t\partial x_{3}}\
 ,\ \\&&\cfrac{\partial^{2} u_{2}}{\partial
  x_{1}\partial x_{2}}\
 ,\ \cfrac{\partial^{2} u_{2}}{\partial
  x_{1}\partial x_{3}}\
 ,\ \cfrac{\partial^{2} u_{2}}{\partial
  x_{2}\partial x_{3}}\
 ,\ \cfrac{\partial^{2} u_{3}}{\partial
 t^{2}}\
 ,\ \cfrac{\partial^{2} u_{3}}{\partial
 x_{1}^{2}}\
 ,\ \cfrac{\partial^{2} u_{3}}{\partial
 x_{2}^{2}}\
 ,\ \cfrac{\partial^{2} u_{3}}{\partial
  x_{3}^{2}}\
 ,\ \cfrac{\partial^{2} u_{3}}{\partial
 t\partial x_{1} }\
 ,\ \cfrac{\partial^{2} u_{3}}{\partial
 t\partial x_{2}}\
 ,\ \cfrac{\partial^{2} u_{3}}{\partial
 t\partial x_{3}}\
 ,\ \\& &\cfrac{\partial^{2} u_{3}}{\partial
  x_{1}\partial x_{2}}\
 ,\ \cfrac{\partial^{2} u_{3}}{\partial
  x_{1}\partial x_{3}}\
 ,\ \cfrac{\partial^{2} u_{3}}{\partial
  x_{2}\partial x_{3}}\
 ,\ u_{1} \cfrac{\partial u_{1}}{\partial
x_{1}}\ ,\ u_{2} \cfrac{\partial u_{1}}{\partial x_{2}}\ ,\
u_{3}\cfrac{\partial u_{1}}{\partial x_{3}}\ ,\ u_{1}
\cfrac{\partial
 u_{2}}{\partial x_{1}}\
 ,\ u_{2} \cfrac{\partial u_{2}}{\partial x_{2}}\
 ,\
u_{3}\cfrac{\partial u_{2}}{\partial x_{3}}\ ,\ \\&&u_{1}
\cfrac{\partial u_{3} }{\partial x_{1}}\ ,\ u_{2} \cfrac{\partial
u_{3}}{\partial x_{2}}\ ,\  u_{3}\cfrac{\partial u_{3}}{\partial
x_{3}}\ )^{T}\ ,\
\end{eqnarray*}$$\beta=(0,\ F_{1},\ F_{2},\ F_{3})\
.$$ We can see that\ $A=(E_{4},A_{1})$\ ,\ $A_{1}$\ is a\
$4\times55$\ matrix. We assume\ $A_{1}=(\alpha_{1},\
                                    \alpha_{2},\ \alpha_{3},\
                                    \alpha_{4})^{T}$\ , where\begin{eqnarray*}&&\alpha_{1}=(0_{4},\ 1,\ 0_{4},\  1,\  0,\  0,\  0,\  0,\ 0,\  0,\  0_{39})^{T},
                                    \\&& \alpha_{2}=(0_{12},\  \tau,\  0_{4},\  -\mu,\ -\mu,\  -\mu,\
0_{26},\  1,\  1,\  1,\  0_{6} )^{T},\\&& \alpha_{3}=(0_{13},\
\tau,\   0_{13},\   -\mu,\   -\mu,\   -\mu,\   0_{19},\   1,\   1,\
1,\  0_{3} )^{T},\\&&  \alpha_{4}=(0_{14},\   \tau,\   0_{22},\
-\mu,\   -\mu,\   -\mu,\   0_{11},\   0,\   1,\ 1,\   1)^{T}\
.\end{eqnarray*} We let\ $X_{0}=(0,\ F_1,\ F_2,\ F_3,\ 0_{55})^{T}$\
,\ $A(\eta)=\left(
                                      \begin{array}{c}
                                        -A_{1} \\
                                        E_{55} \\
                                      \end{array}
                                    \right)_{59\times55}\ ,$
\ then\ $AX_{0}=\beta\ ,$\ $AA(\eta)=0\ ,$\ and the columns of\
$A(\eta)$\ are linear independent, so according to (2) in the
theorem 2.2, we can get the general solution of\ $AX=\beta$\ as
follows,
$$X=X_{0}+A(\eta)Z\ .$$We write the\ $59$\ components of\ $X$\ in
details.
\begin{eqnarray*}&&\cfrac{\partial u_{1}}{\partial x_{1}}=-\alpha_{1}^{T}Z\
,\ \cfrac{\partial u_{1}}{\partial t}=F_{1}-\alpha_{2}^{T}Z\ ,\
\cfrac{\partial u_{2} }{\partial t}=F_{2}-\alpha_{3}^{T}Z\ ,\
\cfrac{\partial u_{3}}{\partial t}=F_{3}-\alpha_{4}^{T}Z \ ,\
\cfrac{\partial u_{1}}{\partial x_{2}}=e_{1}^{T}Z\ ,\\&&
\cfrac{\partial u_{1}}{\partial x_{3}}=e_{2}^{T}Z\ ,\
u_{1}=e_{3}^{T}Z\ ,\ \cfrac{\partial u_{2} }{\partial
x_{1}}=e_{4}^{T}Z\ ,\ \cfrac{\partial u_{2}}{\partial
x_{2}}=e_{5}^{T}Z\ ,\ \cfrac{\partial u_{2}}{\partial
x_{3}}=e_{6}^{T}Z\ ,\ u_{2}=e_{7}^{T}Z\ ,\\&& \cfrac{\partial u_{3}
}{\partial x_{1}}=e_{8}^{T}Z\ ,\ \cfrac{\partial u_{3}}{\partial
x_{2}}=e_{9}^{T}Z\ ,\ \cfrac{\partial u_{3}}{\partial
x_{3}}=e_{10}^{T}Z\ ,\ u_{3}=e_{11}^{T}Z\ ,\ \cfrac{\partial p
}{\partial t}=e_{12}^{T}Z\ ,\ \cfrac{\partial p }{\partial
x_{1}}=e_{13}^{T}Z\ ,\\&& \cfrac{\partial p}{\partial
x_{2}}=e_{14}^{T}Z\ ,\ \cfrac{\partial p}{\partial
x_{3}}=e_{15}^{T}Z\ ,\ p=e_{16}^{T}Z\ ,\ \cfrac{\partial^{2}
u_{1}}{\partial
 t^{2}}=e_{17}^{T}Z\
 ,\ \cfrac{\partial^{2} u_{1}}{\partial
 x_{1}^{2}}=e_{18}^{T}Z\
 ,\ \cfrac{\partial^{2} u_{1}}{\partial
 x_{2}^{2}}=e_{19}^{T}Z\
 ,\\&& \cfrac{\partial^{2} u_{1}}{\partial
  x_{3}^{2}}=e_{20}^{T}Z\
 ,\ \cfrac{\partial^{2} u_{1}}{\partial
 t\partial x_{1} }=e_{21}^{T}Z\
 ,\ \cfrac{\partial^{2} u_{1}}{\partial
 t\partial x_{2}}=e_{22}^{T}Z\
 ,\ \cfrac{\partial^{2} u_{1}}{\partial
 t\partial x_{3}}=e_{23}^{T}Z\
 ,\ \cfrac{\partial^{2} u_{1}}{\partial
  x_{1}\partial x_{2}}=e_{24}^{T}Z\
 ,\\&& \cfrac{\partial^{2} u_{1}}{\partial
  x_{1}\partial x_{3}}=e_{25}^{T}Z\
 ,\ \cfrac{\partial^{2} u_{1}}{\partial
  x_{2}\partial x_{3}}=e_{26}^{T}Z\
 ,\ \cfrac{\partial^{2} u_{2}}{\partial
 t^{2}}=e_{27}^{T}Z\
 ,\ \cfrac{\partial^{2} u_{2}}{\partial
 x_{1}^{2}}=e_{28}^{T}Z\
 ,\ \cfrac{\partial^{2} u_{2}}{\partial
 x_{2}^{2}}=e_{29}^{T}Z\
 ,\\
 && \cfrac{\partial^{2} u_{2}}{\partial
  x_{3}^{2}}=e_{30}^{T}Z\
 ,\ \cfrac{\partial^{2} u_{2}}{\partial
 t\partial x_{1} }=e_{31}^{T}Z\
 ,\ \cfrac{\partial^{2} u_{2}}{\partial
 t\partial x_{2}}=e_{32}^{T}Z\
 ,\ \cfrac{\partial^{2} u_{2}}{\partial
 t\partial x_{3}}=e_{33}^{T}Z\
 ,\ \cfrac{\partial^{2} u_{2}}{\partial
  x_{1}\partial x_{2}}=e_{34}^{T}Z\
 , \\&&\cfrac{\partial^{2} u_{2}}{\partial
  x_{1}\partial x_{3}}=e_{35}^{T}Z\
 ,\ \cfrac{\partial^{2} u_{2}}{\partial
  x_{2}\partial x_{3}}=e_{36}^{T}Z\
 ,\ \cfrac{\partial^{2} u_{3}}{\partial
 t^{2}}=e_{37}^{T}Z\
 ,\ \cfrac{\partial^{2} u_{3}}{\partial
 x_{1}^{2}}=e_{38}^{T}Z\
 , \cfrac{\partial^{2} u_{3}}{\partial
 x_{2}^{2}}=e_{39}^{T}Z\
 ,\end{eqnarray*}\begin{eqnarray*}&& \cfrac{\partial^{2} u_{3}}{\partial
  x_{3}^{2}}=e_{40}^{T}Z\
 ,\ \cfrac{\partial^{2} u_{3}}{\partial
 t\partial x_{1} }=e_{41}^{T}Z\
 ,\ \cfrac{\partial^{2} u_{3}}{\partial
 t\partial x_{2}}=e_{42}^{T}Z\
 ,\ \cfrac{\partial^{2} u_{3}}{\partial
 t\partial x_{3}}=e_{43}^{T}Z\
 ,\ \cfrac{\partial^{2} u_{3}}{\partial
  x_{1}\partial x_{2}}=e_{44}^{T}Z\
 ,\\&& \cfrac{\partial^{2} u_{3}}{\partial
  x_{1}\partial x_{3}}=e_{45}^{T}Z\
 ,\ \cfrac{\partial^{2} u_{3}}{\partial
  x_{2}\partial x_{3}}=e_{46}^{T}Z\
 ,\ u_{1} \cfrac{\partial u_{1}}{\partial
x_{1}}=e_{47}^{T}Z\ ,\ u_{2} \cfrac{\partial u_{1}}{\partial
x_{2}}=e_{48}^{T}Z\ ,\ u_{3}\cfrac{\partial u_{1}}{\partial
x_{3}}=e_{49}^{T}Z\ ,\\&& u_{1} \cfrac{\partial
 u_{2}}{\partial x_{1}}=e_{50}^{T}Z\
 ,\ u_{2} \cfrac{\partial u_{2}}{\partial x_{2}}=e_{51}^{T}Z\
 ,\
u_{3}\cfrac{\partial u_{2}}{\partial x_{3}}=e_{52}^{T}Z\ ,\ u_{1}
\cfrac{\partial u_{3} }{\partial x_{1}}=e_{53}^{T}Z\ ,\ u_{2}
\cfrac{\partial u_{3}}{\partial x_{2}}=e_{54}^{T}Z\ ,\\&&
u_{3}\cfrac{\partial u_{3}}{\partial x_{3}}=e_{55}^{T}Z\ ,
\end{eqnarray*}here\ $e_{i}$\ is the\ $ith$\ $55$\
dimensional unit coordinate vector,\ $1 \leq i \leq 55\ .$\ There
are\ $9$\ quadratic items, hence we can get\ $9$\ quadratic
constraints which\ $Z$\ need to satisfy just as following:
\[\begin{cases} e_{47}^{T}Z=(e_{3}^{T}Z)(-\alpha_{1}^{T}Z),\
e_{48}^{T}Z=(e_{7}^{T}Z)(e_{1}^{T}Z),\
e_{49}^{T}Z=(e_{11}^{T}Z)(e_{2}^{T}Z),\\
e_{50}^{T}Z=(e_{3}^{T}Z)(e_{4}^{T}Z),\
e_{51}^{T}Z=(e_{7}^{T}Z)(e_{5}^{T}Z),\
e_{52}^{T}Z=(e_{11}^{T}Z)(e_{6}^{T}Z),\\
e_{53}^{T}Z=(e_{3}^{T}Z)(e_{8}^{T}Z),\
e_{54}^{T}Z=(e_{7}^{T}Z)(e_{9}^{T}Z),\
e_{55}^{T}Z=(e_{11}^{T}Z)(e_{10}^{T}Z).\
\end{cases}\]Because\ $u_{1},\ u_{2},\ u_{3}\in C^{2}\ ,\ p\in
C^{1}$\ , we can get\ $46$\ differential constraints which\ $Z$\
need to satisfy just as following:
\begin{eqnarray*}&&\cfrac{\partial
(e_{3}^{T}Z)}{\partial x_{1}}=-\alpha_{1}^{T}Z\ ,\ \cfrac{\partial
(e_{3}^{T}Z)}{\partial t}=F_{1}-\alpha_{2}^{T}Z\ ,\ \cfrac{\partial
(e_{7}^{T}Z) }{\partial t}=F_{2}-\alpha_{3}^{T}Z\ ,\ \cfrac{\partial
(e_{11}^{T}Z)}{\partial t}=F_{3}-\alpha_{4}^{T}Z \ ,\\&&
\cfrac{\partial (e_{3}^{T}Z)}{\partial x_{2}}=e_{1}^{T}Z\ ,\
\cfrac{\partial (e_{3}^{T}Z)}{\partial x_{3}}=e_{2}^{T}Z\ ,\
\cfrac{\partial (e_{7}^{T}Z) }{\partial x_{1}}=e_{4}^{T}Z\ ,\
\cfrac{\partial (e_{7}^{T}Z)}{\partial x_{2}}=e_{5}^{T}Z\ ,
\cfrac{\partial (e_{7}^{T}Z)}{\partial x_{3}}=e_{6}^{T}Z\ ,\
\\&& \cfrac{\partial (e_{11}^{T}Z) }{\partial
x_{1}}=e_{8}^{T}Z\ ,\ \cfrac{\partial (e_{11}^{T}Z)}{\partial
x_{2}}=e_{9}^{T}Z\ ,\ \cfrac{\partial (e_{11}^{T}Z)}{\partial
x_{3}}=e_{10}^{T}Z\ ,\ \cfrac{\partial (e_{16}^{T}Z) }{\partial
t}=e_{12}^{T}Z\ ,\ \cfrac{\partial (e_{16}^{T}Z) }{\partial
x_{1}}=e_{13}^{T}Z\ ,\\&& \cfrac{\partial (e_{16}^{T}Z)}{\partial
x_{2}}=e_{14}^{T}Z\ ,\ \cfrac{\partial (e_{16}^{T}Z)}{\partial
x_{3}}=e_{15}^{T}Z\ ,\ \cfrac{\partial^{2} (e_{3}^{T}Z)}{\partial
 t^{2}}=e_{17}^{T}Z\
 ,\ \cfrac{\partial^{2} (e_{3}^{T}Z)}{\partial
 x_{1}^{2}}=e_{18}^{T}Z\
 ,\ \cfrac{\partial^{2} (e_{3}^{T}Z)}{\partial
 x_{2}^{2}}=e_{19}^{T}Z\
 ,\\&& \cfrac{\partial^{2} (e_{3}^{T}Z)}{\partial
  x_{3}^{2}}=e_{20}^{T}Z\
 ,\ \cfrac{\partial^{2} (e_{3}^{T}Z)}{\partial
 t\partial x_{1} }=e_{21}^{T}Z\
 ,\ \cfrac{\partial^{2} (e_{3}^{T}Z)}{\partial
 t\partial x_{2}}=e_{22}^{T}Z\
 ,\ \cfrac{\partial^{2} (e_{3}^{T}Z)}{\partial
 t\partial x_{3}}=e_{23}^{T}Z\
 ,\ \cfrac{\partial^{2} (e_{3}^{T}Z)}{\partial
  x_{1}\partial x_{2}}=e_{24}^{T}Z\
 ,\\&& \cfrac{\partial^{2} (e_{3}^{T}Z)}{\partial
  x_{1}\partial x_{3}}=e_{25}^{T}Z\
 ,\ \cfrac{\partial^{2} (e_{3}^{T}Z)}{\partial
  x_{2}\partial x_{3}}=e_{26}^{T}Z\
 ,\ \cfrac{\partial^{2} (e_{7}^{T}Z)}{\partial
 t^{2}}=e_{27}^{T}Z\
 ,\ \cfrac{\partial^{2} (e_{7}^{T}Z)}{\partial
 x_{1}^{2}}=e_{28}^{T}Z\
 ,\ \cfrac{\partial^{2} (e_{7}^{T}Z)}{\partial
 x_{2}^{2}}=e_{29}^{T}Z\
 ,\\&& \cfrac{\partial^{2} (e_{7}^{T}Z)}{\partial
  x_{3}^{2}}=e_{30}^{T}Z\
 ,\ \cfrac{\partial^{2} (e_{7}^{T}Z)}{\partial
 t\partial x_{1} }=e_{31}^{T}Z\
 ,\ \cfrac{\partial^{2} (e_{7}^{T}Z)}{\partial
 t\partial x_{2}}=e_{32}^{T}Z\
 ,\ \cfrac{\partial^{2} (e_{7}^{T}Z)}{\partial
 t\partial x_{3}}=e_{33}^{T}Z\
 ,\ \cfrac{\partial^{2} (e_{7}^{T}Z)}{\partial
  x_{1}\partial x_{2}}=e_{34}^{T}Z\
 , \\&&\cfrac{\partial^{2} (e_{7}^{T}Z)}{\partial
  x_{1}\partial x_{3}}=e_{35}^{T}Z\
 ,\ \cfrac{\partial^{2} (e_{7}^{T}Z)}{\partial
  x_{2}\partial x_{3}}=e_{36}^{T}Z\
 ,\ \cfrac{\partial^{2} (e_{11}^{T}Z)}{\partial
 t^{2}}=e_{37}^{T}Z\
 ,\ \cfrac{\partial^{2} (e_{11}^{T}Z)}{\partial
 x_{1}^{2}}=e_{38}^{T}Z\
 , \cfrac{\partial^{2} (e_{11}^{T}Z)}{\partial
 x_{2}^{2}}=e_{39}^{T}Z\
 ,\\&& \cfrac{\partial^{2} (e_{11}^{T}Z)}{\partial
  x_{3}^{2}}=e_{40}^{T}Z\
 ,\ \cfrac{\partial^{2} (e_{11}^{T}Z)}{\partial
 t\partial x_{1} }=e_{41}^{T}Z\
 ,\ \cfrac{\partial^{2} (e_{11}^{T}Z)}{\partial
 t\partial x_{2}}=e_{42}^{T}Z\
 ,\ \cfrac{\partial^{2} (e_{11}^{T}Z)}{\partial
 t\partial x_{3}}=e_{43}^{T}Z\
 ,\ \cfrac{\partial^{2} (e_{11}^{T}Z)}{\partial
  x_{1}\partial x_{2}}=e_{44}^{T}Z\
 ,\\&& \cfrac{\partial^{2} (e_{11}^{T}Z)}{\partial
  x_{1}\partial x_{3}}=e_{45}^{T}Z\
 ,\ \cfrac{\partial^{2} (e_{11}^{T}Z)}{\partial
  x_{2}\partial x_{3}}=e_{46}^{T}Z\
 . \end{eqnarray*}Because the second order partial derivative can be
taken as the partial derivative of the first order partial
derivative, and\ $u_{1},\ u_{2},\ u_{3}\in C^{2}$\ , we can learn
that their second order mixed partial derivatives are equal, the
above differential constraints are equivalent to\ $64$\ first order
differential constraints as
follows:\begin{eqnarray*}&&\cfrac{\partial (e_{3}^{T}Z)}{\partial
x_{1}}=-\alpha_{1}^{T}Z\ ,\ \cfrac{\partial (e_{3}^{T}Z)}{\partial
t}=F_{1}-\alpha_{2}^{T}Z\ ,\ \cfrac{\partial (e_{7}^{T}Z) }{\partial
t}=F_{2}-\alpha_{3}^{T}Z\ ,\ \cfrac{\partial (e_{11}^{T}Z)}{\partial
t}=F_{3}-\alpha_{4}^{T}Z \ ,\\&& \cfrac{\partial
(e_{3}^{T}Z)}{\partial x_{2}}=e_{1}^{T}Z\ ,\ \cfrac{\partial
(e_{3}^{T}Z)}{\partial x_{3}}=e_{2}^{T}Z\ ,\ \cfrac{\partial
(e_{7}^{T}Z) }{\partial x_{1}}=e_{4}^{T}Z\ ,\ \cfrac{\partial
(e_{7}^{T}Z)}{\partial x_{2}}=e_{5}^{T}Z\ , \cfrac{\partial
(e_{7}^{T}Z)}{\partial x_{3}}=e_{6}^{T}Z\ ,\
\\&& \cfrac{\partial (e_{11}^{T}Z) }{\partial
x_{1}}=e_{8}^{T}Z\ ,\ \cfrac{\partial (e_{11}^{T}Z)}{\partial
x_{2}}=e_{9}^{T}Z\ ,\ \cfrac{\partial (e_{11}^{T}Z)}{\partial
x_{3}}=e_{10}^{T}Z\ ,\ \cfrac{\partial (e_{16}^{T}Z) }{\partial
t}=e_{12}^{T}Z\ ,\ \cfrac{\partial (e_{16}^{T}Z) }{\partial
x_{1}}=e_{13}^{T}Z\ ,\\&& \cfrac{\partial (e_{16}^{T}Z)}{\partial
x_{2}}=e_{14}^{T}Z\ ,\ \cfrac{\partial (e_{16}^{T}Z)}{\partial
x_{3}}=e_{15}^{T}Z\ ,\ \cfrac{\partial
(F_{1}-\alpha_{2}^{T}Z)}{\partial
 t}=e_{17}^{T}Z\
 ,\ \cfrac{\partial (-\alpha_{1}^{T}Z)}{\partial
 x_{1}}=e_{18}^{T}Z\
 ,\\&& \cfrac{\partial (e_{1}^{T}Z)}{\partial
 x_{2}}=e_{19}^{T}Z\
 ,\ \cfrac{\partial (e_{2}^{T}Z)}{\partial
  x_{3}}=e_{20}^{T}Z\
 ,\ \cfrac{\partial (F_{1}-\alpha_{2}^{T}Z)}{\partial
  x_{1} }=\cfrac{\partial (-\alpha_{1}^{T}Z)}{\partial
  t }=e_{21}^{T}Z\
 ,\\&& \cfrac{\partial (F_{1}-\alpha_{2}^{T}Z)}{\partial
  x_{2} }=\cfrac{\partial (e_{1}^{T}Z)}{\partial
  t }=e_{22}^{T}Z\
 ,\ \cfrac{\partial (F_{1}-\alpha_{2}^{T}Z)}{\partial
  x_{3} }=\cfrac{\partial (e_{2}^{T}Z)}{\partial
  t }=e_{23}^{T}Z\
 ,\\&& \cfrac{\partial (-\alpha_{1}^{T}Z)}{\partial
  x_{2} }=\cfrac{\partial (e_{1}^{T}Z)}{\partial
  x_{1}}=e_{24}^{T}Z\
 ,\ \cfrac{\partial (-\alpha_{1}^{T}Z)}{\partial
  x_{3} }= \cfrac{\partial (e_{2}^{T}Z)}{\partial
  x_{1}}=e_{25}^{T}Z\
 ,\ \cfrac{\partial (e_{1}^{T}Z)}{\partial
   x_{3}}=\cfrac{\partial (e_{2}^{T}Z)}{\partial
  x_{2}}=e_{26}^{T}Z\
 ,\\&& \cfrac{\partial (F_{2}-\alpha_{3}^{T}Z)}{\partial
 t}=e_{27}^{T}Z\
 ,\ \cfrac{\partial (e_{4}^{T}Z)}{\partial
 x_{1}}=e_{28}^{T}Z\
 ,\ \cfrac{\partial (e_{5}^{T}Z)}{\partial
 x_{2}}=e_{29}^{T}Z\
 ,\ \cfrac{\partial (e_{6}^{T}Z)}{\partial
  x_{3}}=e_{30}^{T}Z\
 ,\\&& \cfrac{\partial (F_{2}-\alpha_{3}^{T}Z)}{\partial
 x_{1}}=\cfrac{\partial (e_{4}^{T}Z)}{\partial
 t }=e_{31}^{T}Z\
 ,\ \cfrac{\partial (F_{2}-\alpha_{3}^{T}Z)}{\partial
 x_{2}}=\cfrac{\partial (e_{5}^{T}Z)}{\partial
 t}=e_{32}^{T}Z\
 ,\\&& \cfrac{\partial (F_{2}-\alpha_{3}^{T}Z)}{\partial
 x_{3}}=\cfrac{\partial (e_{6}^{T}Z)}{\partial
 t}=e_{33}^{T}Z\
 ,\ \cfrac{\partial (e_{4}^{T}Z)}{\partial
 x_{2} }=\cfrac{\partial (e_{5}^{T}Z)}{\partial
  x_{1}}=e_{34}^{T}Z\
 ,\ \cfrac{\partial (e_{4}^{T}Z)}{\partial
 x_{3} }=\cfrac{\partial (e_{6}^{T}Z)}{\partial
  x_{1}}=e_{35}^{T}Z\
 ,\\&& \cfrac{\partial (e_{5}^{T}Z)}{\partial
  x_{3}}=\cfrac{\partial (e_{6}^{T}Z)}{\partial
  x_{2}}=e_{36}^{T}Z\
 ,\ \cfrac{\partial (F_{3}-\alpha_{4}^{T}Z)}{\partial
 t}=e_{37}^{T}Z\
 ,\ \cfrac{\partial (e_{8}^{T}Z)}{\partial
 x_{1}}=e_{38}^{T}Z\
 , \cfrac{\partial (e_{9}^{T}Z)}{\partial
 x_{2}}=e_{39}^{T}Z\
 ,\\&& \cfrac{\partial (e_{10}^{T}Z)}{\partial
  x_{3}}=e_{40}^{T}Z\
 ,\ \cfrac{\partial (F_{3}-\alpha_{4}^{T}Z)}{\partial
 x_{1}}=\cfrac{\partial (e_{8}^{T}Z)}{\partial
 t }=e_{41}^{T}Z\
 ,\ \cfrac{\partial (F_{3}-\alpha_{4}^{T}Z)}{\partial
 x_{2}}=\cfrac{\partial (e_{9}^{T}Z)}{\partial
 t}=e_{42}^{T}Z\
 ,\\&& \cfrac{\partial (F_{3}-\alpha_{4}^{T}Z)}{\partial
 x_{3}}=\cfrac{\partial (e_{10}^{T}Z)}{\partial
 t}=e_{43}^{T}Z\
 ,\ \cfrac{\partial (e_{8}^{T}Z)}{\partial
 x_{2} }=\cfrac{\partial (e_{9}^{T}Z)}{\partial
  x_{1}}=e_{44}^{T}Z\
 ,\ \cfrac{\partial (e_{8}^{T}Z)}{\partial
 x_{3} }=\cfrac{\partial (e_{10}^{T}Z)}{\partial
  x_{1}}=e_{45}^{T}Z\
 ,\\&& \cfrac{\partial (e_{9}^{T}Z)}{\partial
 x_{3} }=\cfrac{\partial (e_{10}^{T}Z)}{\partial
  x_{2}}=e_{46}^{T}Z\
 .
\end{eqnarray*}We
write them into the equations, we assume
$$U=(\cfrac{\partial u_{1}}{\partial
x_{1}}\ ,\ \cfrac{\partial u_{1}}{\partial t}\ ,\ \cfrac{\partial
u_{2} }{\partial t}\ ,\ \cfrac{\partial u_{3}}{\partial t}\ ,\
\cfrac{\partial u_{1}}{\partial x_{2}}\ ,\ \cfrac{\partial
u_{1}}{\partial x_{3}}\ ,\ u_{1}\ ,\ \cfrac{\partial u_{2}
}{\partial x_{1}}\ ,\ \cfrac{\partial u_{2}}{\partial x_{2}}\ ,\
\cfrac{\partial u_{2}}{\partial x_{3}}\ ,\ u_{2}\ ,\ \cfrac{\partial
u_{3} }{\partial x_{1}}\ ,\ \cfrac{\partial u_{3}}{\partial x_{2}}\
,\ \cfrac{\partial u_{3}}{\partial x_{3}}\ ,\ u_{3}\ ,\  p\ )^{T}\
,$$ $$H=(-\alpha_{1},\ -\alpha_{2},\ -\alpha_{3},\ -\alpha_{4},\
e_{1},\ e_{2},\ e_{3},\ e_{4},\ e_{5},\ e_{6},\ e_{7},\ e_{8},\
e_{9},\ e_{10},\ e_{11},\ e_{16})^{T}\ ,$$ $h=(0,\ F_{1},\ F_{2},\
F_{3},\ 0_{12})^{T}\ ,$ then we can get\ $U=HZ+h$\ , moreover we
assume
\begin{eqnarray*}&&H_{0}=(e_{21},\ e_{17},\ e_{27},\ e_{37},\ e_{22},\ e_{23},\
-\alpha_{2},\ e_{31},\ e_{32},\ e_{33},\ -\alpha_{3},\ e_{41},\
e_{42},\ e_{43},\ -\alpha_{4},\ e_{12})^{T}\ ,\\&& H_{1}=(e_{18},\
e_{21},\ e_{31},\ e_{41},\ e_{24},\ e_{25},\ -\alpha_{1},\ e_{28},\
e_{34},\ e_{35},\ e_{4},\ e_{38},\ e_{44},\ e_{45},\ e_{8},\
e_{13})^{T}\ ,\\&& H_{2}=(e_{24},\ e_{22},\ e_{32},\ e_{42},\
e_{19},\ e_{26},\ e_{1},\ e_{34},\ e_{29},\ e_{36},\ e_{5},\
e_{44},\ e_{39},\ e_{46},\ e_{9},\ e_{14})^{T}\ ,\\&&
H_{3}=(e_{25},\ e_{23},\ e_{33},\ e_{43},\ e_{26},\ e_{20},\ e_{2},\
e_{35},\ e_{36},\ e_{30},\ e_{6},\ e_{45},\ e_{46},\ e_{40},\
e_{10},\ e_{15})^{T}\ ,\end{eqnarray*} $h_{0}=(0_{6},\ F_{1},\
0_{3},\ F_{2},\ 0_{3},\ F_{3},\ 0)^{T}$\ , then we can get\[\label
{A-Problem}\begin{cases}\cfrac{\partial U}{\partial
t}=H_{0}Z+h_{0}\\\cfrac{\partial U}{\partial x_{i}}=H_{i}Z\ ,\ i=1,\
2,\ 3,\end{cases}
\]or
\[\label
{A-Problem}\begin{cases}\cfrac{\partial (HZ+h)}{\partial
t}=H_{0}Z+h_{0}\\\cfrac{\partial (HZ+h)}{\partial x_{i}}=H_{i}Z\ ,\
i=1,\ 2,\ 3.\end{cases}
\]Now we have converted the Navier-Stokes
equations into the simultaneous of the first order linear partial
differential equations with constant coefficients (3.3) and\ $9$\
quadratic polynomial equations (3.1).\\ And we get a necessary
condition for the existence of the solution of the Navier-Stokes
equations as follows,\ $\exists\ 55$\ dimensional four variables
functional vector\ $Z$\ , such that\ $Z$\ satisfies the first order
linear partial differential equations with constant coefficients
(3.3) and\ $9$\ quadratic polynomial equations (3.1).\\In fact, this
condition is also sufficient. If\ $Z$\ satisfies the above first
order linear partial differential equations with constant
coefficients (3.3) and\ $9$\ quadratic polynomial equations (3.1),
we let$$(u_{1},\ u_{2},\ u_{3},\ p)^{T}=(e_{3},\ e_{7},\ e_{11},\
e_{16})^{T}Z\ ,$$then from (3.3) we can get\ $U=HZ+h$\ ,
where$$U=(\cfrac{\partial u_{1}}{\partial x_{1}}\ ,\ \cfrac{\partial
u_{1}}{\partial t}\ ,\ \cfrac{\partial u_{2} }{\partial t}\ ,\
\cfrac{\partial u_{3}}{\partial t}\ ,\ \cfrac{\partial
u_{1}}{\partial x_{2}}\ ,\ \cfrac{\partial u_{1}}{\partial x_{3}}\
,\ u_{1}\ ,\ \cfrac{\partial u_{2} }{\partial x_{1}}\ ,\
\cfrac{\partial u_{2}}{\partial x_{2}}\ ,\ \cfrac{\partial
u_{2}}{\partial x_{3}}\ ,\ u_{2}\ ,\ \cfrac{\partial u_{3}
}{\partial x_{1}}\ ,\ \cfrac{\partial u_{3}}{\partial x_{2}}\ ,\
\cfrac{\partial u_{3}}{\partial x_{3}}\ ,\ u_{3}\ ,\  p\ )^{T}\ ,$$
again from (3.3) we can learn that\ $\cfrac{\partial U}{\partial
t}=H_{0}Z+h_{0}\ ,\ \cfrac{\partial U}{\partial x_{i}}=H_{i}Z\ ,\
i=1,\ 2,\ 3$\ , and from (3.1) we can get\ $X=X_{0}+A(\eta)Z$\ ,
where\ $X$\ includes\ $u_{1},\ u_{2},\ u_{3},\ p$\ and all their
first order partial derivative, and all the second order partial
derivative of\ $u_{1},\ u_{2},\ u_{3}$\ , and all the products which
they are in the Navier-Stokes equations, then we can get\
$AX=\beta$\ , this is just the Navier-Stokes equations.
\\Hence we can get the corollary as follows.
 \begin{corollary} \label{corollary1} If we assume\ $u_{1},\ u_{2},\ u_{3}\in C^{2}\ ,\ p\in
C^{1}\ ,\ F_{1},\ F_{2},\ F_{3}\in C^{1}$\ , then a necessary and
sufficient condition for the existence of the solution for the
Navier-Stokes equations is that\ $\exists\ 55$\ dimensional four
variables functional vector\ $Z$\ , such that\ $Z$\ satisfies the
first order linear partial differential equations with constant
coefficients (3.3) and\ $9$\ quadratic equations
(3.1).\end{corollary}Under this circumstance,\ $(u_{1},\ u_{2},\
u_{3},\ p)^{T}=(e_{3},\ e_{7},\ e_{11},\ e_{16})^{T}Z$\ is the
solution of the Navier-Stokes equations, here\ $e_{i}$\ is the\
$ith$\ $55$\
dimensional unit coordinate vector,\ $1\leq i\leq55$\ .\\\\
{\it Proof of theorem2-3}. (1) Under the assumption (1.1) and from
the theorem 2.2, we can learn that when\ $t$\ is not in\ $[0\ ,\
T]$\ , or\ $(x_{1},\ x_{2},\ x_{3})$\ is not in\ $K_{1}$\ ,\
$Z\equiv0$\ , hence\ $Z$\ can do the Fourier transformation with\
$t\ ,\ x_{1}\ ,\ x_{2}\ ,\ x_{3}$\ , we use the Fourier
transformation to convert\ (2.1)\ into the linear equations as
follows.\[\label {A-Problem}\left(
             \begin{array}{c}
               i\xi_{0}H-H_{0}\\
               i\xi_{1}H-H_{1} \\
               i\xi_{2}H-H_{2}\\
               i\xi_{3}H-H_{3} \\
             \end{array}
           \right)F(Z)=\left(
                         \begin{array}{c}
                           F(h_{0})-i\xi_{0}F(h) \\
                           -i\xi_{1}F(h) \\
                           -i\xi_{2}F(h) \\
                           -i\xi_{3}F(h) \\
                         \end{array}
                       \right)\ ,\]
                       where\begin{eqnarray*}&&F(Z)=\int_{R^{4}}Ze^{-i\xi_{0}t-i\sum_{j=1}^{3}\xi_{j}x_{j}}dtdx_{1}dx_{2}dx_{3}\
                       ,\\&& F(h)=\int_{R^{4}}he^{-i\xi_{0}t-i\sum_{j=1}^{3}\xi_{j}x_{j}}dtdx_{1}dx_{2}dx_{3}\
                       ,\\
                       &&F(h_{0})=\int_{R^{4}}h_{0}e^{-i\xi_{0}t-i\sum_{j=1}^{3}\xi_{j}x_{j}}dtdx_{1}dx_{2}dx_{3}\
                       .\end{eqnarray*}
                       We assume that$$B=\left(
             \begin{array}{c}
               i\xi_{0}H-H_{0}\\
               i\xi_{1}H-H_{1} \\
               i\xi_{2}H-H_{2}\\
               i\xi_{3}H-H_{3} \\
             \end{array}
           \right)_{64\times55}\ ,\ G=\left(
                         \begin{array}{c}
                           F(h_{0})-i\xi_{0}F(h) \\
                           -i\xi_{1}F(h) \\
                           -i\xi_{2}F(h) \\
                           -i\xi_{3}F(h) \\
                         \end{array}
                       \right)_{64\times1}\ ,$$next we solve the linear
                       equations\ $BY=G$\ , where\
                       $Y=(y_{j})_{55\times1}$\ . We can get\
                       $rank(B)=46$\ , when\ $(\xi_{1}\ ,\ \xi_{2}\ ,\
                       \xi_{3})\neq 0$\ , moreover the first\ $46$\
                       columns are linear independent.\\We write out
                       all the rows of the matrix
                       B.
\begin{eqnarray*}
&&-i\xi_{0}\alpha_{1}^{T}-e_{21}^{T}\ ,\
-i\xi_{1}\alpha_{1}^{T}-e_{18}^{T}\ ,\
-i\xi_{2}\alpha_{1}^{T}-e_{24}^{T}\ ,\
-i\xi_{3}\alpha_{1}^{T}-e_{25}^{T}\ ,
\\&&-i\xi_{0}\alpha_{2}^{T}-e_{17}^{T}\ ,\ -i\xi_{1}\alpha_{2}^{T}-e_{21}^{T}\ ,\ -i\xi_{2}\alpha_{2}^{T}-e_{22}^{T}\ ,\ -i\xi_{3}\alpha_{2}^{T}-e_{23}^{T}\ ,
\\&&-i\xi_{0}\alpha_{3}^{T}-e_{27}^{T}\ ,\
-i\xi_{1}\alpha_{3}^{T}-e_{31}^{T}\ ,\
-i\xi_{2}\alpha_{3}^{T}-e_{32}^{T}\ ,\
-i\xi_{3}\alpha_{3}^{T}-e_{33}^{T}\ ,
\\&&-i\xi_{0}\alpha_{4}^{T}-e_{37}^{T}\ ,\ -i\xi_{1}\alpha_{4}^{T}-e_{41}^{T}\ ,\ -i\xi_{2}\alpha_{4}^{T}-e_{42}^{T}\ ,\ -i\xi_{3}\alpha_{4}^{T}-e_{43}^{T}\ ,
\\&&i\xi_{0}e_{1}^{T}-e_{22}^{T}\ ,\ i\xi_{1}e_{1}^{T}-e_{24}^{T}\ ,\
i\xi_{2}e_{1}^{T}-e_{19}^{T}\ ,\ i\xi_{3}e_{1}^{T}-e_{26}^{T}\ ,
\\&&i\xi_{0}e_{2}^{T}-e_{23}^{T}\ ,\ i\xi_{1}e_{2}^{T}-e_{25}^{T}\ ,\ i\xi_{2}e_{2}^{T}-e_{26}^{T}\ ,\ i\xi_{3}e_{2}^{T}-e_{20}^{T}\ ,
\\&&i\xi_{0}e_{3}^{T}+\alpha_{2}^{T}\ ,\ i\xi_{1}e_{3}^{T}+\alpha_{1}^{T}\ ,\ i\xi_{2}e_{3}^{T}-e_{1}^{T}\ ,\ i\xi_{3}e_{3}^{T}-e_{2}^{T}\ ,
\\&&i\xi_{0}e_{4}^{T}-e_{31}^{T}\ ,\ i\xi_{1}e_{4}^{T}-e_{28}^{T}\ ,\ i\xi_{2}e_{4}^{T}-e_{34}^{T}\ ,\ i\xi_{3}e_{4}^{T}-e_{35}^{T}\ ,
\\&&i\xi_{0}e_{5}^{T}-e_{32}^{T}\ ,\ i\xi_{1}e_{5}^{T}-e_{34}^{T}\ ,\ i\xi_{2}e_{5}^{T}-e_{29}^{T}\ ,\ i\xi_{3}e_{5}^{T}-e_{36}^{T}\ ,
\\&&i\xi_{0}e_{6}^{T}-e_{33}^{T}\ ,\ i\xi_{1}e_{6}^{T}-e_{35}^{T}\ ,\ i\xi_{2}e_{6}^{T}-e_{36}^{T}\ ,\ i\xi_{3}e_{6}^{T}-e_{30}^{T}\ ,
\\&&i\xi_{0}e_{7}^{T}+\alpha_{3}^{T}\ ,\ i\xi_{1}e_{7}^{T}-e_{4}^{T}\ ,\ i\xi_{2}e_{7}^{T}-e_{5}^{T}\ ,\ i\xi_{3}e_{7}^{T}-e_{6}^{T}\ ,
\\&&i\xi_{0}e_{8}^{T}-e_{41}^{T}\ ,\ i\xi_{1}e_{8}^{T}-e_{38}^{T}\ ,\ i\xi_{2}e_{8}^{T}-e_{44}^{T}\ ,\ i\xi_{3}e_{8}^{T}-e_{45}^{T}\ ,
\\&&i\xi_{0}e_{9}^{T}-e_{42}^{T}\ ,\ i\xi_{1}e_{9}^{T}-e_{44}^{T}\ ,\
i\xi_{2}e_{9}^{T}-e_{39}^{T}\ ,\ i\xi_{3}e_{9}^{T}-e_{46}^{T}\ ,
\\&&i\xi_{0}e_{10}^{T}-e_{43}^{T}\ ,\ i\xi_{1}e_{10}^{T}-e_{45}^{T}\ ,\ i\xi_{2}e_{10}^{T}-e_{46}^{T}\ ,\ i\xi_{3}e_{10}^{T}-e_{40}^{T}\ ,
\\&&i\xi_{0}e_{11}^{T}+\alpha_{4}^{T}\ ,\ i\xi_{1}e_{11}^{T}-e_{8}^{T}\ ,\ i\xi_{2}e_{11}^{T}-e_{9}^{T}\ ,\ i\xi_{3}e_{11}^{T}-e_{10}^{T}\ ,
\\&&i\xi_{0}e_{16}^{T}-e_{12}^{T}\ ,\ i\xi_{1}e_{16}^{T}-e_{13}^{T}\ ,\ i\xi_{2}e_{16}^{T}-e_{14}^{T}\ ,\ i\xi_{3}e_{16}^{T}-e_{15}^{T}\
,
\end{eqnarray*}where\ $e_{i}$\ is the\ $ith$\ $55$\ dimensional unit
coordinate vector,\ $1 \leq i \leq 55$\ ,
and\begin{eqnarray*}&&\alpha_{1}=(0_{4},\ 1,\ 0_{4},\  1,\  0,\ 0,\
0,\  0,\ 0,\ 0,\ 0_{39})^{T},\\&& \alpha_{2}=(0_{12},\  \tau,\
0_{4},\  -\mu,\ -\mu,\  -\mu,\ 0_{26},\  1,\  1,\  1,\  0_{6}
)^{T},\\&& \alpha_{3}=(0_{13},\ \tau,\   0_{13},\   -\mu,\   -\mu,\
-\mu,\   0_{19},\   1,\   1,\ 1,\  0_{3} )^{T},\\&&
\alpha_{4}=(0_{14},\   \tau,\   0_{22},\ -\mu,\   -\mu,\   -\mu,\
0_{11},\   0,\   1,\ 1,\   1)^{T}\ .\end{eqnarray*}First we let\
$y_{47+j-1}=0\ ,\ 1\leq j\leq 9$\ , we will show that the solution
of\ $BY=0$\ is only\ $0$\ , when\ $(\xi_{1},\ \xi_{2},\ \xi_{3})\neq
0$\ , hence the first\ $46$\
                       columns are linear independent, and\ $rank(B)\geq
                       46$\ .\\From the\ $7$th row we can get\
                       $y_{1}=i\xi_{2}y_{3},\ y_{2}=i\xi_{3}y_{3},\ y_{5}+y_{10}=-i\xi_{1}y_{3}$\
                       , and from the first, the second, the\ $5$th
                       , the\ $6$th rows, we can get\
                       $y_{18}=(i\xi_{1})^{2}y_{3},\ y_{19}=(i\xi_{2})^{2}y_{3},\
                       y_{20}=(i\xi_{3})^{2}y_{3},\ y_{13}=ay_{3}\ ,$
                       where
                       $$a=\cfrac{\mu[(i\xi_{1})^{2}+(i\xi_{2})^{2}+(i\xi_{3})^{2}]-i\xi_{0}}{\tau}\neq
                       0\ , $$when\ $(\xi_{1},\ \xi_{2},\ \xi_{3})\neq
0$\ , moreover\ $y_{17}=(i\xi_{0})^{2}y_{3}\ ,$
$y_{21}=i\xi_{0}i\xi_{1}y_{3},\ y_{22}=i\xi_{0}i\xi_{2}y_{3},\
y_{23}=i\xi_{0}i\xi_{3}y_{3},$\\$ y_{24}=i\xi_{1}i\xi_{2}y_{3},\
y_{25}=i\xi_{1}i\xi_{3}y_{3},\ y_{26}=i\xi_{2}i\xi_{3}y_{3}$\
.\\From the\ $11$th row we can get\
                       $y_{4}=i\xi_{1}y_{7},\ y_{5}=i\xi_{2}y_{7},\ y_{6}=i\xi_{3}y_{7}$\
                       , and from the\ $8$th, the third, the\ $9$th
                       , the\ $10$th rows, we can get\
                       $y_{28}=(i\xi_{1})^{2}y_{7},\ y_{29}=(i\xi_{2})^{2}y_{7},\
                       y_{30}=(i\xi_{3})^{2}y_{7},\ y_{14}=ay_{7}\ ,$
                       \ moreover\ $y_{27}=(i\xi_{0})^{2}y_{7}\ ,$
$y_{31}=i\xi_{0}i\xi_{1}y_{7},\ y_{32}=i\xi_{0}i\xi_{2}y_{7},\
y_{33}=i\xi_{0}i\xi_{3}y_{7},\ y_{34}=i\xi_{1}i\xi_{2}y_{7},\
y_{35}=i\xi_{1}i\xi_{3}y_{7}$,\\$ y_{36}=i\xi_{2}i\xi_{3}y_{7}$\
.\\From the\ $15$th row we can get\
                       $y_{8}=i\xi_{1}y_{11},\ y_{9}=i\xi_{2}y_{11},\ y_{10}=i\xi_{3}y_{11}$\
                       , and from the\ $12$th, the\ $13$th, the\ $14$th
                       ,  the\ $4$th rows, we can get\
                       $y_{38}=(i\xi_{1})^{2}y_{11},\ y_{39}=(i\xi_{2})^{2}y_{11},\
                       y_{40}=(i\xi_{3})^{2}y_{11},\ y_{15}=ay_{11}\ ,$
                       \ moreover\ $y_{37}=(i\xi_{0})^{2}y_{11}\ ,$
$y_{41}=i\xi_{0}i\xi_{1}y_{11},\ y_{42}=i\xi_{0}i\xi_{2}y_{11},\
y_{43}=i\xi_{0}i\xi_{3}y_{11},\ y_{44}=i\xi_{1}i\xi_{2}y_{11}$,\\$
y_{45}=i\xi_{1}i\xi_{3}y_{11},\ y_{46}=i\xi_{2}i\xi_{3}y_{11}$.\\And
from the last row, we can get\ $y_{12}=i\xi_{0}y_{16},\
y_{13}=i\xi_{1}y_{16},\ y_{14}=i\xi_{2}y_{16},\
y_{15}=i\xi_{3}y_{16}$\ .\\ Because\ $y_{13}=ay_{3},\
y_{14}=ay_{7},\ y_{15}=ay_{11}$\ , we can get
$$y_{3}=\cfrac{i\xi_{1}}{a}y_{16},\ y_{7}=\cfrac{i\xi_{2}}{a}y_{16},\
y_{11}=\cfrac{i\xi_{3}}{a}y_{16}\ .$$ From\
$(y_{5}+y_{10})=y_{5}+y_{10}$\ , we can get\
$-i\xi_{1}y_{3}=i\xi_{2}y_{7}+i\xi_{3}y_{11}$\ . This is equal to
the
following.$$[\cfrac{(i\xi_{1})^{2}}{a}+\cfrac{(i\xi_{2})^{2}}{a}+\cfrac{(i\xi_{3})^{2}}{a}]y_{16}=0\
.$$Hence\ $y_{16}=0$\ , when\ $(\xi_{1},\ \xi_{2},\ \xi_{3})\neq 0$\
, and we can get\ $y_{3}=0,\ y_{7}=0,\ y_{11}=0$\ , the solution of\
$BY=0$\ is only\ $0$\ , the first\ $46$\
                       columns are linear independent, and\ $rank(B)\geq
                       46$\ .\\Next we will in turn let\ $y_{47+j-1}=1,\
                       y_{k}=0,\ 47\leq k\leq 55,\ k\neq 47+j-1,\
                       1\leq j\leq 9$\ , we can work out\ $9$\ linear
                       independent solutions of\
$BY=0$\ , we assume they are\ $\eta_{j},\ 1\leq j\leq 9$\ . This
means that\ $rank(B)\leq
                       46$\ , hence\ $rank(B)=
                       46\ ,$ when\ $(\xi_{1},\ \xi_{2},\ \xi_{3})\neq 0$\
. \\If we let\ $y_{47+j-1}=1,\
                       y_{k}=0,\ 47\leq k\leq 55,\ k\neq 47+j-1,\
                       1\leq j\leq 3$\ , we only need to notice
                       that\ $\alpha_{2}^{T}Y$\ will change into\
                       $\alpha_{2}^{T}Y+1$\ , and\
                       $y_{13}=i\xi_{1}y_{16}=ay_{3}-\cfrac{1}{\tau}$\
                       , then we can get\
                       $y_{3}=\cfrac{i\xi_{1}}{a}y_{16}+\cfrac{1}{a\tau}$\
                       , and\ $y_{7}=\cfrac{i\xi_{2}}{a}y_{16}$\ ,\ $y_{11}=\cfrac{i\xi_{3}}{a}y_{16}$\
                       . Hence we can get the
                       following,$$by_{16}+\cfrac{i\xi_{1}}{a\tau}=0\
                       ,\ \mbox{or}\
                       y_{16}=-\cfrac{i\xi_{1}}{ab\tau}\ .$$Moreover we
                       can work out\ $y_{3}\ ,\ y_{7}\ ,\ y_{11}$\
                       and\ $\eta_{j}\ ,\ 1\leq j\leq 3$\ .
\\If we let\ $y_{47+j-1}=1,\
                       y_{k}=0,\ 47\leq k\leq 55,\ k\neq 47+j-1,\
                       4\leq j\leq 6$\ , we only need to notice
                       that\ $\alpha_{3}^{T}Y$\ will change into\
                       $\alpha_{3}^{T}Y+1$\ , and\
                       $y_{14}=i\xi_{2}y_{16}=ay_{7}-\cfrac{1}{\tau}$\
                       , then we can get\
                       $y_{7}=\cfrac{i\xi_{2}}{a}y_{16}+\cfrac{1}{a\tau}$\
                       , and\ $y_{3}=\cfrac{i\xi_{1}}{a}y_{16}$\ ,\ $y_{11}=\cfrac{i\xi_{3}}{a}y_{16}$\
                       . Hence we can get the
                       following,$$by_{16}+\cfrac{i\xi_{2}}{a\tau}=0\
                       ,\ \mbox{or}\
                       y_{16}=-\cfrac{i\xi_{2}}{ab\tau}\ .$$Moreover we
                       can work out\ $y_{3}\ ,\ y_{7}\ ,\ y_{11}$\
                       and\ $\eta_{j}\ ,\ 4\leq j\leq 6$\ .
\\If we let\ $y_{47+j-1}=1,\
                       y_{k}=0,\ 47\leq k\leq 55,\ k\neq 47+j-1,\
                       7\leq j\leq 9$\ , we only need to notice
                       that\ $\alpha_{4}^{T}Y$\ will change into\
                       $\alpha_{4}^{T}Y+1$\ , and\
                       $y_{15}=i\xi_{3}y_{16}=ay_{11}-\cfrac{1}{\tau}$\
                       , then we can get\
                       $y_{11}=\cfrac{i\xi_{3}}{a}y_{16}+\cfrac{1}{a\tau}$\
                       , and\ $y_{3}=\cfrac{i\xi_{1}}{a}y_{16}$\ ,\ $y_{7}=\cfrac{i\xi_{2}}{a}y_{16}$\
                       . Hence we can get the
                       following,$$by_{16}+\cfrac{i\xi_{3}}{a\tau}=0\
                       ,\ \mbox{or}\
                       y_{16}=-\cfrac{i\xi_{3}}{ab\tau}\ .$$Moreover we
                       can work out\ $y_{3}\ ,\ y_{7}\ ,\ y_{11}$\
                       and\ $\eta_{j}\ ,\ 7\leq j\leq 9$\ .\\
We write\ $\eta_{j}\ ,\ 1\leq j\leq 9$\ , as
follows.\begin{eqnarray*}\eta_{j}&=&(i\xi_{2}y_{3},\ i\xi_{3}y_{3},\
y_{3},\ i\xi_{1}y_{7},\ i\xi_{2}y_{7},\ i\xi_{3}y_{7},\ y_{7},\
i\xi_{1}y_{11},\ i\xi_{2}y_{11},\ i\xi_{3}y_{11},\ y_{11},\\&&
i\xi_{0}y_{16},\ i\xi_{1}y_{16},\ i\xi_{2}y_{16},\ i\xi_{3}y_{16},\
y_{16},\ (i\xi_{0})^{2}y_{3},\ (i\xi_{1})^{2}y_{3},\
(i\xi_{2})^{2}y_{3},\ (i\xi_{3})^{2}y_{3},\
i\xi_{0}i\xi_{1}y_{3},\\&& i\xi_{0}i\xi_{2}y_{3},\
i\xi_{0}i\xi_{3}y_{3},\ i\xi_{1}i\xi_{2}y_{3},\
i\xi_{1}i\xi_{3}y_{3},\ i\xi_{2}i\xi_{3}y_{3},\
(i\xi_{0})^{2}y_{7},\ (i\xi_{1})^{2}y_{7},\ (i\xi_{2})^{2}y_{7},\
(i\xi_{3})^{2}y_{7},\\&& i\xi_{0}i\xi_{1}y_{7},\
i\xi_{0}i\xi_{2}y_{7},\ i\xi_{0}i\xi_{3}y_{7},\
i\xi_{1}i\xi_{2}y_{7},\ i\xi_{1}i\xi_{3}y_{7},\
i\xi_{2}i\xi_{3}y_{7},\ (i\xi_{0})^{2}y_{11},\
(i\xi_{1})^{2}y_{11},\ (i\xi_{2})^{2}y_{11},\\&&
(i\xi_{3})^{2}y_{11},\ i\xi_{0}i\xi_{1}y_{11},\
i\xi_{0}i\xi_{2}y_{11},\ i\xi_{0}i\xi_{3}y_{11},\
i\xi_{1}i\xi_{2}y_{11},\ i\xi_{1}i\xi_{3}y_{11},\
i\xi_{2}i\xi_{3}y_{11},\ e_{j}^{T})^{T}\ ,\end{eqnarray*}here\
$e_{j}$\ is the\ $jth$\ $9$\ dimensional unit coordinate vector,\ $1
\leq j \leq 9$\ , moreover when\ $j=1,\ 2,\ 3$\
,$$y_{3}=\cfrac{\xi_{1}^{2}}{a^{2}b\tau}+\cfrac{1}{a\tau}\ ,\
y_{7}=\cfrac{-i\xi_{1}i\xi_{2}}{a^{2}b\tau}\ ,\
y_{11}=\cfrac{-i\xi_{1}i\xi_{3}}{a^{2}b\tau}\ ,\
y_{16}=\cfrac{-i\xi_{1}}{ab\tau}\ ,$$where\
$b=\cfrac{(i\xi_{1})^{2}+(i\xi_{2})^{2}+(i\xi_{3})^{2}}{a}\neq 0$\ ,
when\ $(\xi_{1},\ \xi_{2},\ \xi_{3})\neq 0$\ , when\ $j=4,\ 5,\ 6$\
,$$y_{7}=\cfrac{\xi_{2}^{2}}{a^{2}b\tau}+\cfrac{1}{a\tau}\ ,\
y_{3}=\cfrac{-i\xi_{1}i\xi_{2}}{a^{2}b\tau}\ ,\
y_{11}=\cfrac{-i\xi_{2}i\xi_{3}}{a^{2}b\tau}\ ,\
y_{16}=\cfrac{-i\xi_{2}}{ab\tau}\ ,$$when\ $j=7,\ 8,\ 9$\
,$$y_{11}=\cfrac{\xi_{3}^{2}}{a^{2}b\tau}+\cfrac{1}{a\tau}\ ,\
y_{3}=\cfrac{-i\xi_{1}i\xi_{3}}{a^{2}b\tau}\ ,\
y_{7}=\cfrac{-i\xi_{2}i\xi_{3}}{a^{2}b\tau}\ ,\
y_{16}=\cfrac{-i\xi_{3}}{ab\tau}\ .$$Finally we work out a
particular solution\ $Y_{1}$\ of\ $BY=G$\ . We can see
that\begin{eqnarray*}G&=&(0,\ -i\xi_{0}F(F_{1}),\
-i\xi_{0}F(F_{2}),\ -i\xi_{0}F(F_{3}),\ 0_{2},\ F(F_{1}),\ 0_{3},\
F(F_{2}),\ 0_{3},\ F(F_{3}),\ 0,\ 0,\\&& i\xi_{1}F(F_{1}),\
-i\xi_{1}F(F_{2}),\ -i\xi_{1}F(F_{3}),\ 0_{12},\ 0,\
i\xi_{2}F(F_{1}),\ -i\xi_{2}F(F_{2}),\ -i\xi_{2}F(F_{3}),\
0_{12},\\&& 0,\ i\xi_{3}F(F_{1}),\ -i\xi_{3}F(F_{2}),\
-i\xi_{3}F(F_{3}),\ 0_{12})^{T}\ ,\end{eqnarray*} If we let\
$y_{47+j-1}=0,\ 1\leq j\leq 9$\ , then we can get\
$y_{13}=i\xi_{1}y_{16}=ay_{3}+\cfrac{F(F_{1})}{\tau}$\ ,\
$y_{14}=i\xi_{2}y_{16}=ay_{7}+\cfrac{F(F_{2})}{\tau}$\ , and\
$y_{15}=i\xi_{3}y_{16}=ay_{11}+\cfrac{F(F_{3})}{\tau}$\ . Hence we
can get the
following,$$by_{16}=\cfrac{i\xi_{1}F(F_{1})+i\xi_{2}F(F_{2})+i\xi_{3}F(F_{3})}{a\tau}\
,\ \mbox{or}\
y_{16}=\cfrac{i\xi_{1}F(F_{1})+i\xi_{2}F(F_{2})+i\xi_{3}F(F_{3})}{ab\tau}\
.$$Moreover we can work out\ $y_{3}\ ,\ y_{7}\ ,\ y_{11}$\
                       and\ $Y_{1}$\ .\\
We write\ $Y_{1}$\ as follows.
\begin{eqnarray*}Y_{1}&=&(i\xi_{2}y_{3},\ i\xi_{3}y_{3},\ y_{3},\
i\xi_{1}y_{7},\ i\xi_{2}y_{7},\ i\xi_{3}y_{7},\ y_{7},\
i\xi_{1}y_{11},\ i\xi_{2}y_{11},\ i\xi_{3}y_{11},\ y_{11},\\&&
i\xi_{0}y_{16},\ i\xi_{1}y_{16},\ i\xi_{2}y_{16},\ i\xi_{3}y_{16},\
y_{16},\ (i\xi_{0})^{2}y_{3},\ (i\xi_{1})^{2}y_{3},\
(i\xi_{2})^{2}y_{3},\ (i\xi_{3})^{2}y_{3},\
i\xi_{0}i\xi_{1}y_{3},\\&& i\xi_{0}i\xi_{2}y_{3},\
i\xi_{0}i\xi_{3}y_{3},\ i\xi_{1}i\xi_{2}y_{3},\
i\xi_{1}i\xi_{3}y_{3},\ i\xi_{2}i\xi_{3}y_{3},\
(i\xi_{0})^{2}y_{7},\ (i\xi_{1})^{2}y_{7},\ (i\xi_{2})^{2}y_{7},\
(i\xi_{3})^{2}y_{7},\\&& i\xi_{0}i\xi_{1}y_{7},\
i\xi_{0}i\xi_{2}y_{7},\ i\xi_{0}i\xi_{3}y_{7},\
i\xi_{1}i\xi_{2}y_{7},\ i\xi_{1}i\xi_{3}y_{7},\
i\xi_{2}i\xi_{3}y_{7},\ (i\xi_{0})^{2}y_{11},\
(i\xi_{1})^{2}y_{11},\ (i\xi_{2})^{2}y_{11},\\&&
(i\xi_{3})^{2}y_{11},\ i\xi_{0}i\xi_{1}y_{11},\
i\xi_{0}i\xi_{2}y_{11},\ i\xi_{0}i\xi_{3}y_{11},\
i\xi_{1}i\xi_{2}y_{11},\ i\xi_{1}i\xi_{3}y_{11},\
i\xi_{2}i\xi_{3}y_{11},\ 0_{9\times1}^{T})^{T}\
,\end{eqnarray*}where\begin{eqnarray*}&&
y_{3}=\cfrac{i\xi_{1}(i\xi_{1}F(F_{1})+i\xi_{2}F(F_{2})+i\xi_{3}F(F_{3}))}{a^{2}b\tau}-\cfrac{F(F_{1})}{a\tau}\
,\\&&y_{7}=\cfrac{i\xi_{2}(i\xi_{1}F(F_{1})+i\xi_{2}F(F_{2})+i\xi_{3}F(F_{3}))}{a^{2}b\tau}-\cfrac{F(F_{2})}{a\tau}\
,
\\&&y_{11}=\cfrac{i\xi_{3}(i\xi_{1}F(F_{1})+i\xi_{2}F(F_{2})+i\xi_{3}F(F_{3}))}{a^{2}b\tau}-\cfrac{F(F_{3})}{a\tau}\
,\\&&y_{16}=\cfrac{i\xi_{1}F(F_{1})+i\xi_{2}F(F_{2})+i\xi_{3}F(F_{3})}{ab\tau}\
.\end{eqnarray*}And we can test such\ $Y_{1}$\ satisfies\
$BY_{1}=G$\ . \\If we assume\ $A_{1}(\eta)=(\eta_{1} ,\ \eta_{2} ,\
\eta_{3} ,\ \eta_{4} ,\ \eta_{5},\ \eta_{6} ,\ \eta_{7} ,\ \eta_{8}
,\ \eta_{9})$\ , then we can get the explicit general solutions of
(3.4) as follows,$$F(Z)=Y_{1}+A_{1}(\eta)Z_{1}\ , \mbox{or}\
                       Z=F^{-1}(Y_{1}+A_{1}(\eta)Z_{1})\ ,$$
 where\ $Z_{1}=(Z_{1j})_{9\times1}$\ , moreover we assume$$F^{-1}(Y_{1}+A_{1}(\eta)Z_{1}I_{\{(\xi_{1},\ \xi_{2},\
\xi_{3})\neq 0\}})=F^{-1}(Y_{1}+A_{1}(\eta)Z_{1})\ ,$$ and\
$Z_{1}\in\ \Omega_{1}$\ . Because\ $Z$\ need to
satisfy$$HZ=HZI_{K_{1}^{\prime}}\in C^{1}(K_{1}^{\prime})\ ,$$ we
can get
\begin{eqnarray*}&&\Omega_{1}=\{Z_{1}|H[F^{-1}(Y_{1}+A_{1}(\eta)Z_{1})]=H[F^{-1}(Y_{1}+A_{1}(\eta)Z_{1})]I_{K_{1}^{\prime}}\in
C^{1}(K_{1}^{\prime}).\}\ , \mbox{where}\\&&H=(-\alpha_{1},\
-\alpha_{2},\ -\alpha_{3},\ -\alpha_{4},\ e_{1},\ e_{2},\ e_{3},\
e_{4},\ e_{5},\ e_{6},\ e_{7},\ e_{8},\ e_{9},\ e_{10},\ e_{11},\
e_{16})^{T}\ .
\end{eqnarray*} (2)We can learn from (1) that\
$F(e_{47+j-1}^{T}Z)=Z_{1j}\ ,\ 1\leq j\leq 9$\ , hence we can
convert\ (2.2)\ into the integral equations as follows,
\[\label {A-Problem}\begin{cases}Z_{11}=F[F^{-1}(e_{3}^{T}(Y_{1}+A_{1}(\eta)Z_{1}))F^{-1}(-\alpha_{1}^{T}(Y_{1}+A_{1}(\eta)Z_{1}))]=f_{1}(Z_{1})\ ,
\\Z_{12}=F[F^{-1}(e_{7}^{T}(Y_{1}+A_{1}(\eta)Z_{1}))F^{-1}(e_{1}^{T}(Y_{1}+A_{1}(\eta)Z_{1}))]=f_{2}(Z_{1})\ ,\\
Z_{13}=F[F^{-1}(e_{11}^{T}(Y_{1}+A_{1}(\eta)Z_{1}))F^{-1}(e_{2}^{T}(Y_{1}+A_{1}(\eta)Z_{1}))]=f_{3}(Z_{1})\
,\\
Z_{14}=F[F^{-1}(e_{3}^{T}(Y_{1}+A_{1}(\eta)Z_{1}))F^{-1}(e_{4}^{T}(Y_{1}+A_{1}(\eta)Z_{1}))]=f_{4}(Z_{1})\
,\\
Z_{15}=F[F^{-1}(e_{7}^{T}(Y_{1}+A_{1}(\eta)Z_{1}))F^{-1}(e_{5}^{T}(Y_{1}+A_{1}(\eta)Z_{1}))]=f_{5}(Z_{1})\
,\\
Z_{16}=F[F^{-1}(e_{11}^{T}(Y_{1}+A_{1}(\eta)Z_{1}))F^{-1}(e_{6}^{T}(Y_{1}+A_{1}(\eta)Z_{1}))]=f_{6}(Z_{1})\
,\\
Z_{17}=F[F^{-1}(e_{3}^{T}(Y_{1}+A_{1}(\eta)Z_{1}))F^{-1}(e_{8}^{T}(Y_{1}+A_{1}(\eta)Z_{1}))]=f_{7}(Z_{1})\
,\\
Z_{18}=F[F^{-1}(e_{7}^{T}(Y_{1}+A_{1}(\eta)Z_{1}))F^{-1}(e_{9}^{T}(Y_{1}+A_{1}(\eta)Z_{1}))]=f_{8}(Z_{1})\
,\\
Z_{19}=F[F^{-1}(e_{11}^{T}(Y_{1}+A_{1}(\eta)Z_{1}))F^{-1}(e_{10}^{T}(Y_{1}+A_{1}(\eta)Z_{1}))]=f_{9}(Z_{1})\
.\end{cases}
\]This is also the question to find the fixed-points
of\ $f(Z_{1})$\ , where\ $f(Z_{1})=(f_{j}(Z_{1}))_{9\times1}$\ . But
we can not see immediately\ $f(Z_{1})\in\Omega_{1}$\ , even if\
$Z_{1}\in\Omega_{1}$\ . We need to attain this point first. We
assume\ $\Omega_{2}$\ as follows.$$\Omega_{2}=\{Z_{1}|\exists\
h_{1}=h_{1}I_{K_{1}^{\prime}}\in C^{1}(K_{1}^{\prime})\ ,\
\mbox{such that}\ Z_{1}=F(h_{1})\ .\}\ ,$$ where
$$h_{1}=(h_{1j})_{9\times 1}\ ,\ F(h_{1})=\int_{R^{4}}h_{1}e^{-i\xi_{0}t-i\sum_{j=1}^{3}\xi_{j}x_{j}}dtdx_{1}dx_{2}dx_{3}\
.$$ We will prove that\ $\Omega_{2}\subset\Omega_{1}$\ and\
$f(Z_{1})\in\Omega_{2}$\ , if\ $Z_{1}\in\Omega_{2}$\ . We look at a
lemma as follows.
\begin{lemma} \label{Lemma-3.1} $\forall\
h_{1}=h_{1}I_{K_{1}^{\prime}}\in\ C(K_{1}^{\prime})$\ , and\
$h_{1}$\ satisfies the H\"{o}lder condition,\ $\exists\
h_{2}=h_{2}I_{K_{1}^{\prime}}\in C^{1}[0,\ T]\bigcap
C^{\infty}(K_{1})$\ , such that\ $F(h_{1})=a^{2}b\tau F(h_{2})$,
where$$a=\cfrac{\mu[(i\xi_{1})^{2}+(i\xi_{2})^{2}+(i\xi_{3})^{2}]-i\xi_{0}}{\tau}\
,\ b=\cfrac{(i\xi_{1})^{2}+(i\xi_{2})^{2}+(i\xi_{3})^{2}}{a}\
,$$and\ $h_{1}$\ satisfies the H\"{o}lder condition means that\ $
\exists\ C> 0\ , 0<\alpha< 1$\ , such
that$$|h_{1}(x_{t})-h_{1}(x_{t}^{\prime})|=\max_{1\leq j\leq
9}|h_{1j}(x_{t})-h_{1j}(x_{t}^{\prime})|\leq
C|x_{t}-x_{t}^{\prime}|^{\alpha}\ ,\ \forall\ x_{t}\ ,\
x_{t}^{\prime}\ \in\ K_{1}^{\prime}\ .$$
\end{lemma}{\it Proof of lemma3-1}. We see that\
$a^{2}b\tau=(\mu[(i\xi_{1})^{2}+(i\xi_{2})^{2}+(i\xi_{3})^{2}]-i\xi_{0})((i\xi_{1})^{2}+(i\xi_{2})^{2}+(i\xi_{3})^{2})$\
, hence this lemma is equivalent to\ $\forall\
h_{1}=h_{1}I_{K_{1}^{\prime}}\in\ C(K_{1}^{\prime})$\ , and\
$h_{1}$\ satisfies the H\"{o}lder condition,\ $\exists\
h_{2}=h_{2}I_{K_{1}^{\prime}}\in C^{1}[0,\ T]\bigcap
C^{\infty}(K_{1})$\ , such
that\[\mu\sum_{j=1}^{3}\sum_{k=1}^{3}\cfrac{\partial^{4}
h_{2}}{\partial x_{j}^{2}\partial
x_{k}^{2}}-\sum_{j=1}^{3}\cfrac{\partial^{3} h_{2}}{\partial
t\partial x_{j}^{2}}=h_{1}\ ,\ a.e.\]If we let\ $v=\mu\triangle
h_{2}-\partial h_{2}/\partial t$\ , then we convert (3.6) into
follows,\[\begin{cases}\triangle v=h_{1}\ ,\\\mu\triangle
h_{2}-\cfrac{\partial h_{2}}{\partial t}=v\ .\end{cases}\]We see
that (3.7) are the simultaneous of Poisson's equation and the
heat-conduct equation, they are all classical mathematical-physics
equations and\ $h_{1}$\ satisfies the H\"{o}lder condition, the
boundary of\ $K_{1}$\ ,\ $\partial K_{1}$\ satisfies the exterior
ball condition, we know their solutions are all exist. We can write
\ $h_{2}$\ as
follows.\begin{eqnarray*}&&v(t,M_{0})=-\frac{1}{4\pi}\int_{K_{1}}\frac{h_{1}(t,
M)}{r_{MM_{0}}} dx_{1}dx_{2}dx_{3}\ ,\ \mbox{where}\
M=(x_{1},x_{2},x_{3})\ ,\ M_{0}=(x_{10},x_{20},x_{30})\
,\\&&r_{MM_{0}}=\sqrt{(x_{1}-x_{10})^{2}+
(x_{2}-x_{20})^{2}+(x_{3}-x_{30})^{2}}\ ,\\&&h_{2}(t,
M)=(\frac{1}{2\sqrt{\pi\mu}})^{3}\int_{0}^{t}\int_{R^{3}}\frac{v(\tau_{1},
y_{1}, y_{2} , y_{3})}{(\sqrt{t-\tau_{1}})^{3}}\
e^{-\frac{(x_{1}-y_{1})^{2}+(x_{2}-y_{2})^{2}+(x_{3}-y_{3})^{2}}{4\mu(t-\tau_{1})}}
dy_{1}dy_{2}dy_{3}d\tau_{1}\ . \end{eqnarray*} We assume the measure
of the boundary of\ $K_{1}$\ is\ $0$\ , then we can get that\
$h_{2}$\ satisfies (3.6) . Next we will prove\
$h_{2}=h_{2}I_{K_{1}^{\prime}}\in C^{1}[0,\ T]\bigcap
C^{\infty}(K_{1})$\ . \\First we can see\ $v$\ is continuous on the
region\ $K_{1}^{\prime}$\ , and from$$\cfrac{\partial^{k}h_{2}(t,
M)}{\partial x_{1}^{a_{1}}\partial x_{2}^{a_{2}}\partial
x_{3}^{a_{3}}}=(\frac{1}{2\sqrt{\pi\mu}})^{3}\int_{0}^{t}\int_{R^{3}}\frac{v(\tau_{1},
y_{1}, y_{2} , y_{3})}{(\sqrt{t-\tau_{1}})^{3}}\
\cfrac{\partial^{k}e^{-\frac{(x_{1}-y_{1})^{2}+(x_{2}-y_{2})^{2}+(x_{3}-y_{3})^{2}}{4\mu(t-\tau_{1})}}}
{\partial x_{1}^{a_{1}}\partial x_{2}^{a_{2}}\partial
x_{3}^{a_{3}}}dy_{1}dy_{2}dy_{3}d\tau_{1}\ ,$$where\
$k=a_{1}+a_{2}+a_{3}\ ,\ a_{1}\ ,\ a_{2}\ ,\ a_{3}$\ are all
nonnegative integral numbers, we know\
$h_{2}=h_{2}I_{K_{1}^{\prime}}\in C^{\infty}(K_{1})$\ . And because
$$\cfrac{\partial h_{2}}{\partial t}=\mu\triangle h_{2}-v\ ,$$ we know\
$h_{2}=h_{2}I_{K_{1}^{\prime}}\in C^{1}[0,\ T]\bigcap
C^{\infty}(K_{1})$\ .\\We know\ $h_{1}$\ satisfies the H\"{o}lder
condition if\ $h_{1}=h_{1}I_{K_{1}^{\prime}}\in
C^{1}(K_{1}^{\prime})$\ , hence we can get the corollary as follows.
\begin{corollary} \label{corollary1}$\forall\
h_{1}=h_{1}I_{K_{1}^{\prime}}\in C^{1}(K_{1}^{\prime})$\ , \
$\exists\ h_{2}=h_{2}I_{K_{1}^{\prime}}\in C^{2}[0,\ T]\bigcap
C^{\infty}(K_{1})$\ , such that\\ $F(h_{1})=a^{2}b\tau F(h_{2})$,
where$$a=\cfrac{\mu[(i\xi_{1})^{2}+(i\xi_{2})^{2}+(i\xi_{3})^{2}]-i\xi_{0}}{\tau}\
,\ b=\cfrac{(i\xi_{1})^{2}+(i\xi_{2})^{2}+(i\xi_{3})^{2}}{a}\
.$$\end{corollary} From the corollary 3.2 and the assumption 1.1 ,
we can get\ $\exists\ F_{j1}=F_{j1}I_{K_{1}^{\prime}}\in C^{2}[0,\
T]\bigcap C^{\infty}(K_{1})$\ , such that\ $F(F_{j})=a^{2}b\tau
F(F_{j1})$\ ,\ $1\leq j\leq 3$\ . And we can get the
following.\begin{eqnarray*}&&F_{j1}(t,M)=(\frac{1}{2\sqrt{\pi\mu}})^{3}\int_{0}^{t}\int_{R^{3}}\frac{F_{jv}(\tau_{1},
y_{1}, y_{2} , y_{3})}{(\sqrt{t-\tau_{1}})^{3}}\
e^{-\frac{(x_{1}-y_{1})^{2}+(x_{2}-y_{2})^{2}+(x_{3}-y_{3})^{2}}{4\mu(t-\tau_{1})}}
dy_{1}dy_{2}dy_{3}d\tau_{1}\
,\\&&F_{jv}(t,M_{0})=-\frac{1}{4\pi}\int_{K_{1}}\frac{F_{j}(t,M)}{r_{MM_{0}}}
dx_{1}dx_{2}dx_{3}\ ,\ 1\leq\ j\leq 3\ .\end{eqnarray*}Next we see\
$H[F^{-1}(Y_{1})]$\ , where\begin{eqnarray*}&&H=(-\alpha_{1},\
-\alpha_{2},\ -\alpha_{3},\ -\alpha_{4},\ e_{1},\ e_{2},\ e_{3},\
e_{4},\ e_{5},\ e_{6},\ e_{7},\ e_{8},\ e_{9},\ e_{10},\ e_{11},\
e_{16})^{T}\ ,\\&&\alpha_{1}=(0_{4},\ 1,\ 0_{4},\  1,\  0,\ 0,\  0,\
0,\ 0,\ 0,\ 0_{39})^{T},\\&& \alpha_{2}=(0_{12},\  \tau,\  0_{4},\
-\mu,\ -\mu,\ -\mu,\ 0_{26},\  1,\  1,\  1,\  0_{6} )^{T},\\&&
\alpha_{3}=(0_{13},\ \tau,\   0_{13},\   -\mu,\   -\mu,\   -\mu,\
0_{19},\   1,\   1,\ 1,\  0_{3} )^{T},\\&&  \alpha_{4}=(0_{14},\
\tau,\   0_{22},\ -\mu,\   -\mu,\   -\mu,\   0_{11},\   0,\   1,\
1,\   1)^{T}\ ,\end{eqnarray*} and\ $e_{i}$\ is the\ $ith$\ $55$\
dimensional unit coordinate vector,\ $1 \leq i \leq 55$\
,\begin{eqnarray*}Y_{1}&=&(i\xi_{2}y_{3},\ i\xi_{3}y_{3},\ y_{3},\
i\xi_{1}y_{7},\ i\xi_{2}y_{7},\ i\xi_{3}y_{7},\ y_{7},\
i\xi_{1}y_{11},\ i\xi_{2}y_{11},\ i\xi_{3}y_{11},\ y_{11},\\&&
i\xi_{0}y_{16},\ i\xi_{1}y_{16},\ i\xi_{2}y_{16},\ i\xi_{3}y_{16},\
y_{16},\ (i\xi_{0})^{2}y_{3},\ (i\xi_{1})^{2}y_{3},\
(i\xi_{2})^{2}y_{3},\ (i\xi_{3})^{2}y_{3},\
i\xi_{0}i\xi_{1}y_{3},\\&& i\xi_{0}i\xi_{2}y_{3},\
i\xi_{0}i\xi_{3}y_{3},\ i\xi_{1}i\xi_{2}y_{3},\
i\xi_{1}i\xi_{3}y_{3},\ i\xi_{2}i\xi_{3}y_{3},\
(i\xi_{0})^{2}y_{7},\ (i\xi_{1})^{2}y_{7},\ (i\xi_{2})^{2}y_{7},\
(i\xi_{3})^{2}y_{7},\\&& i\xi_{0}i\xi_{1}y_{7},\
i\xi_{0}i\xi_{2}y_{7},\ i\xi_{0}i\xi_{3}y_{7},\
i\xi_{1}i\xi_{2}y_{7},\ i\xi_{1}i\xi_{3}y_{7},\
i\xi_{2}i\xi_{3}y_{7},\ (i\xi_{0})^{2}y_{11},\
(i\xi_{1})^{2}y_{11},\ (i\xi_{2})^{2}y_{11},\\&&
(i\xi_{3})^{2}y_{11},\ i\xi_{0}i\xi_{1}y_{11},\
i\xi_{0}i\xi_{2}y_{11},\ i\xi_{0}i\xi_{3}y_{11},\
i\xi_{1}i\xi_{2}y_{11},\ i\xi_{1}i\xi_{3}y_{11},\
i\xi_{2}i\xi_{3}y_{11},\ 0_{9\times1}^{T})^{T}\
,\end{eqnarray*}and\begin{eqnarray*}
y_{3}&=&\cfrac{i\xi_{1}(i\xi_{1}F(F_{1})+i\xi_{2}F(F_{2})+i\xi_{3}F(F_{3}))}{a^{2}b\tau}-\cfrac{F(F_{1})}{a\tau}
\\&=&i\xi_{1}(i\xi_{1}F(F_{11})+i\xi_{2}F(F_{21})+i\xi_{3}F(F_{31}))-a
b F(F_{11})\
,\\y_{7}&=&\cfrac{i\xi_{2}(i\xi_{1}F(F_{1})+i\xi_{2}F(F_{2})+i\xi_{3}F(F_{3}))}{a^{2}b\tau}-\cfrac{F(F_{2})}{a\tau}
\\&=&i\xi_{2}(i\xi_{1}F(F_{11})+i\xi_{2}F(F_{21})+i\xi_{3}F(F_{31}))-a
b F(F_{21})\ ,
\\y_{11}&=&\cfrac{i\xi_{3}(i\xi_{1}F(F_{1})+i\xi_{2}F(F_{2})+i\xi_{3}F(F_{3}))}{a^{2}b\tau}-\cfrac{F(F_{3})}{a\tau}
\\&=&i\xi_{3}(i\xi_{1}F(F_{11})+i\xi_{2}F(F_{21})+i\xi_{3}F(F_{31}))-a
b F(F_{31})\
,\\y_{16}&=&\cfrac{i\xi_{1}F(F_{1})+i\xi_{2}F(F_{2})+i\xi_{3}F(F_{3})}{ab\tau}
\\&=&a(i\xi_{1}F(F_{11})+i\xi_{2}F(F_{21})+i\xi_{3}F(F_{31}))\ .\end{eqnarray*}
We assume\ $H[F^{-1}(Y_{1})]=W_{1}(F_{11},\ F_{21},\ F_{31})$\ ,
where\ $W_{1}(F_{11},\ F_{21},\ F_{31})=(w_{1j})_{16\times 1}$\ ,
and we can get the following.\begin{eqnarray*}&&
w_{11}=-(e_{5}^{T}+e_{10}^{T})F^{-1}(Y_{1})=-F^{-1}(i\xi_{2}y_{7}+i\xi_{3}y_{11})\
,\\&&w_{12}=-[\tau
e_{13}^{T}-\mu(e_{18}^{T}+e_{19}^{T}+e_{20}^{T})+e_{47}^{T}+e_{48}^{T}+e_{49}^{T}]F^{-1}(Y_{1})
\\&&=-F^{-1}\{\tau
i\xi_{1}y_{16}-\mu[(i\xi_{1})^{2}+(i\xi_{2})^{2}+(i\xi_{3})^{2}]y_{3}\}\
,\\&&w_{13}=-[\tau
e_{14}^{T}-\mu(e_{28}^{T}+e_{29}^{T}+e_{30}^{T})+e_{50}^{T}+e_{51}^{T}+e_{52}^{T}]F^{-1}(Y_{1})
\\&&=-F^{-1}\{\tau
i\xi_{2}y_{16}-\mu[(i\xi_{1})^{2}+(i\xi_{2})^{2}+(i\xi_{3})^{2}]y_{7}\}\
,\\&&w_{14}=-[\tau
e_{15}^{T}-\mu(e_{38}^{T}+e_{39}^{T}+e_{40}^{T})+e_{53}^{T}+e_{54}^{T}+e_{55}^{T}]F^{-1}(Y_{1})
\\&&=-F^{-1}\{\tau
i\xi_{1}y_{16}-\mu[(i\xi_{1})^{2}+(i\xi_{2})^{2}+(i\xi_{3})^{2}]y_{11}\}\
,\\&&w_{15}=e_{1}^{T}F^{-1}(Y_{1})=F^{-1}(i\xi_{2}y_{3})\ ,\
w_{16}=e_{2}^{T}F^{-1}(Y_{1})=F^{-1}(i\xi_{3}y_{3})\
,\\&&w_{17}=e_{3}^{T}F^{-1}(Y_{1})=F^{-1}(y_{3})\ ,\
w_{18}=e_{4}^{T}F^{-1}(Y_{1})=F^{-1}(i\xi_{1}y_{7})\
,\\&&w_{19}=e_{5}^{T}F^{-1}(Y_{1})=F^{-1}(i\xi_{2}y_{7})\ ,\
w_{110}=e_{6}^{T}F^{-1}(Y_{1})=F^{-1}(i\xi_{3}y_{7})\
,\\&&w_{111}=e_{7}^{T}F^{-1}(Y_{1})=F^{-1}(y_{7})\ ,\
w_{112}=e_{8}^{T}F^{-1}(Y_{1})=F^{-1}(i\xi_{1}y_{11})\
,\\&&w_{113}=e_{9}^{T}F^{-1}(Y_{1})=F^{-1}(i\xi_{2}y_{11})\ ,\
w_{114}=e_{10}^{T}F^{-1}(Y_{1})=F^{-1}(i\xi_{3}y_{11})\
,\\&&w_{115}=e_{11}^{T}F^{-1}(Y_{1})=F^{-1}(y_{11})\ ,\
w_{116}=e_{16}^{T}F^{-1}(Y_{1})=F^{-1}(y_{16})\ .\end{eqnarray*}
Because\ $F_{j1}=F_{j1}I_{K_{1}^{\prime}}\in C^{2}[0,\ T]\bigcap
C^{\infty}(K_{1})$\ ,\ $1\leq j\leq 3$\ , we can get as
follows,$$F^{-1}(i\xi_{0}F(F_{j1}))=F^{-1}(F(\cfrac{\partial
F_{j1}}{\partial t}))=\cfrac{\partial F_{j1}}{\partial t}\ ,\ 1\leq
j\leq 3\ ,$$for the same reason we can get the
following,$$F^{-1}(i\xi_{1}F(F_{j1}))=\cfrac{\partial
F_{j1}}{\partial x_{1}}\ ,\
F^{-1}(i\xi_{2}F(F_{j1}))=\cfrac{\partial F_{j1}}{\partial x_{2}}\
,\ F^{-1}(i\xi_{3}F(F_{j1}))=\cfrac{\partial F_{j1}}{\partial
x_{3}}\ ,\ 1\leq j\leq 3\ .$$Hence we can get the
following,\begin{eqnarray*}F^{-1}(y_{3})&=&
\cfrac{\partial^{2}F_{21}}{\partial x_{1}\partial
x_{2}}+\cfrac{\partial^{2}F_{31}}{\partial x_{1}\partial
x_{3}}-\cfrac{\partial^{2}F_{11}}{\partial
x_{2}^{2}}-\cfrac{\partial^{2}F_{11}}{\partial x_{3}^{2}}\
,\\F^{-1}(y_{7})&=& \cfrac{\partial^{2}F_{11}}{\partial
x_{1}\partial x_{2}}+\cfrac{\partial^{2}F_{31}}{\partial
x_{2}\partial x_{3}}-\cfrac{\partial^{2}F_{21}}{\partial
x_{1}^{2}}-\cfrac{\partial^{2}F_{21}}{\partial x_{3}^{2}}\
,\\F^{-1}(y_{11})&=& \cfrac{\partial^{2}F_{11}}{\partial
x_{1}\partial x_{3}}+\cfrac{\partial^{2}F_{21}}{\partial
x_{2}\partial x_{3}}-\cfrac{\partial^{2}F_{31}}{\partial
x_{1}^{2}}-\cfrac{\partial^{2}F_{31}}{\partial x_{2}^{2}}\
,\\F^{-1}(y_{16})&=&\cfrac{1}{\tau}\ [\ \mu
\triangle(\cfrac{\partial F_{11}}{\partial x_{1}}+\cfrac{\partial
F_{21}}{\partial x_{2}}+\cfrac{\partial F_{31}}{\partial
x_{3}})-\cfrac{\partial^{2} F_{11}}{\partial t\partial
x_{1}}-\cfrac{\partial^{2} F_{21}}{\partial t\partial
x_{2}}-\cfrac{\partial^{2} F_{31}}{\partial t\partial x_{3}}\ ]\ .
\end{eqnarray*}And we can get the
following,\begin{eqnarray*}&&w_{11}=-\cfrac{\partial
F^{-1}(y_{7})}{\partial x_{2}}-\cfrac{\partial
F^{-1}(y_{11})}{\partial x_{3}}=\cfrac{\partial^{3}F_{21}}{\partial
x_{1}^{2}\partial x_{2}}+\cfrac{\partial^{3}F_{31}}{\partial
x_{1}^{2}\partial x_{3}}-\cfrac{\partial^{3}F_{11}}{\partial
x_{1}\partial x_{2}^{2}}-\cfrac{\partial^{3}F_{11}}{\partial
x_{1}\partial x_{3}^{2}}\ ,\\&&w_{12}=-\cfrac{\partial \tau
F^{-1}(y_{16})}{\partial x_{1}}+\mu \triangle F^{-1}(y_{3})\ ,\
w_{13}=-\cfrac{\partial \tau F^{-1}(y_{16})}{\partial x_{2}}+\mu
\triangle F^{-1}(y_{7})\ ,\\&& w_{14}=-\cfrac{\partial \tau
F^{-1}(y_{16})}{\partial x_{3}}+\mu \triangle F^{-1}(y_{11})\ ,\\&&
w_{15}=\cfrac{\partial F^{-1}(y_{3})}{\partial
x_{2}}=\cfrac{\partial^{3}F_{21}}{\partial x_{1}\partial
x_{2}^{2}}+\cfrac{\partial^{3}F_{31}}{\partial x_{1}\partial
x_{2}\partial x_{3}}-\cfrac{\partial^{3}F_{11}}{\partial
x_{2}^{3}}-\cfrac{\partial^{3}F_{11}}{\partial x_{2}\partial
x_{3}^{2}}\ ,\\&& w_{16}=\cfrac{\partial F^{-1}(y_{3})}{\partial
x_{3}}=\cfrac{\partial^{3}F_{21}}{\partial x_{1}\partial
x_{2}\partial x_{3}}+\cfrac{\partial^{3}F_{31}}{\partial
x_{1}\partial x_{3}^{2}}-\cfrac{\partial^{3}F_{11}}{\partial
x_{2}^{2}\partial x_{3}}-\cfrac{\partial^{3}F_{11}}{ \partial
x_{3}^{3}}\ ,\\&&
w_{17}=F^{-1}(y_{3})=\cfrac{\partial^{2}F_{21}}{\partial
x_{1}\partial x_{2}}+\cfrac{\partial^{2}F_{31}}{\partial
x_{1}\partial x_{3}}-\cfrac{\partial^{2}F_{11}}{\partial
x_{2}^{2}}-\cfrac{\partial^{2}F_{11}}{\partial x_{3}^{2}}\ ,\\&&
w_{18}=\cfrac{\partial F^{-1}(y_{7})}{\partial
x_{1}}=\cfrac{\partial^{3}F_{11}}{\partial x_{1}^{2}\partial
x_{2}}+\cfrac{\partial^{3}F_{31}}{\partial x_{1}\partial
x_{2}\partial x_{3}}-\cfrac{\partial^{3}F_{21}}{\partial
x_{1}^{3}}-\cfrac{\partial^{3}F_{21}}{\partial x_{1}\partial
x_{3}^{2}}\ ,\\&& w_{19}=\cfrac{\partial F^{-1}(y_{7})}{\partial
x_{2}}=\cfrac{\partial^{3}F_{11}}{\partial x_{1}\partial
x_{2}^{2}}+\cfrac{\partial^{3}F_{31}}{\partial x_{2}^{2}\partial
x_{3}}-\cfrac{\partial^{3}F_{21}}{\partial x_{1}^{2}\partial
x_{2}}-\cfrac{\partial^{3}F_{21}}{\partial x_{2}\partial x_{3}^{2}}\
,\\&& w_{110}=\cfrac{\partial F^{-1}(y_{7})}{\partial
x_{3}}=\cfrac{\partial^{3}F_{11}}{\partial x_{1}\partial
x_{2}\partial x_{3}}+\cfrac{\partial^{3}F_{31}}{ \partial
x_{2}\partial x_{3}^{2}}-\cfrac{\partial^{3}F_{21}}{\partial
x_{1}^{2}\partial x_{3}}-\cfrac{\partial^{3}F_{21}}{ \partial
x_{3}^{3}}\ ,\\&&w_{111}=F^{-1}(y_{7})=
\cfrac{\partial^{2}F_{11}}{\partial x_{1}\partial
x_{2}}+\cfrac{\partial^{2}F_{31}}{\partial x_{2}\partial
x_{3}}-\cfrac{\partial^{2}F_{21}}{\partial
x_{1}^{2}}-\cfrac{\partial^{2}F_{21}}{\partial x_{3}^{2}}\ ,\\&&
w_{112}=\cfrac{\partial F^{-1}(y_{11})}{\partial
x_{1}}=\cfrac{\partial^{3}F_{11}}{\partial x_{1}^{2}\partial
x_{3}}+\cfrac{\partial^{3}F_{21}}{\partial x_{1}\partial
x_{2}\partial x_{3}}-\cfrac{\partial^{3}F_{31}}{\partial
x_{1}^{3}}-\cfrac{\partial^{3}F_{31}}{\partial x_{1}\partial
x_{2}^{2}}\ ,\\&& w_{113}=\cfrac{\partial F^{-1}(y_{11})}{\partial
x_{2}}=\cfrac{\partial^{3}F_{11}}{\partial x_{1}\partial
x_{2}\partial x_{3}}+\cfrac{\partial^{3}F_{21}}{ \partial
x_{2}^{2}\partial x_{3}}-\cfrac{\partial^{3}F_{31}}{\partial
x_{1}^{2}\partial x_{2}}-\cfrac{\partial^{3}F_{31}}{ \partial
x_{2}^{3}}\ ,\\&& w_{114}=\cfrac{\partial F^{-1}(y_{11})}{\partial
x_{3}}=\cfrac{\partial^{3}F_{11}}{\partial x_{1} \partial
x_{3}^{2}}+\cfrac{\partial^{3}F_{21}}{ \partial x_{2}\partial
x_{3}^{2}}-\cfrac{\partial^{3}F_{31}}{\partial x_{1}^{2}\partial
x_{3}}-\cfrac{\partial^{3}F_{31}}{ \partial x_{2}^{2}\partial
x_{3}}\
,\end{eqnarray*}\begin{eqnarray*}&&w_{115}=F^{-1}(y_{11})=\cfrac{\partial^{2}F_{11}}{\partial
x_{1}\partial x_{3}}+\cfrac{\partial^{2}F_{21}}{\partial
x_{2}\partial x_{3}}-\cfrac{\partial^{2}F_{31}}{\partial
x_{1}^{2}}-\cfrac{\partial^{2}F_{31}}{\partial x_{2}^{2}}\
,\\&&w_{116}=F^{-1}(y_{16})=\cfrac{1}{\tau}\ [\ \mu
\triangle(\cfrac{\partial F_{11}}{\partial x_{1}}+\cfrac{\partial
F_{21}}{\partial x_{2}}+\cfrac{\partial F_{31}}{\partial
x_{3}})-\cfrac{\partial^{2} F_{11}}{\partial t\partial
x_{1}}-\cfrac{\partial^{2} F_{21}}{\partial t\partial
x_{2}}-\cfrac{\partial^{2} F_{31}}{\partial t\partial x_{3}}\ ]\
.\end{eqnarray*} And we can see that\ $W_{1}$\ is the function to do
the partial derivation with\ $F_{11},\ F_{21},\ F_{31}$\ no more
than the fourth order and their linear combination, moreover no more
than the first order with the variable\ $t$\ . Hence\
$H[F^{-1}(Y_{1})]=H[F^{-1}(Y_{1})]I_{K_{1}^{\prime}}\in
C^{1}(K_{1}^{\prime})$\ .\\ If\ $Z_{1}\in \Omega_{2}$\ , then there
exists\ $h_{1}=h_{1}I_{K_{1}^{\prime}}\in C^{1}(K_{1}^{\prime})$\ ,
such that\ $Z_{1}=F(h_{1})$\ . Again from the corollary 3.2 , we can
get\ $\exists\ h_{2}=h_{2}I_{K_{1}^{\prime}}\in C^{2}[0,\ T]\bigcap
C^{\infty}(K_{1})$\ , such that\ $F(h_{1})=a^{2}b\tau F(h_{2})$\ .
Next we see\ $H[F^{-1}(A_{1}(\eta)Z_{1})]$\ , where\
$A_{1}(\eta)=(\eta_{1} ,\ \eta_{2} ,\ \eta_{3} ,\ \eta_{4} ,\
\eta_{5},\ \eta_{6} ,\ \eta_{7} ,\ \eta_{8} ,\ \eta_{9})$\
,\begin{eqnarray*}\eta_{j}&=&(i\xi_{2}y_{3},\ i\xi_{3}y_{3},\
y_{3},\ i\xi_{1}y_{7},\ i\xi_{2}y_{7},\ i\xi_{3}y_{7},\ y_{7},\
i\xi_{1}y_{11},\ i\xi_{2}y_{11},\ i\xi_{3}y_{11},\ y_{11},\\&&
i\xi_{0}y_{16},\ i\xi_{1}y_{16},\ i\xi_{2}y_{16},\ i\xi_{3}y_{16},\
y_{16},\ (i\xi_{0})^{2}y_{3},\ (i\xi_{1})^{2}y_{3},\
(i\xi_{2})^{2}y_{3},\ (i\xi_{3})^{2}y_{3},\
i\xi_{0}i\xi_{1}y_{3},\\&& i\xi_{0}i\xi_{2}y_{3},\
i\xi_{0}i\xi_{3}y_{3},\ i\xi_{1}i\xi_{2}y_{3},\
i\xi_{1}i\xi_{3}y_{3},\ i\xi_{2}i\xi_{3}y_{3},\
(i\xi_{0})^{2}y_{7},\ (i\xi_{1})^{2}y_{7},\ (i\xi_{2})^{2}y_{7},\
(i\xi_{3})^{2}y_{7},\\&& i\xi_{0}i\xi_{1}y_{7},\
i\xi_{0}i\xi_{2}y_{7},\ i\xi_{0}i\xi_{3}y_{7},\
i\xi_{1}i\xi_{2}y_{7},\ i\xi_{1}i\xi_{3}y_{7},\
i\xi_{2}i\xi_{3}y_{7},\ (i\xi_{0})^{2}y_{11},\
(i\xi_{1})^{2}y_{11},\ (i\xi_{2})^{2}y_{11},\\&&
(i\xi_{3})^{2}y_{11},\ i\xi_{0}i\xi_{1}y_{11},\
i\xi_{0}i\xi_{2}y_{11},\ i\xi_{0}i\xi_{3}y_{11},\
i\xi_{1}i\xi_{2}y_{11},\ i\xi_{1}i\xi_{3}y_{11},\
i\xi_{2}i\xi_{3}y_{11},\ e_{j}^{T})^{T}\ ,\end{eqnarray*}here\
$e_{j}$\ is the\ $jth$\ $9$\ dimensional unit coordinate vector,\ $1
\leq j \leq 9$\ , moreover \\when\ $j=1,\ 2,\ 3$\
,$$y_{3}=\cfrac{1}{a^{2}b\tau}(\xi_{1}^{2}+ab)\ ,\
y_{7}=\cfrac{1}{a^{2}b\tau}(-i\xi_{1}i\xi_{2})\ ,\
y_{11}=\cfrac{1}{a^{2}b\tau}(-i\xi_{1}i\xi_{3})\ ,\
y_{16}=\cfrac{1}{a^{2}b\tau}(-i\xi_{1}a)\ ,$$where\
$b=\cfrac{(i\xi_{1})^{2}+(i\xi_{2})^{2}+(i\xi_{3})^{2}}{a}\neq 0$\ ,
and\ $(\xi_{1},\ \xi_{2},\ \xi_{3})\neq (0\ ,\ 0\ ,\ 0)$\ ,\\ when\
$j=4,\ 5,\ 6$\ ,$$y_{7}=\cfrac{1}{a^{2}b\tau}(\xi_{2}^{2}+ab)\ ,\
y_{3}=\cfrac{1}{a^{2}b\tau}(-i\xi_{1}i\xi_{2})\ ,\
y_{11}=\cfrac{1}{a^{2}b\tau}(-i\xi_{2}i\xi_{3})\ ,\
y_{16}=\cfrac{1}{a^{2}b\tau}(-i\xi_{2}a)\ ,$$when\ $j=7,\ 8,\ 9$\
,$$y_{11}=\cfrac{1}{a^{2}b\tau}(\xi_{3}^{2}+ab)\ ,\
y_{3}=\cfrac{1}{a^{2}b\tau}(-i\xi_{1}i\xi_{3})\ ,\
y_{7}=\cfrac{1}{a^{2}b\tau}(-i\xi_{2}i\xi_{3})\ ,\
y_{16}=\cfrac{1}{a^{2}b\tau}(-i\xi_{3}a)\ .$$ Hence we can
get\begin{eqnarray*}&&H[F^{-1}(A_{1}(\eta)Z_{1})]=H[F^{-1}(A_{1}(\eta)F(h_{1}))]\\&&=H[F^{-1}(A_{1}(\eta)a^{2}b\tau
F(h_{2}))]=H\{F^{-1}[(a^{2}b\tau A_{1}(\eta))F(h_{2})]\}\
,\end{eqnarray*} and we can see that\ $a^{2}b\tau A_{1}(\eta)$\ is a
polynomial matrix. We assume$$H\{F^{-1}[(a^{2}b\tau
A_{1}(\eta))F(h_{2})]\}=W_{2}(h_{2})\ ,$$ where\
$W_{2}(h_{2})=(w_{2j})_{16\times 1}$\ , and we assume\
$h_{2}=(h_{2j})_{9\times 1}$\ , then we can get the
following.\begin{eqnarray*}
w_{21}&=&-(e_{5}^{T}+e_{10}^{T})F^{-1}[(a^{2}b\tau
A_{1}(\eta))F(h_{2})]
\\&=&-F^{-1}\{i\xi_{2}[-i\xi_{1}i\xi_{2}(F(h_{21})+F(h_{22})+F(h_{23}))+(\xi_{2}^{2}+a
b)(F(h_{24})+F(h_{25})+F(h_{26}))\\&&-i\xi_{1}i\xi_{2}(F(h_{27})+F(h_{28})+F(h_{29}))]
+i\xi_{3}[-i\xi_{1}i\xi_{3}(F(h_{21})+F(h_{22})+F(h_{23}))
\\&&-i\xi_{2}i\xi_{3}(F(h_{24})+F(h_{25})+F(h_{26}))+(\xi_{3}^{2}+a
b)(F(h_{27})+F(h_{28})+F(h_{29}))]\}\
,\end{eqnarray*}\begin{eqnarray*}w_{22}&=&-[\tau
e_{13}^{T}-\mu(e_{18}^{T}+e_{19}^{T}+e_{20}^{T})+e_{47}^{T}+e_{48}^{T}+e_{49}^{T}]F^{-1}[(a^{2}b\tau
A_{1}(\eta))F(h_{2})]
 \\&=&-F^{-1}\{[\tau
i\xi_{1}(-i\xi_{1}a)-\mu((i\xi_{1})^{2}+(i\xi_{2})^{2}+(i\xi_{3})^{2})(\xi_{1}^{2}+a
b)+a^{2}b\tau](F(h_{21})+F(h_{22})+F(h_{23}))
\\&&+[\tau
i\xi_{1}(-i\xi_{2}a)-\mu((i\xi_{1})^{2}+(i\xi_{2})^{2}+(i\xi_{3})^{2})(-i\xi_{1}i\xi_{2})](F(h_{24})+F(h_{25})+F(h_{26}))
\\&&+[\tau
i\xi_{1}(-i\xi_{3}a)-\mu((i\xi_{1})^{2}+(i\xi_{2})^{2}+(i\xi_{3})^{2})(-i\xi_{1}i\xi_{3})](F(h_{27})+F(h_{28})+F(h_{29}))\}\
,\\w_{23}&=&-[\tau
e_{14}^{T}-\mu(e_{28}^{T}+e_{29}^{T}+e_{30}^{T})+e_{50}^{T}+e_{51}^{T}+e_{52}^{T}]F^{-1}[(a^{2}b\tau
A_{1}(\eta))F(h_{2})]
 \\&=&-F^{-1}\{[\tau
i\xi_{2}(-i\xi_{1}a)-\mu((i\xi_{1})^{2}+(i\xi_{2})^{2}+(i\xi_{3})^{2})(-i\xi_{1}i\xi_{2})](F(h_{21})+F(h_{22})+F(h_{23}))
\\&&+[\tau
i\xi_{2}(-i\xi_{2}a)-\mu((i\xi_{1})^{2}+(i\xi_{2})^{2}+(i\xi_{3})^{2})(\xi_{2}^{2}+a
b)+a^{2}b\tau](F(h_{24})+F(h_{25})+F(h_{26}))
\\&&+[\tau
i\xi_{2}(-i\xi_{3}a)-\mu((i\xi_{1})^{2}+(i\xi_{2})^{2}+(i\xi_{3})^{2})(-i\xi_{2}i\xi_{3})](F(h_{27})+F(h_{28})+F(h_{29}))\}\
,\\w_{24}&=&-[\tau
e_{15}^{T}-\mu(e_{38}^{T}+e_{39}^{T}+e_{40}^{T})+e_{53}^{T}+e_{54}^{T}+e_{55}^{T}]F^{-1}[(a^{2}b\tau
A_{1}(\eta))F(h_{2})]
 \\&=&-F^{-1}\{[\tau
i\xi_{3}(-i\xi_{1}a)-\mu((i\xi_{1})^{2}+(i\xi_{2})^{2}+(i\xi_{3})^{2})(-i\xi_{1}i\xi_{3})](F(h_{21})+F(h_{22})+F(h_{23}))
\\&&+[\tau
i\xi_{3}(-i\xi_{2}a)-\mu((i\xi_{1})^{2}+(i\xi_{2})^{2}+(i\xi_{3})^{2})(-i\xi_{2}i\xi_{3})](F(h_{24})+F(h_{25})+F(h_{26}))
\\&&+[\tau
i\xi_{3}(-i\xi_{3}a)-\mu((i\xi_{1})^{2}+(i\xi_{2})^{2}+(i\xi_{3})^{2})(\xi_{3}^{2}+a
b)+a^{2}b\tau](F(h_{27})+F(h_{28})+F(h_{29}))\}\ ,
\\w_{25}&=&e_{1}^{T}F^{-1}[(a^{2}b\tau
A_{1}(\eta))F(h_{2})]=F^{-1}\{i\xi_{2}[(\xi_{1}^{2}+a
b)(F(h_{21})+F(h_{22})+F(h_{23}))\\&&-i\xi_{1}i\xi_{2}(F(h_{24})+F(h_{25})+F(h_{26}))-i\xi_{1}i\xi_{3}(F(h_{27})+F(h_{28})+F(h_{29}))]\}\
,\\w_{26}&=&e_{2}^{T}F^{-1}[(a^{2}b\tau
A_{1}(\eta))F(h_{2})]=F^{-1}\{i\xi_{3}[(\xi_{1}^{2}+a
b)(F(h_{21})+F(h_{22})+F(h_{23}))\\&&-i\xi_{1}i\xi_{2}(F(h_{24})+F(h_{25})+F(h_{26}))-i\xi_{1}i\xi_{3}(F(h_{27})+F(h_{28})+F(h_{29}))]\}\
,\\w_{27}&=&e_{3}^{T}F^{-1}[(a^{2}b\tau
A_{1}(\eta))F(h_{2})]=F^{-1}\{(\xi_{1}^{2}+a
b)(F(h_{21})+F(h_{22})+F(h_{23}))\\&&-i\xi_{1}i\xi_{2}(F(h_{24})+F(h_{25})+F(h_{26}))-i\xi_{1}i\xi_{3}(F(h_{27})+F(h_{28})+F(h_{29}))\}\
,\\w_{28}&=&e_{4}^{T}F^{-1}[(a^{2}b\tau
A_{1}(\eta))F(h_{2})]=F^{-1}\{i\xi_{1}[-i\xi_{1}i\xi_{2}(F(h_{21})+F(h_{22})+F(h_{23}))\\&&+(\xi_{2}^{2}+a
b)(F(h_{24})+F(h_{25})+F(h_{26}))-i\xi_{2}i\xi_{3}(F(h_{27})+F(h_{28})+F(h_{29}))]\}\
,\\w_{29}&=&e_{5}^{T}F^{-1}[(a^{2}b\tau
A_{1}(\eta))F(h_{2})]=F^{-1}\{i\xi_{2}[-i\xi_{1}i\xi_{2}(F(h_{21})+F(h_{22})+F(h_{23}))\\&&+(\xi_{2}^{2}+a
b)(F(h_{24})+F(h_{25})+F(h_{26}))-i\xi_{2}i\xi_{3}(F(h_{27})+F(h_{28})+F(h_{29}))]\}\
,\\w_{210}&=&e_{6}^{T}F^{-1}[(a^{2}b\tau
A_{1}(\eta))F(h_{2})]=F^{-1}\{i\xi_{3}[-i\xi_{1}i\xi_{2}(F(h_{21})+F(h_{22})+F(h_{23}))\\&&+(\xi_{2}^{2}+a
b)(F(h_{24})+F(h_{25})+F(h_{26}))-i\xi_{2}i\xi_{3}(F(h_{27})+F(h_{28})+F(h_{29}))]\}\
,\\w_{211}&=&e_{7}^{T}F^{-1}[(a^{2}b\tau
A_{1}(\eta))F(h_{2})]=F^{-1}\{-i\xi_{1}i\xi_{2}(F(h_{21})+F(h_{22})+F(h_{23}))\\&&+(\xi_{2}^{2}+a
b)(F(h_{24})+F(h_{25})+F(h_{26}))-i\xi_{2}i\xi_{3}(F(h_{27})+F(h_{28})+F(h_{29}))\}\
,\\w_{212}&=&e_{8}^{T}F^{-1}[(a^{2}b\tau
A_{1}(\eta))F(h_{2})]=F^{-1}\{i\xi_{1}[-i\xi_{1}i\xi_{3}(F(h_{21})+F(h_{22})+F(h_{23}))\\&&-i\xi_{2}i\xi_{3}(F(h_{24})+F(h_{25})+F(h_{26}))+(\xi_{3}^{2}+a
b)(F(h_{27})+F(h_{28})+F(h_{29}))]\}\
,\\w_{213}&=&e_{9}^{T}F^{-1}[(a^{2}b\tau
A_{1}(\eta))F(h_{2})]=F^{-1}\{i\xi_{2}[-i\xi_{1}i\xi_{3}(F(h_{21})+F(h_{22})+F(h_{23}))\\&&-i\xi_{2}i\xi_{3}(F(h_{24})+F(h_{25})+F(h_{26}))+(\xi_{3}^{2}+a
b)(F(h_{27})+F(h_{28})+F(h_{29}))]\}\
,\\w_{214}&=&e_{10}^{T}F^{-1}[(a^{2}b\tau
A_{1}(\eta))F(h_{2})]=F^{-1}\{i\xi_{3}[-i\xi_{1}i\xi_{3}(F(h_{21})+F(h_{22})+F(h_{23}))\\&&-i\xi_{2}i\xi_{3}(F(h_{24})+F(h_{25})+F(h_{26}))+(\xi_{3}^{2}+a
b)(F(h_{27})+F(h_{28})+F(h_{29}))]\}\
,\\w_{215}&=&e_{11}^{T}F^{-1}[(a^{2}b\tau
A_{1}(\eta))F(h_{2})]=F^{-1}\{-i\xi_{1}i\xi_{3}(F(h_{21})+F(h_{22})+F(h_{23}))\\&&-i\xi_{2}i\xi_{3}(F(h_{24})+F(h_{25})+F(h_{26}))+(\xi_{3}^{2}+a
b)(F(h_{27})+F(h_{28})+F(h_{29}))\}\
,\\w_{216}&=&e_{16}^{T}F^{-1}[(a^{2}b\tau
A_{1}(\eta))F(h_{2})]=F^{-1}\{-i\xi_{1}a(F(h_{21})+F(h_{22})+F(h_{23}))\\&&-i\xi_{2}a(F(h_{24})+F(h_{25})+F(h_{26}))-i\xi_{3}a
(F(h_{27})+F(h_{28})+F(h_{29}))\}\ .
\end{eqnarray*}
Because\ $h_{2}=h_{2}I_{K_{1}^{\prime}}\in C^{2}[0,\ T]\bigcap
C^{\infty}(K_{1})$\ , we can get
$$F^{-1}(i\xi_{0}F(h_{2}))=F^{-1}(F(\cfrac{\partial h_{2}}{\partial
t}))=\cfrac{\partial h_{2}}{\partial t}\ ,$$ for the same reason we
can get$$F^{-1}(i\xi_{1}F(h_{2}))=\cfrac{\partial h_{2}}{\partial
x_{1}}\ ,\ F^{-1}(i\xi_{2}F(h_{2}))=\cfrac{\partial h_{2}}{\partial
x_{2}}\ ,\ F^{-1}(i\xi_{3}F(h_{2}))=\cfrac{\partial h_{2}}{\partial
x_{3}}\ .$$ If we assume\ $h_{31}=h_{21}+h_{22}+h_{23}\ ,\
h_{32}=h_{24}+h_{25}+h_{26}\ ,\ h_{33}=h_{27}+h_{28}+h_{29}$\ , then
we can get the following.
\begin{eqnarray*}&&w_{21}=\cfrac{\partial^{3}h_{31}}{\partial
x_{1}\partial x_{2}^{2}}+\cfrac{\partial^{3}h_{31}}{\partial
x_{1}\partial x_{3}^{2}}-\cfrac{\partial^{3}h_{32}}{\partial
x_{1}^{2}\partial x_{2}}-\cfrac{\partial^{3}h_{33}}{\partial
x_{1}^{2}\partial x_{3}}\ ,\
w_{22}=\cfrac{\partial^{3}h_{31}}{\partial t\partial
x_{2}^{2}}+\cfrac{\partial^{3}h_{31}}{\partial t\partial
x_{3}^{2}}-\cfrac{\partial^{3}h_{32}}{\partial t\partial
x_{1}\partial x_{2}}-\cfrac{\partial^{3}h_{33}}{\partial t\partial
x_{1}\partial x_{3}}\
,\\&&w_{23}=\cfrac{\partial^{3}h_{32}}{\partial t\partial
x_{1}^{2}}+\cfrac{\partial^{3}h_{32}}{\partial t\partial
x_{3}^{2}}-\cfrac{\partial^{3}h_{31}}{\partial t\partial
x_{1}\partial x_{2}}-\cfrac{\partial^{3}h_{33}}{\partial t\partial
x_{2}\partial x_{3}}\ ,\ w_{24}=\cfrac{\partial^{3}h_{33}}{\partial
t\partial x_{1}^{2}}+\cfrac{\partial^{3}h_{33}}{\partial t\partial
x_{2}^{2}}-\cfrac{\partial^{3}h_{31}}{\partial t\partial
x_{1}\partial x_{3}}-\cfrac{\partial^{3}h_{32}}{\partial t\partial
x_{2}\partial x_{3}}\
,\\&&w_{25}=\cfrac{\partial^{3}h_{31}}{\partial
x_{2}^{3}}+\cfrac{\partial^{3}h_{31}}{\partial x_{2}\partial
x_{3}^{2}}-\cfrac{\partial^{3}h_{32}}{\partial x_{1}\partial
x_{2}^{2}}-\cfrac{\partial^{3}h_{33}}{\partial x_{1}\partial
x_{2}\partial x_{3}}\ ,\ w_{26}=\cfrac{\partial^{3}h_{31}}{\partial
x_{2}^{2}\partial x_{3}}+\cfrac{\partial^{3}h_{31}}{\partial
x_{3}^{3}}-\cfrac{\partial^{3}h_{32}}{\partial x_{1}\partial
x_{2}\partial x_{3}}-\cfrac{\partial^{3}h_{33}}{\partial
x_{1}\partial x_{3}^{2}}\
,\\&&w_{27}=\cfrac{\partial^{2}h_{31}}{\partial
x_{2}^{2}}+\cfrac{\partial^{2}h_{31}}{\partial
x_{3}^{2}}-\cfrac{\partial^{2}h_{32}}{\partial x_{1}\partial
x_{2}}-\cfrac{\partial^{2}h_{33}}{\partial x_{1}\partial x_{3}}\ ,\
w_{28}=\cfrac{\partial^{3}h_{32}}{\partial
x_{1}^{3}}+\cfrac{\partial^{3}h_{32}}{\partial x_{1}\partial
x_{3}^{2}}-\cfrac{\partial^{3}h_{31}}{\partial x_{1}^{2}\partial
x_{2}}-\cfrac{\partial^{3}h_{33}}{\partial x_{1}\partial
x_{2}\partial x_{3}}\
,\\&&w_{29}=\cfrac{\partial^{3}h_{32}}{\partial x_{1}^{2}\partial
x_{2}}+\cfrac{\partial^{3}h_{32}}{\partial x_{2}\partial
x_{3}^{2}}-\cfrac{\partial^{3}h_{31}}{\partial x_{1}\partial
x_{2}^{2}}-\cfrac{\partial^{3}h_{33}}{\partial x_{2}^{2}\partial
x_{3}}\ ,\ w_{210}=\cfrac{\partial^{3}h_{32}}{\partial
x_{1}^{2}\partial x_{3}}+\cfrac{\partial^{3}h_{32}}{\partial
x_{3}^{3}}-\cfrac{\partial^{3}h_{31}}{\partial x_{1}\partial
x_{2}\partial x_{3}}-\cfrac{\partial^{3}h_{33}}{\partial
x_{2}\partial x_{3}^{2}}\
,\\&&w_{211}=\cfrac{\partial^{2}h_{32}}{\partial
x_{1}^{2}}+\cfrac{\partial^{2}h_{32}}{\partial
x_{3}^{2}}-\cfrac{\partial^{2}h_{31}}{\partial x_{1}\partial
x_{2}}-\cfrac{\partial^{2}h_{33}}{\partial x_{2}\partial x_{3}}\ ,\
w_{212}=\cfrac{\partial^{3}h_{33}}{\partial
x_{1}^{3}}+\cfrac{\partial^{3}h_{33}}{\partial x_{1}\partial
x_{2}^{2}}-\cfrac{\partial^{3}h_{31}}{\partial x_{1}^{2}\partial
x_{3}}-\cfrac{\partial^{3}h_{32}}{\partial x_{1}\partial
x_{2}\partial x_{3}}\
,\\&&w_{213}=\cfrac{\partial^{3}h_{33}}{\partial x_{1}^{2}\partial
x_{2}}+\cfrac{\partial^{3}h_{33}}{\partial x_{2}\partial
x_{3}^{2}}-\cfrac{\partial^{3}h_{31}}{\partial x_{1}\partial
x_{2}\partial x_{3}}-\cfrac{\partial^{3}h_{32}}{\partial
x_{2}^{2}\partial x_{3}}\ ,\
w_{214}=\cfrac{\partial^{3}h_{33}}{\partial x_{1}^{2}\partial
x_{3}}+\cfrac{\partial^{3}h_{33}}{\partial
x_{3}^{3}}-\cfrac{\partial^{3}h_{31}}{\partial x_{1}\partial
x_{3}^{2}}-\cfrac{\partial^{3}h_{32}}{\partial x_{2}\partial
x_{3}^{2}}\ ,\\&&w_{215}=\cfrac{\partial^{2}h_{33}}{\partial
x_{1}^{2}}+\cfrac{\partial^{2}h_{33}}{\partial
x_{2}^{2}}-\cfrac{\partial^{2}h_{31}}{\partial x_{1}\partial
x_{3}}-\cfrac{\partial^{3}h_{32}}{\partial x_{2}\partial x_{3}}\
,\\&& w_{216}=\cfrac{1}{\tau}\ (\cfrac{\partial^{2}h_{31}}{\partial
t\partial x_{1}}+\cfrac{\partial^{2}h_{32}}{\partial t\partial
x_{2}}+\cfrac{\partial^{2}h_{33}}{\partial t\partial
x_{3}})-\cfrac{\mu}{\tau}\ (\cfrac{\partial\triangle
h_{31}}{\partial x_{1}}+\cfrac{\partial\triangle h_{32}}{\partial
x_{2}}+\cfrac{\partial\triangle h_{33}}{\partial x_{3}})\ .
\end{eqnarray*}
We can see that\ $W_{2}$\ is the function to do the partial
derivation with the components of\ $h_{2}$\ no more than the third
order and their linear combination, moreover no more than the first
order with the variable\ $t$\ . Hence\
$H[F^{-1}(A_{1}(\eta)Z_{1})]=H[F^{-1}(A_{1}(\eta)Z_{1})]I_{K_{1}^{\prime}}\in
C^{1}(K_{1}^{\prime})$\ .\\Now we can get\
$H[F^{-1}(Y_{1}+A_{1}(\eta)Z_{1})]=H[F^{-1}(Y_{1}+A_{1}(\eta)Z_{1})]I_{K_{1}^{\prime}}\in
C^{1}(K_{1}^{\prime})$\ , if\ $Z_{1}\in \Omega_{2}$\ , this means
that\ $\Omega_{2}\subset\Omega_{1}$\ , moreover the components of\
$F^{-1}(f(Z_{1}))$\ are all in the\ \
$H[F^{-1}(Y_{1}+A_{1}(\eta)Z_{1})]\ ,$ hence\
$f(Z_{1})\in\Omega_{2}$\ , if\ $Z_{1}\in\Omega_{2}$\ .
And we can get if\ $Z_{1}\in\Omega_{1}$\ and\ $Z_{1}=f(Z_{1})$\ , then\ $Z_{1}\in\Omega_{2}$\ .\\
Secondly we prove the fixed-point of\ $f(Z_{1})$\ is exist, where\
$Z_{1}\in\Omega_{2}$\ . In order to discuss more conveniently, we
assume\
$f_{j}(Z_{1})=F[F^{-1}(\alpha_{j1}^{T}(Y_{1}+A_{1}(\eta)Z_{1}))F^{-1}(\alpha_{j2}^{T}(Y_{1}+A_{1}(\eta)Z_{1}))]\
,\ 1\leq j\leq 9$\ , where\ $\alpha_{11}=e_{3}\ ,\
\alpha_{21}=e_{7}\ ,\ \alpha_{31}=e_{11}\ ,\ \alpha_{41}=e_{3}\ ,\
\alpha_{51}=e_{7}\ ,\ \alpha_{61}=e_{11}\ ,\ \alpha_{71}=e_{3}\ ,\
\alpha_{81}=e_{7}\ ,$\\$ \alpha_{91}=e_{11}\ ,$\
$\alpha_{12}=-\alpha_{1}\ ,\ \alpha_{22}=e_{1}\ ,\
\alpha_{32}=e_{2}\ ,\ \alpha_{42}=e_{4}\ ,\ \alpha_{52}=e_{5}\ ,\
\alpha_{62}=e_{6}\ ,\ \alpha_{72}=e_{8}\ ,$\\$ \alpha_{82}=e_{9}\ ,\
\alpha_{92}=e_{10}\ ,$ and\ $\alpha_{1}=e_{5}+e_{10}$\ ,\ $e_{i}$\
is the\ $ith$\ $55$\ dimensional unit coordinate vector,\ $1 \leq i
\leq 55$\ .\\ And we assume the set\ $\Omega_{C}$\ as
follows.\begin{eqnarray*}\Omega_{C}&=&\{h_{1}|\ h_{1}\in
C(K_{1}^{\prime})\ ,\ \exists\ M>0\ ,\ \mbox{such that}\ |h_{1}|\leq
M\ ,\ \mbox{moreover}\ \exists\ C >0\ ,\ 0<\alpha< 1\ ,\\&&
\mbox{such that}\ \forall\ x_{t}\ ,\ x_{t}^{\prime}\in
K_{1}^{\prime}\ ,\ |h_{1}(x_{t})-h_{1}(x_{t}^{\prime})|\leq
C|x_{t}-x_{t}^{\prime}|^{\alpha}\ .\}\ ,\\&& \mbox{where}\
h_{1}=(h_{1j})_{9\times1}\ ,\ |h_{1}|=\max_{1\leq j\leq 9}|h_{1j}|\
,\ |h_{1j}|=\max_{x_{t}\in K_{1}^{\prime}}|h_{1j}(x_{t})|\ ,\
x_{t}=(t,\ x_{1},\ x_{2},\ x_{3})^{T}\ .
\end{eqnarray*}
From the Arzela-Ascoli theorem, we know\ $\Omega_{C}$\ is a compact
set. And we assume\ $M$\ in the\ $\Omega_{C}$\ satisfies the
following.\begin{eqnarray*}&&\max_{1\leq j\leq 9}\{\
|F^{-1}(\alpha_{j1}^{T}Y_{1})|\ ,|F^{-1}(\alpha_{j2}^{T}Y_{1})|\
,|F^{-1}(\alpha_{j1}^{T}Y_{1})F^{-1}(\alpha_{j2}^{T}Y_{1})|\}\leq
\theta M\ ,\ \mbox{where}\ 0<\theta<1\ ,\\&&
|F^{-1}(\alpha_{j1}^{T}Y_{1})|=\max_{x_{t}\in
K_{1}^{\prime}}|F^{-1}(\alpha_{j1}^{T}Y_{1})(x_{t})|\ ,\
|F^{-1}(\alpha_{j2}^{T}Y_{1})|=\max_{x_{t}\in
K_{1}^{\prime}}|F^{-1}(\alpha_{j2}^{T}Y_{1})(x_{t})|\
,\\&&|F^{-1}(\alpha_{j1}^{T}Y_{1})F^{-1}(\alpha_{j2}^{T}Y_{1})|=\max_{x_{t}\in
K_{1}^{\prime}}|F^{-1}(\alpha_{j1}^{T}Y_{1})(x_{t})F^{-1}(\alpha_{j2}^{T}Y_{1})(x_{t})|\
,\ 1\leq j\leq 9\ .\end{eqnarray*}We assume\ $\forall\ x_{t}\ ,\
x_{t}^{\prime}\in K_{1}^{\prime}\ ,\ \forall\ j\ ,\ 1\leq j\leq 9\
,\ \exists\ C_{1}> 0$\ , such that\begin{eqnarray*}&&
|F^{-1}(\alpha_{j1}^{T}Y_{1})(x_{t})-F^{-1}(\alpha_{j1}^{T}Y_{1})(x_{t}^{\prime})|\leq
C_{1}|x_{t}-x_{t}^{\prime}|^{\alpha}\ ,\\&&
|F^{-1}(\alpha_{j2}^{T}Y_{1})(x_{t})-F^{-1}(\alpha_{j2}^{T}Y_{1})(x_{t}^{\prime})|\leq
C_{1}|x_{t}-x_{t}^{\prime}|^{\alpha}\ ,\\&&
|F^{-1}(\alpha_{j1}^{T}Y_{1})(x_{t})F^{-1}(\alpha_{j2}^{T}Y_{1})(x_{t})-F^{-1}(\alpha_{j1}^{T}Y_{1})(x_{t}^{\prime})F^{-1}(\alpha_{j2}^{T}Y_{1})(x_{t}^{\prime})|\leq
C_{1}|x_{t}-x_{t}^{\prime}|^{\alpha}\ .\end{eqnarray*}Next we
assume\ $g(h_{1})=F^{-1}[f(F(h_{1}))]$, and\
$g_{j}(h_{1})=F^{-1}[f_{j}(F(h_{1}))],\ 1\leq j\leq 9$\ , where\
$h_{1}\in \Omega_{C}$. We will prove that the fixed-point of\
$g(h_{1})$\ is exist. Moreover we can get\ $h_{1}\in
C^{1}(K_{1}^{\prime})$\ , if\ $h_{1}=g(h_{1})$\ and\ $h_{1}\in
\Omega_{C}$\ .\\We only need to show\ $|g_{j}(h_{1})|\leq M$\ ,
moreover\ $\forall\ x_{t}\ ,\ x_{t}^{\prime}\in
K_{1}^{\prime}$\ , we can get\\
$|g_{j}(x_{t})-g_{j}(x_{t}^{\prime})|\leq
C|x_{t}-x_{t}^{\prime}|^{\alpha}\ ,\ 1\leq j\leq 9$\ . We can see
that
\begin{eqnarray*}
g_{j}(h_{1})&=&F^{-1}[\alpha_{j1}^{T}(Y_{1}+A_{1}(\eta)F(h_{1}))]F^{-1}[\alpha_{j2}^{T}(Y_{1}+A_{1}(\eta)F(h_{1}))]\
,\\&=&F^{-1}(\alpha_{j1}^{T}Y_{1})F^{-1}(\alpha_{j2}^{T}Y_{1})+F^{-1}(\alpha_{j1}^{T}Y_{1})F^{-1}(\alpha_{j2}^{T}A_{1}(\eta)F(h_{1}))+
\\&&F^{-1}(\alpha_{j2}^{T}Y_{1})F^{-1}(\alpha_{j1}^{T}A_{1}(\eta)F(h_{1}))+
\\&&F^{-1}(\alpha_{j1}^{T}A_{1}(\eta)F(h_{1}))F^{-1}(\alpha_{j2}^{T}A_{1}(\eta)F(h_{1}))\
,\end{eqnarray*} and from the lemma 3.1 , we can get\ $\exists\
h_{2}=h_{2}I_{K_{1}^{\prime}}\in C^{1}[0,\ T]\bigcap
C^{\infty}(K_{1})$\ , such that\begin{eqnarray*}
g_{j}(h_{1})&=&F^{-1}(\alpha_{j1}^{T}Y_{1})F^{-1}(\alpha_{j2}^{T}Y_{1})+F^{-1}(\alpha_{j1}^{T}Y_{1})F^{-1}(\alpha_{j2}^{T}A_{1}(\eta)a^{2}b\tau
F(h_{2}))+
\\&&F^{-1}(\alpha_{j2}^{T}Y_{1})F^{-1}(\alpha_{j1}^{T}A_{1}(\eta)a^{2}b\tau F(h_{2}))+
\\&&F^{-1}(\alpha_{j1}^{T}A_{1}(\eta)a^{2}b\tau F(h_{2}))F^{-1}(\alpha_{j2}^{T}A_{1}(\eta)a^{2}b\tau F(h_{2}))
\\&=&F^{-1}(\alpha_{j1}^{T}Y_{1})F^{-1}(\alpha_{j2}^{T}Y_{1})+F^{-1}(\alpha_{j1}^{T}Y_{1})W_{j2}(h_{2})+
\\&&F^{-1}(\alpha_{j2}^{T}Y_{1})W_{j1}(h_{2})+W_{j1}(h_{2})W_{j2}(h_{2})\ ,\end{eqnarray*}
where\ $W_{j1}(h_{2})=F^{-1}(\alpha_{j1}^{T}(a^{2}b\tau
A_{1}(\eta))F(h_{2}))\ ,\
W_{j2}(h_{2})=F^{-1}(\alpha_{j2}^{T}(a^{2}b\tau
A_{1}(\eta))F(h_{2}))\ ,\ 1\leq j\leq 9\ ,$\ and\
$\alpha_{11}=e_{3}\ ,\ \alpha_{21}=e_{7}\ ,\ \alpha_{31}=e_{11}\ ,\
\alpha_{41}=e_{3}\ ,\ \alpha_{51}=e_{7}\ ,\ \alpha_{61}=e_{11}\ ,\
\alpha_{71}=e_{3}\ ,\ \alpha_{81}=e_{7}\ ,$\\$ \alpha_{91}=e_{11}\
,$\ $\alpha_{12}=-\alpha_{1}\ ,\ \alpha_{22}=e_{1}\ ,\
\alpha_{32}=e_{2}\ ,\ \alpha_{42}=e_{4}\ ,\ \alpha_{52}=e_{5}\ ,\
\alpha_{62}=e_{6}\ ,\ \alpha_{72}=e_{8}\ ,$\\$ \alpha_{82}=e_{9}\ ,\
\alpha_{92}=e_{10}\ ,$ and\ $\alpha_{1}=e_{5}+e_{10}$\ ,\ $e_{i}$\
is the\ $ith$\ $55$\ dimensional unit coordinate vector,\ $1 \leq i
\leq 55$\ . And we can get the following.
\begin{eqnarray*}F^{-1}(\alpha_{11}^{T}Y_{1})&=&w_{17}=\cfrac{\partial^{2}F_{21}}{\partial
x_{1}\partial x_{2}}+\cfrac{\partial^{2}F_{31}}{\partial
x_{1}\partial x_{3}}-\cfrac{\partial^{2}F_{11}}{\partial
x_{2}^{2}}-\cfrac{\partial^{2}F_{11}}{\partial x_{3}^{2}}\
,\\F^{-1}(\alpha_{12}^{T}Y_{1})&=&w_{11}=\cfrac{\partial^{3}F_{21}}{\partial
x_{1}^{2}\partial x_{2}}+\cfrac{\partial^{3}F_{31}}{\partial
x_{1}^{2}\partial x_{3}}-\cfrac{\partial^{3}F_{11}}{\partial
x_{1}\partial x_{2}^{2}}-\cfrac{\partial^{3}F_{11}}{\partial
x_{1}\partial x_{3}^{2}}\
,\\F^{-1}(\alpha_{21}^{T}Y_{1})&=&w_{111}=\cfrac{\partial^{2}F_{11}}{\partial
x_{1}\partial x_{2}}+\cfrac{\partial^{2}F_{31}}{\partial
x_{2}\partial x_{3}}-\cfrac{\partial^{2}F_{21}}{\partial
x_{1}^{2}}-\cfrac{\partial^{2}F_{21}}{\partial x_{3}^{2}}\
,\\F^{-1}(\alpha_{22}^{T}Y_{1})&=&w_{15}=\cfrac{\partial^{3}F_{21}}{\partial
x_{1}\partial x_{2}^{2}}+\cfrac{\partial^{3}F_{31}}{\partial
x_{1}\partial x_{2}\partial
x_{3}}-\cfrac{\partial^{3}F_{11}}{\partial
x_{2}^{3}}-\cfrac{\partial^{3}F_{11}}{\partial x_{2}\partial
x_{3}^{2}}\
,\\F^{-1}(\alpha_{31}^{T}Y_{1})&=&w_{115}=\cfrac{\partial^{2}F_{11}}{\partial
x_{1}\partial x_{3}}+\cfrac{\partial^{2}F_{21}}{\partial
x_{2}\partial x_{3}}-\cfrac{\partial^{2}F_{31}}{\partial
x_{1}^{2}}-\cfrac{\partial^{2}F_{31}}{\partial x_{2}^{2}}\
,\\F^{-1}(\alpha_{32}^{T}Y_{1})&=&w_{16}=\cfrac{\partial^{3}F_{21}}{\partial
x_{1}\partial x_{2}\partial
x_{3}}+\cfrac{\partial^{3}F_{31}}{\partial x_{1}\partial
x_{3}^{2}}-\cfrac{\partial^{3}F_{11}}{\partial x_{2}^{2}\partial
x_{3}}-\cfrac{\partial^{3}F_{11}}{ \partial x_{3}^{3}}\
,\\F^{-1}(\alpha_{41}^{T}Y_{1})&=&w_{17}=\cfrac{\partial^{2}F_{21}}{\partial
x_{1}\partial x_{2}}+\cfrac{\partial^{2}F_{31}}{\partial
x_{1}\partial x_{3}}-\cfrac{\partial^{2}F_{11}}{\partial
x_{2}^{2}}-\cfrac{\partial^{2}F_{11}}{\partial x_{3}^{2}}\
,\\F^{-1}(\alpha_{42}^{T}Y_{1})&=&w_{18}=\cfrac{\partial^{3}F_{11}}{\partial
x_{1}^{2}\partial x_{2}}+\cfrac{\partial^{3}F_{31}}{\partial
x_{1}\partial x_{2}\partial
x_{3}}-\cfrac{\partial^{3}F_{21}}{\partial
x_{1}^{3}}-\cfrac{\partial^{3}F_{21}}{\partial x_{1}\partial
x_{3}^{2}}\
,\\F^{-1}(\alpha_{51}^{T}Y_{1})&=&w_{111}=\cfrac{\partial^{2}F_{11}}{\partial
x_{1}\partial x_{2}}+\cfrac{\partial^{2}F_{31}}{\partial
x_{2}\partial x_{3}}-\cfrac{\partial^{2}F_{21}}{\partial
x_{1}^{2}}-\cfrac{\partial^{2}F_{21}}{\partial x_{3}^{2}}\
,\\F^{-1}(\alpha_{52}^{T}Y_{1})&=&w_{19}=\cfrac{\partial^{3}F_{11}}{\partial
x_{1}\partial x_{2}^{2}}+\cfrac{\partial^{3}F_{31}}{\partial
x_{2}^{2}\partial x_{3}}-\cfrac{\partial^{3}F_{21}}{\partial
x_{1}^{2}\partial x_{2}}-\cfrac{\partial^{3}F_{21}}{\partial
x_{2}\partial x_{3}^{2}}\
,\\F^{-1}(\alpha_{61}^{T}Y_{1})&=&w_{115}=\cfrac{\partial^{2}F_{11}}{\partial
x_{1}\partial x_{3}}+\cfrac{\partial^{2}F_{21}}{\partial
x_{2}\partial x_{3}}-\cfrac{\partial^{2}F_{31}}{\partial
x_{1}^{2}}-\cfrac{\partial^{2}F_{31}}{\partial x_{2}^{2}}\
,\\F^{-1}(\alpha_{62}^{T}Y_{1})&=&w_{110}=\cfrac{\partial^{3}F_{11}}{\partial
x_{1}\partial x_{2}\partial x_{3}}+\cfrac{\partial^{3}F_{31}}{
\partial x_{2}\partial
x_{3}^{2}}-\cfrac{\partial^{3}F_{21}}{\partial x_{1}^{2}\partial
x_{3}}-\cfrac{\partial^{3}F_{21}}{ \partial x_{3}^{3}}\
,\\F^{-1}(\alpha_{71}^{T}Y_{1})&=&w_{17}=\cfrac{\partial^{2}F_{21}}{\partial
x_{1}\partial x_{2}}+\cfrac{\partial^{2}F_{31}}{\partial
x_{1}\partial x_{3}}-\cfrac{\partial^{2}F_{11}}{\partial
x_{2}^{2}}-\cfrac{\partial^{2}F_{11}}{\partial x_{3}^{2}}\
,\\F^{-1}(\alpha_{72}^{T}Y_{1})&=&w_{112}=\cfrac{\partial^{3}F_{11}}{\partial
x_{1}^{2}\partial x_{3}}+\cfrac{\partial^{3}F_{21}}{\partial
x_{1}\partial x_{2}\partial
x_{3}}-\cfrac{\partial^{3}F_{31}}{\partial
x_{1}^{3}}-\cfrac{\partial^{3}F_{31}}{\partial x_{1}\partial
x_{2}^{2}}\
,\\F^{-1}(\alpha_{81}^{T}Y_{1})&=&w_{111}=\cfrac{\partial^{2}F_{11}}{\partial
x_{1}\partial x_{2}}+\cfrac{\partial^{2}F_{31}}{\partial
x_{2}\partial x_{3}}-\cfrac{\partial^{2}F_{21}}{\partial
x_{1}^{2}}-\cfrac{\partial^{2}F_{21}}{\partial x_{3}^{2}}\
,\\F^{-1}(\alpha_{82}^{T}Y_{1})&=&w_{113}=\cfrac{\partial^{3}F_{11}}{\partial
x_{1}\partial x_{2}\partial x_{3}}+\cfrac{\partial^{3}F_{21}}{
\partial x_{2}^{2}\partial
x_{3}}-\cfrac{\partial^{3}F_{31}}{\partial x_{1}^{2}\partial
x_{2}}-\cfrac{\partial^{3}F_{31}}{ \partial x_{2}^{3}}\
,\\F^{-1}(\alpha_{91}^{T}Y_{1})&=&w_{115}=\cfrac{\partial^{2}F_{11}}{\partial
x_{1}\partial x_{3}}+\cfrac{\partial^{2}F_{21}}{\partial
x_{2}\partial x_{3}}-\cfrac{\partial^{2}F_{31}}{\partial
x_{1}^{2}}-\cfrac{\partial^{2}F_{31}}{\partial x_{2}^{2}}\
,\\F^{-1}(\alpha_{92}^{T}Y_{1})&=&w_{114}=\cfrac{\partial^{3}F_{11}}{\partial
x_{1} \partial x_{3}^{2}}+\cfrac{\partial^{3}F_{21}}{ \partial
x_{2}\partial x_{3}^{2}}-\cfrac{\partial^{3}F_{31}}{\partial
x_{1}^{2}\partial x_{3}}-\cfrac{\partial^{3}F_{31}}{ \partial
x_{2}^{2}\partial x_{3}}\
,\end{eqnarray*}\begin{eqnarray*}F_{j1}(t,M)&=&(\frac{1}{2\sqrt{\pi\mu}})^{3}\int_{0}^{t}\int_{R^{3}}\frac{F_{jv}(\tau_{1},
y_{1}, y_{2} , y_{3})}{(\sqrt{t-\tau_{1}})^{3}}\
e^{-\frac{(x_{1}-y_{1})^{2}+(x_{2}-y_{2})^{2}+(x_{3}-y_{3})^{2}}{4\mu(t-\tau_{1})}}
dy_{1}dy_{2}dy_{3}d\tau_{1}\
,\\F_{jv}(t,M_{0})&=&-\frac{1}{4\pi}\int_{K_{1}}\frac{F_{j}(t,M)}{r_{MM_{0}}}
dx_{1}dx_{2}dx_{3}\ ,\ 1\leq\ j\leq 3\ ,\\
W_{11}(h_{2})&=&w_{27}=\cfrac{\partial^{2}h_{31}}{\partial
x_{2}^{2}}+\cfrac{\partial^{2}h_{31}}{\partial
x_{3}^{2}}-\cfrac{\partial^{2}h_{32}}{\partial x_{1}\partial
x_{2}}-\cfrac{\partial^{2}h_{33}}{\partial x_{1}\partial x_{3}}\
,\\W_{12}(h_{2})&=&w_{21}=\cfrac{\partial^{3}h_{31}}{\partial
x_{1}\partial x_{2}^{2}}+\cfrac{\partial^{3}h_{31}}{\partial
x_{1}\partial x_{3}^{2}}-\cfrac{\partial^{3}h_{32}}{\partial
x_{1}^{2}\partial x_{2}}-\cfrac{\partial^{3}h_{33}}{\partial
x_{1}^{2}\partial x_{3}}\
,\\W_{21}(h_{2})&=&w_{211}=\cfrac{\partial^{2}h_{32}}{\partial
x_{1}^{2}}+\cfrac{\partial^{2}h_{32}}{\partial
x_{3}^{2}}-\cfrac{\partial^{2}h_{31}}{\partial x_{1}\partial
x_{2}}-\cfrac{\partial^{2}h_{33}}{\partial x_{2}\partial x_{3}}\
,\\W_{22}(h_{2})&=&w_{25}=\cfrac{\partial^{3}h_{31}}{\partial
x_{2}^{3}}+\cfrac{\partial^{3}h_{31}}{\partial x_{2}\partial
x_{3}^{2}}-\cfrac{\partial^{3}h_{32}}{\partial x_{1}\partial
x_{2}^{2}}-\cfrac{\partial^{3}h_{33}}{\partial x_{1}\partial
x_{2}\partial x_{3}}\
,\\W_{31}(h_{2})&=&w_{215}=\cfrac{\partial^{2}h_{33}}{\partial
x_{1}^{2}}+\cfrac{\partial^{2}h_{33}}{\partial
x_{2}^{2}}-\cfrac{\partial^{2}h_{31}}{\partial x_{1}\partial
x_{3}}-\cfrac{\partial^{3}h_{32}}{\partial x_{2}\partial x_{3}}\
,\\W_{32}(h_{2})&=&w_{26}=\cfrac{\partial^{3}h_{31}}{\partial
x_{2}^{2}\partial x_{3}}+\cfrac{\partial^{3}h_{31}}{\partial
x_{3}^{3}}-\cfrac{\partial^{3}h_{32}}{\partial x_{1}\partial
x_{2}\partial x_{3}}-\cfrac{\partial^{3}h_{33}}{\partial
x_{1}\partial x_{3}^{2}}\
,\\W_{41}(h_{2})&=&w_{27}=\cfrac{\partial^{2}h_{31}}{\partial
x_{2}^{2}}+\cfrac{\partial^{2}h_{31}}{\partial
x_{3}^{2}}-\cfrac{\partial^{2}h_{32}}{\partial x_{1}\partial
x_{2}}-\cfrac{\partial^{2}h_{33}}{\partial x_{1}\partial x_{3}}\
,\\W_{42}(h_{2})&=&w_{28}=\cfrac{\partial^{3}h_{32}}{\partial
x_{1}^{3}}+\cfrac{\partial^{3}h_{32}}{\partial x_{1}\partial
x_{3}^{2}}-\cfrac{\partial^{3}h_{31}}{\partial x_{1}^{2}\partial
x_{2}}-\cfrac{\partial^{3}h_{33}}{\partial x_{1}\partial
x_{2}\partial x_{3}}\
,\\W_{51}(h_{2})&=&w_{211}=\cfrac{\partial^{2}h_{32}}{\partial
x_{1}^{2}}+\cfrac{\partial^{2}h_{32}}{\partial
x_{3}^{2}}-\cfrac{\partial^{2}h_{31}}{\partial x_{1}\partial
x_{2}}-\cfrac{\partial^{2}h_{33}}{\partial x_{2}\partial x_{3}}\
,\\W_{52}(h_{2})&=&w_{29}=\cfrac{\partial^{3}h_{32}}{\partial
x_{1}^{2}\partial x_{2}}+\cfrac{\partial^{3}h_{32}}{\partial
x_{2}\partial x_{3}^{2}}-\cfrac{\partial^{3}h_{31}}{\partial
x_{1}\partial x_{2}^{2}}-\cfrac{\partial^{3}h_{33}}{\partial
x_{2}^{2}\partial x_{3}}\
,\\W_{61}(h_{2})&=&w_{215}=\cfrac{\partial^{2}h_{33}}{\partial
x_{1}^{2}}+\cfrac{\partial^{2}h_{33}}{\partial
x_{2}^{2}}-\cfrac{\partial^{2}h_{31}}{\partial x_{1}\partial
x_{3}}-\cfrac{\partial^{3}h_{32}}{\partial x_{2}\partial x_{3}}\
,\\W_{62}(h_{2})&=&w_{210}=\cfrac{\partial^{3}h_{32}}{\partial
x_{1}^{2}\partial x_{3}}+\cfrac{\partial^{3}h_{32}}{\partial
x_{3}^{3}}-\cfrac{\partial^{3}h_{31}}{\partial x_{1}\partial
x_{2}\partial x_{3}}-\cfrac{\partial^{3}h_{33}}{\partial
x_{2}\partial x_{3}^{2}}\
,\\W_{71}(h_{2})&=&w_{27}=\cfrac{\partial^{2}h_{31}}{\partial
x_{2}^{2}}+\cfrac{\partial^{2}h_{31}}{\partial
x_{3}^{2}}-\cfrac{\partial^{2}h_{32}}{\partial x_{1}\partial
x_{2}}-\cfrac{\partial^{2}h_{33}}{\partial x_{1}\partial x_{3}}\
,\\W_{72}(h_{2})&=&w_{212}=\cfrac{\partial^{3}h_{33}}{\partial
x_{1}^{3}}+\cfrac{\partial^{3}h_{33}}{\partial x_{1}\partial
x_{2}^{2}}-\cfrac{\partial^{3}h_{31}}{\partial x_{1}^{2}\partial
x_{3}}-\cfrac{\partial^{3}h_{32}}{\partial x_{1}\partial
x_{2}\partial x_{3}}\
,\\W_{81}(h_{2})&=&w_{211}=\cfrac{\partial^{2}h_{32}}{\partial
x_{1}^{2}}+\cfrac{\partial^{2}h_{32}}{\partial
x_{3}^{2}}-\cfrac{\partial^{2}h_{31}}{\partial x_{1}\partial
x_{2}}-\cfrac{\partial^{2}h_{33}}{\partial x_{2}\partial x_{3}}\
,\\W_{82}(h_{2})&=&w_{213}=\cfrac{\partial^{3}h_{33}}{\partial
x_{1}^{2}\partial x_{2}}+\cfrac{\partial^{3}h_{33}}{\partial
x_{2}\partial x_{3}^{2}}-\cfrac{\partial^{3}h_{31}}{\partial
x_{1}\partial x_{2}\partial
x_{3}}-\cfrac{\partial^{3}h_{32}}{\partial x_{2}^{2}\partial x_{3}}\
,\\W_{91}(h_{2})&=&w_{215}=\cfrac{\partial^{2}h_{33}}{\partial
x_{1}^{2}}+\cfrac{\partial^{2}h_{33}}{\partial
x_{2}^{2}}-\cfrac{\partial^{2}h_{31}}{\partial x_{1}\partial
x_{3}}-\cfrac{\partial^{3}h_{32}}{\partial x_{2}\partial x_{3}}\
,\\W_{92}(h_{2})&=&w_{214}=\cfrac{\partial^{3}h_{33}}{\partial
x_{1}^{2}\partial x_{3}}+\cfrac{\partial^{3}h_{33}}{\partial
x_{3}^{3}}-\cfrac{\partial^{3}h_{31}}{\partial x_{1}\partial
x_{3}^{2}}-\cfrac{\partial^{3}h_{32}}{\partial x_{2}\partial
x_{3}^{2}}\ ,
\end{eqnarray*}where\ $h_{2}=(h_{2j})_{9\times 1}$\ , and\ $h_{31}=h_{21}+h_{22}+h_{23}\ ,\ h_{32}=h_{24}+h_{25}+h_{26}\
,\ h_{33}=h_{27}+h_{28}+h_{29}$\ .\\ Hence we can get\
$W_{j1}(h_{2})\ ,\ W_{j2}(h_{2})\ ,\ 1\leq j\leq 9$\ , are the
functions to do the partial derivation with the components of\
$h_{2}$\ no more than the third order and their linear combination,
moreover only do the partial derivation with the variables\ $x_{1}\
,\ x_{2}\ ,\ x_{3}$\ . From the lemma 3.1 , we
know\begin{eqnarray*}&&v(t,M_{0})=-\frac{1}{4\pi}\int_{K_{1}}\frac{h_{1}(t,
M)}{r_{MM_{0}}} dx_{1}dx_{2}dx_{3}\ ,\ \mbox{where}\
M=(x_{1},x_{2},x_{3})\ ,\ M_{0}=(x_{10},x_{20},x_{30})\
,\\&&r_{MM_{0}}=\sqrt{(x_{1}-x_{10})^{2}+
(x_{2}-x_{20})^{2}+(x_{3}-x_{30})^{2}}\ ,\\&&h_{2}(t,
M)=(\frac{1}{2\sqrt{\pi\mu}})^{3}\int_{0}^{t}\int_{R^{3}}\frac{v(\tau_{1},
y_{1}, y_{2} , y_{3})}{(\sqrt{t-\tau_{1}})^{3}}\
e^{-\frac{(x_{1}-y_{1})^{2}+(x_{2}-y_{2})^{2}+(x_{3}-y_{3})^{2}}{4\mu(t-\tau_{1})}}
dy_{1}dy_{2}dy_{3}d\tau_{1}\ . \end{eqnarray*}  And we can get the
lemma as follows.
\begin{lemma} \label{Lemma-3.2} $\exists\ M_{T,\ 1}>0\ ,\ \mbox{and\ $M_{T,\ 1}$\ is independent with\ $M\ ,\ K_{1}$}\ ,\ \forall\ h_{1}\in\
\Omega_{C}$\ , we can get the following, $$\max_{1\leq
i\leq3}\{|h_{2}|\ ,\ |\cfrac{\partial h_{2}}{\partial t}\ |\ ,\
|\cfrac{\partial h_{2}}{\partial x_{i}}\ |\ ,\
|\cfrac{\partial\cfrac{\partial h_{2}}{\partial x_{i}}}{\partial t}\
|\}\leq\ MM_{T,\ 1}M(K_{1})\ ,$$ \mbox{where}\begin{eqnarray*}&&
h_{1}=(h_{1m})_{9\times1}\ ,\ h_{2}=h_{2}(t,\ M)=(h_{2m}(t,\
M))_{9\times1}\ ,\ |h_{2}|=\max_{1\leq m\leq 9}\max_{(t,\ M)\in\
K_{1}^{\prime}}|h_{2m}(t,\ M)|\ ,\\&& |\cfrac{\partial
h_{2}}{\partial t}\ |=\max_{1\leq m\leq 9}\max_{(t,\ M)\in\
K_{1}^{\prime}}|\cfrac{\partial h_{2m}(t,\ M)}{\partial t}\ |\ ,\
|\cfrac{\partial h_{2}}{\partial x_{i}}\ |=\max_{1\leq m\leq
9}\max_{(t,\ M)\in\ K_{1}^{\prime}}|\cfrac{\partial h_{2m}(t,\
M)}{\partial x_{i}}\ |\ ,\\&&|\cfrac{\partial\cfrac{\partial
h_{2}}{\partial x_{i}}}{\partial t}\ |=\max_{1\leq m\leq
9}\max_{(t,\ M)\in\ K_{1}^{\prime}}|\cfrac{\partial\cfrac{\partial
h_{2m}(t,\ M)}{\partial x_{i}}}{\partial t}\ |\ ,\
M(K_{1})=\max_{M_{0}\in
K_{1}}\frac{1}{4\pi}\int_{K_{1}}\frac{1}{r_{MM_{0}}}
dx_{1}dx_{2}dx_{3}\ .\end{eqnarray*}
\end{lemma}{\it Proof of lemma3-2}. We denote\ $h_{2m}=h_{2m}(t,\ M)\ ,\ 1\leq m\leq9$\ ,
and we assume\ $v(t,M)=(v_{m})_{9\times 1}$\ , where\
$v_{m}=v_{m}(t,M)$\ , then from the lemma 3.1 , we can get\
$\forall\ m\ ,\ 1\leq m\leq9$\ ,
\begin{eqnarray*}&&v_{m}(t,M_{0})=-\frac{1}{4\pi}\int_{K_{1}}\frac{h_{1m}(t,
M)}{r_{MM_{0}}} dx_{1}dx_{2}dx_{3}\ ,\ \mbox{where}\
M=(x_{1},x_{2},x_{3})\ ,\ M_{0}=(x_{10},x_{20},x_{30})\
,\\&&r_{MM_{0}}=\sqrt{(x_{1}-x_{10})^{2}+
(x_{2}-x_{20})^{2}+(x_{3}-x_{30})^{2}}\ ,\\&&h_{2m}(t,
M)=(\frac{1}{2\sqrt{\pi\mu}})^{3}\int_{0}^{t}\int_{R^{3}}\frac{v_{m}(\tau_{1},
y_{1}, y_{2} , y_{3})}{(\sqrt{t-\tau_{1}})^{3}}\
e^{-\frac{(x_{1}-y_{1})^{2}+(x_{2}-y_{2})^{2}+(x_{3}-y_{3})^{2}}{4\mu(t-\tau_{1})}}
dy_{1}dy_{2}dy_{3}d\tau_{1}\ . \end{eqnarray*}  if\ $h_{1}\in\
\Omega_{C}$\ , then we can get$$|v_{m}|\leq |h_{1m}|\
\frac{1}{4\pi}\int_{K_{1}}\frac{1}{r_{MM_{0}}}
dx_{1}dx_{2}dx_{3}\leq MM(K_{1})\ , $$moreover
$$|h_{2m}|\leq(\frac{1}{2\sqrt{\pi\mu}})^{3}\int_{0}^{t}\int_{R^{3}}\frac{MM(K_{1})}{(\sqrt{t-\tau_{1}})^{3}}\
e^{-\frac{(x_{1}-y_{1})^{2}+(x_{2}-y_{2})^{2}+(x_{3}-y_{3})^{2}}{4\mu(t-\tau_{1})}}
dy_{1}dy_{2}dy_{3}d\tau_{1}\leq TMM(K_{1})\ ,$$
where$$M(K_{1})=\max_{M_{0}\in
K_{1}}\frac{1}{4\pi}\int_{K_{1}}\frac{1}{r_{MM_{0}}}
dx_{1}dx_{2}dx_{3}\ .$$ And we know\begin{eqnarray*}\cfrac{\partial
h_{2m}}{\partial t}&=&\lim_{\tau_{1}\rightarrow
t}(\frac{1}{2\sqrt{\pi\mu}})^{3}\int_{R^{3}}\frac{v_{m}(\tau_{1},
y_{1}, y_{2} , y_{3})}{(\sqrt{t-\tau_{1}})^{3}}\
e^{-\frac{(x_{1}-y_{1})^{2}+(x_{2}-y_{2})^{2}+(x_{3}-y_{3})^{2}}{4\mu(t-\tau_{1})}}
dy_{1}dy_{2}dy_{3}+\\&&(\frac{1}{2\sqrt{\pi\mu}})^{3}\int_{0}^{t}\int_{R^{3}}\cfrac{\partial(\cfrac{v_{m}(\tau_{1},
y_{1}, y_{2} , y_{3})}{(\sqrt{t-\tau_{1}})^{3}}\
e^{-\frac{(x_{1}-y_{1})^{2}+(x_{2}-y_{2})^{2}+(x_{3}-y_{3})^{2}}{4\mu(t-\tau_{1})}})}{\partial
t}dy_{1}dy_{2}dy_{3}d\tau_{1}\ .\end{eqnarray*}If we let\
$y_{j}=x_{j}+2\sqrt{\mu(t-\tau_{1})}\theta_{j}\ ,\ j=1\ ,\ 2\ ,\ 3$\
, then we can get as follows,\begin{eqnarray*}\cfrac{\partial
h_{2m}}{\partial t}&=&\lim_{\tau_{1}\rightarrow
t}(\frac{1}{\sqrt{\pi}})^{3}\int_{R^{3}}v_{m}(\tau_{1},
x_{j}+2(\sqrt{\mu(t-\tau_{1})})\theta_{j},\ j=1,2,3,)\
e^{-(\theta_{1}^{2}+\theta_{2}^{2}+\theta_{3}^{2})}
d\theta_{1}d\theta_{2}d\theta_{3}+\\&&(\frac{1}{2\sqrt{\pi\mu}})^{3}\int_{0}^{t}\int_{R^{3}}\cfrac{-3}{2}\
\frac{v_{m}(\tau_{1}, y_{1}, y_{2} ,
y_{3})}{(\sqrt{t-\tau_{1}})^{5}}\
e^{-\frac{(x_{1}-y_{1})^{2}+(x_{2}-y_{2})^{2}+(x_{3}-y_{3})^{2}}{4\mu(t-\tau_{1})}}dy_{1}dy_{2}dy_{3}d\tau_{1}+
\\&&(\frac{1}{2\sqrt{\pi\mu}})^{3}\int_{0}^{t}\int_{R^{3}}
\frac{v_{m}(\tau_{1}, y_{1}, y_{2} ,
y_{3})}{4\mu(\sqrt{t-\tau_{1}})^{7}}\
\sum_{i=1}^{3}(x_{i}-y_{i})^{2}\
e^{-\frac{\sum_{i=1}^{3}(x_{i}-y_{i})^{2}}{4\mu(t-\tau_{1})}}dy_{1}dy_{2}dy_{3}d\tau_{1}\
.\end{eqnarray*}If we let\
$y_{j}=x_{j}+2\sqrt{\mu}(\sqrt{(t-\tau_{1})})^{5}\theta_{j}\ ,\ j=1\
,\ 2\ ,\ 3\ ,$\ to the second integral, and we let\\
$y_{j}=x_{j}+2\sqrt{\mu}(\sqrt{(t-\tau_{1})})^{7}\theta_{j}\ ,\ j=1\
,\ 2\ ,\ 3\ ,$\ to the third integral, then we can get as follows,
\begin{eqnarray*}\cfrac{\partial
h_{2m}}{\partial t}&=&v_{m}(t,
x_{1},x_{2},x_{3})+(\frac{1}{\sqrt{\pi}})^{3}\int_{0}^{t}\int_{R^{3}}\cfrac{-3}{2}\
v_{m}^{(1)}(t-\tau_{1})^{5}\
e^{-(t-\tau_{1})^{4}(\theta_{1}^{2}+\theta_{2}^{2}+\theta_{3}^{2})}
d\theta_{1}d\theta_{2}d\theta_{3}d\tau_{1}+\\&&
(\frac{1}{\sqrt{\pi}})^{3}\int_{0}^{t}\int_{R^{3}}
v_{m}^{(2)}(t-\tau_{1})^{14}(\theta_{1}^{2}+\theta_{2}^{2}+\theta_{3}^{2})\
e^{-(t-\tau_{1})^{6}(\theta_{1}^{2}+\theta_{2}^{2}+\theta_{3}^{2})}
d\theta_{1}d\theta_{2}d\theta_{3}d\tau_{1}\ ,\end{eqnarray*}
where$$v_{m}^{(1)}=v_{m}(\tau_{1},
x_{j}+2\sqrt{\mu}(\sqrt{(t-\tau_{1})})^{5}\theta_{j},j=1,2,3,)\ ,\
v_{m}^{(2)}=v_{m}(\tau_{1},
x_{j}+2\sqrt{\mu}(\sqrt{(t-\tau_{1})})^{7}\theta_{j},j=1,2,3,)\ .$$
Hence we can get\begin{eqnarray*}|\cfrac{\partial h_{2m}}{\partial
t}\ |&\leq& M
M(K_{1})+(\frac{1}{\sqrt{\pi}})^{3}\int_{0}^{t}\int_{R^{3}}\cfrac{3}{2}\
M M(K_{1})(t-\tau_{1})^{5}\
e^{-(t-\tau_{1})^{4}(\theta_{1}^{2}+\theta_{2}^{2}+\theta_{3}^{2})}
d\theta_{1}d\theta_{2}d\theta_{3}d\tau_{1}+\\&&
(\frac{1}{\sqrt{\pi}})^{3}\int_{0}^{t}\int_{R^{3}}M
M(K_{1})(t-\tau_{1})^{14}(\theta_{1}^{2}+\theta_{2}^{2}+\theta_{3}^{2})\
e^{-(t-\tau_{1})^{6}(\theta_{1}^{2}+\theta_{2}^{2}+\theta_{3}^{2})}
d\theta_{1}d\theta_{2}d\theta_{3}d\tau_{1}\ .\end{eqnarray*} We
assume\ $\varphi_{1}(t,\ \tau_{1})\ ,\ \varphi_{2}(t,\ \tau_{1})$\
as follows,\begin{eqnarray*}&&\varphi_{1}(t,\
\tau_{1})=\int_{R^{3}}(t-\tau_{1})^{5}\
e^{-(t-\tau_{1})^{4}(\theta_{1}^{2}+\theta_{2}^{2}+\theta_{3}^{2})}
d\theta_{1}d\theta_{2}d\theta_{3}\ ,\\&&\varphi_{2}(t,\
\tau_{1})=\int_{R^{3}}(t-\tau_{1})^{14}(\theta_{1}^{2}+\theta_{2}^{2}+\theta_{3}^{2})\
e^{-(t-\tau_{1})^{6}(\theta_{1}^{2}+\theta_{2}^{2}+\theta_{3}^{2})}
d\theta_{1}d\theta_{2}d\theta_{3}\ .\end{eqnarray*}We can see that\
$\varphi_{1}(t,\ \tau_{1}),\ \varphi_{2}(t,\ \tau_{1})$\ are
continuous about\ $t,\ \tau_{1}$\ on the region\ $0\leq \tau_{1}\leq
t\ ,\ 0\leq t\leq T\ ,$ hence they are bounded on the region\ $0\leq
\tau_{1}\leq t\ ,\ 0\leq t\leq T$\ . We assume there exist\
$M^{\prime}>0$\ , such that\ $|\varphi_{1}(t,\ \tau_{1})|\leq
M^{\prime}\ ,\ |\varphi_{2}(t,\ \tau_{1})|\leq M^{\prime}$\ , where\
$0\leq \tau_{1}\leq t\ ,\ 0\leq t\leq T\ .$\ Hence we can
get$$|\cfrac{\partial h_{2m}}{\partial t}\ |\leq M
M(K_{1})+(\frac{1}{\sqrt{\pi}})^{3}\cfrac{5}{2}\
TMM(K_{1})M^{\prime}\ .$$ And we know$$\cfrac{\partial
h_{2m}}{\partial
x_{i}}=(\frac{1}{2\sqrt{\pi\mu}})^{3}\int_{0}^{t}\int_{R^{3}}
\frac{v_{m}(\tau_{1}, y_{1}, y_{2} ,
y_{3})}{4\mu(\sqrt{t-\tau_{1}})^{5}}\ (-2)(x_{i}-y_{i})\
e^{-\frac{\sum_{j=1}^{3}(x_{j}-y_{j})^{2}}{4\mu(t-\tau_{1})}}dy_{1}dy_{2}dy_{3}d\tau_{1}\
.$$If we let\
$y_{j}=x_{j}+2\sqrt{\mu}(\sqrt{(t-\tau_{1})})^{5}\theta_{j}\ ,\ j=1\
,\ 2\ ,\ 3\ ,$\ we can get the
follows,\begin{eqnarray*}\cfrac{\partial h_{2m}}{\partial
x_{i}}&=&(\frac{1}{\sqrt{\pi}})^{3}\int_{0}^{t}\int_{R^{3}}\cfrac{1}{\sqrt{\mu}}\
v_{m}^{(1)}\ (\sqrt{t-\tau_{1}})^{15}\theta_{i}\
e^{-(t-\tau_{1})^{4}(\theta_{1}^{2}+\theta_{2}^{2}+\theta_{3}^{2})}
d\theta_{1}d\theta_{2}d\theta_{3}d\tau_{1}\ .\end{eqnarray*}We
assume\ $\varphi_{3}(t,\ \tau_{1})$\ as follows,$$\varphi_{3}(t,\
\tau_{1})=\int_{R^{3}}(\sqrt{t-\tau_{1}})^{15}\theta_{i}\
e^{-(t-\tau_{1})^{4}(\theta_{1}^{2}+\theta_{2}^{2}+\theta_{3}^{2})}
d\theta_{1}d\theta_{2}d\theta_{3}\ ,$$We can also see that\
$\varphi_{3}(t,\ \tau_{1})$\ is continuous about\ $t,\ \tau_{1}$\ on
the region\ $0\leq \tau_{1}\leq t\ ,\ 0\leq t\leq T\ ,$ hence it is
bounded on the region\ $0\leq \tau_{1}\leq t\ ,\ 0\leq t\leq T$\ .
We assume there exist\ $M^{\prime\prime}>0$\ , such that\
$|\varphi_{3}(t,\ \tau_{1})|\leq M^{\prime\prime}$\ , where\ $0\leq
\tau_{1}\leq t\ ,\ 0\leq t\leq T\ .$\ Hence we can get$$
|\cfrac{\partial h_{2m}}{\partial x_{i}}\ |\leq
(\frac{1}{\sqrt{\pi}})^{3}\cfrac{T}{\sqrt{\mu}}\
MM(K_{1})M^{\prime\prime}\ .$$ At last we see
from\begin{eqnarray*}\cfrac{\partial\cfrac{\partial h_{2m}}{\partial
x_{i}}}{\partial t}&=&\lim_{\tau_{1}\rightarrow
t}(\frac{1}{2\sqrt{\pi\mu}})^{3}\int_{R^{3}} \frac{v_{m}(\tau_{1},
y_{1}, y_{2} , y_{3})}{4\mu(\sqrt{t-\tau_{1}})^{5}}\
(-2)(x_{i}-y_{i})\
e^{-\frac{\sum_{j=1}^{3}(x_{j}-y_{j})^{2}}{4\mu(t-\tau_{1})}}dy_{1}dy_{2}dy_{3}+\\&&
(\frac{1}{2\sqrt{\pi\mu}})^{3}[\int_{0}^{t}\int_{R^{3}}
\frac{v_{m}(\tau_{1}, y_{1}, y_{2} ,
y_{3})}{4\mu(\sqrt{t-\tau_{1}})^{7}}\ (5)(x_{i}-y_{i})\
e^{-\frac{\sum_{j=1}^{3}(x_{j}-y_{j})^{2}}{4\mu(t-\tau_{1})}}dy_{1}dy_{2}dy_{3}d\tau_{1}+
\\&&\int_{0}^{t}\int_{R^{3}}
\frac{v_{m}(\tau_{1}, y_{1}, y_{2} ,
y_{3})}{16\mu^{2}(\sqrt{t-\tau_{1}})^{9}}\
(-2)(x_{i}-y_{i})\sum_{j=1}^{3}(x_{j}-y_{j})^{2}\
e^{-\frac{\sum_{i=1}^{3}(x_{i}-y_{i})^{2}}{4\mu(t-\tau_{1})}}dy_{1}dy_{2}dy_{3}d\tau_{1}]\
.\end{eqnarray*}If we let\
$y_{j}=x_{j}+2\sqrt{\mu}(\sqrt{(t-\tau_{1})})^{5}\theta_{j}\ ,\ j=1\
,\ 2\ ,\ 3\ ,$\ to the first integral, and we let\\
$y_{j}=x_{j}+2\sqrt{\mu}(\sqrt{(t-\tau_{1})})^{7}\theta_{j}\ ,\ j=1\
,\ 2\ ,\ 3\ ,$\ to the second integral,\
$y_{j}=x_{j}+2\sqrt{\mu}(\sqrt{(t-\tau_{1})})^{9}\theta_{j}$\ ,\\$
j=1\ ,\ 2\ ,\ 3\ ,$\ to the third integral, then we can get as
follows,\begin{eqnarray*}\cfrac{\partial\cfrac{\partial
h_{2m}}{\partial x_{i}}}{\partial t}&=&\lim_{\tau_{1}\rightarrow
t}(\frac{1}{\sqrt{\pi}})^{3}\int_{R^{3}}\cfrac{1}{\sqrt{\mu}}\
v_{m}^{(1)}\ (\sqrt{t-\tau_{1}})^{15}\theta_{i}\
e^{-(t-\tau_{1})^{4}(\theta_{1}^{2}+\theta_{2}^{2}+\theta_{3}^{2})}
d\theta_{1}d\theta_{2}d\theta_{3}+\\&&
(\frac{1}{\sqrt{\pi}})^{3}\int_{0}^{t}\int_{R^{3}}
\frac{v_{m}^{(2)}}{2\sqrt{\mu}}(\sqrt{t-\tau_{1}})^{21}\
(5)\theta_{i}\
e^{-(t-\tau_{1})^{6}\sum_{j=1}^{3}\theta_{j}^{2}}dy_{1}dy_{2}dy_{3}d\tau_{1}+
\\&&(\frac{1}{\sqrt{\pi}})^{3}\int_{0}^{t}\int_{R^{3}}
\frac{v_{m}^{(3)}}{\sqrt{\mu}}(\sqrt{t-\tau_{1}})^{45} \theta_{i}
(\theta_{1}^{2}+\theta_{2}^{2}+\theta_{3}^{2})e^{-(t-\tau_{1})^{8}\sum_{j=1}^{3}\theta_{j}^{2}}dy_{1}dy_{2}dy_{3}d\tau_{1}
,\end{eqnarray*}where$$v_{m}^{(3)}=v_{m}(\tau_{1},
x_{j}+2\sqrt{\mu}(\sqrt{(t-\tau_{1})})^{9}\theta_{j},j=1,2,3,)\ .$$
 We assume\ $\varphi_{4}(t,\ \tau_{1})\ ,\ \varphi_{5}(t,\ \tau_{1})$\
as follows,\begin{eqnarray*}&&\varphi_{4}(t,\
\tau_{1})=\int_{R^{3}}(\sqrt{t-\tau_{1}})^{21}\theta_{i}\
e^{-(t-\tau_{1})^{6}(\theta_{1}^{2}+\theta_{2}^{2}+\theta_{3}^{2})}
d\theta_{1}d\theta_{2}d\theta_{3}\ ,\\&&\varphi_{5}(t,\
\tau_{1})=\int_{R^{3}}(\sqrt{t-\tau_{1}})^{45}\theta_{i}(\theta_{1}^{2}+\theta_{2}^{2}+\theta_{3}^{2})\
e^{-(t-\tau_{1})^{8}(\theta_{1}^{2}+\theta_{2}^{2}+\theta_{3}^{2})}
d\theta_{1}d\theta_{2}d\theta_{3}\ .\end{eqnarray*}We can see that\
$\varphi_{4}(t,\ \tau_{1}),\ \varphi_{5}(t,\ \tau_{1})$\ are
continuous about\ $t,\ \tau_{1}$\ on the region\ $0\leq \tau_{1}\leq
t\ ,\ 0\leq t\leq T\ ,$ hence they are bounded on the region\ $0\leq
\tau_{1}\leq t\ ,\ 0\leq t\leq T$\ . We assume there exist\
$M^{\prime\prime\prime}>0$\ , such that\ $|\varphi_{4}(t,\
\tau_{1})|\leq M^{\prime\prime\prime}\ ,\ |\varphi_{5}(t,\
\tau_{1})|\leq M^{\prime\prime\prime}$\ , where\ $0\leq \tau_{1}\leq
t\ ,\ 0\leq t\leq T\ .$\ Hence we can
get$$|\cfrac{\partial\cfrac{\partial h_{2m}}{\partial
x_{i}}}{\partial t}\ |\leq
(\frac{1}{\sqrt{\pi}})^{3}\cfrac{MM(K_{1})}{\sqrt{\mu}}\
M^{\prime\prime}+(\frac{1}{\sqrt{\pi}})^{3}\cfrac{7}{2\sqrt{\mu}}\
TMM(K_{1})M^{\prime\prime\prime}\ .$$ If we let$$ M_{T,\ 1}=\max\{T\
,\ 1+(\frac{1}{\sqrt{\pi}})^{3}\cfrac{5}{2}\ TM^{\prime}\ ,\
(\frac{1}{\sqrt{\pi}})^{3}\cfrac{T}{\sqrt{\mu}}\ M^{\prime\prime}\
,\ (\frac{1}{\sqrt{\pi}})^{3}\cfrac{1}{\sqrt{\mu}}\
M^{\prime\prime}+(\frac{1}{\sqrt{\pi}})^{3}\cfrac{7}{2\sqrt{\mu}}\
TM^{\prime\prime\prime}\}\ ,$$then we can get
$$\max_{1\leq
i\leq3}\{|h_{2}|\ ,\ |\cfrac{\partial h_{2}}{\partial t}\ |\ ,\
|\cfrac{\partial h_{2}}{\partial x_{i}}\ |\ ,\
|\cfrac{\partial\cfrac{\partial h_{2}}{\partial x_{i}}}{\partial t}\
|\}\leq\ MM_{T,\ 1}M(K_{1})\ .$$ Hence the lemma is true.
\begin{corollary} \label{corollary1}$\exists\ M_{T,\ 2}>0\ ,\ \mbox{and\ $M_{T,\ 2}$\ is independent with\ $M\ ,\ K_{1}$}\ ,\forall\ h_{1}\in\
\Omega_{C}$\ , we can get the following, $$\max_{i,\ j,\ k,\ l}\{\
|\cfrac{\partial^{2} h_{2}}{\partial x_{i}\partial x_{j}}\ |\ ,\
|\cfrac{\partial\cfrac{\partial^{2} h_{2}}{\partial x_{i}\partial
x_{j}}}{\partial t}\ |\ ,\ |\cfrac{\partial^{3} h_{2}}{\partial
x_{i}\partial x_{j}\partial x_{k}}\ |\ ,\
|\cfrac{\partial\cfrac{\partial^{3} h_{2}}{\partial x_{i}\partial
x_{j}\partial x_{k}}}{\partial t}\ |\ ,\ |\cfrac{\partial^{4}
h_{2}}{\partial x_{i}\partial x_{j}\partial x_{k}\partial x_{l}}\
|\}\leq\ MM_{T,\ 2}M(K_{1})\ ,$$ \mbox{where}\begin{eqnarray*}&&
h_{1}=(h_{1m})_{9\times1}\ ,\ h_{2}=h_{2}(t,\ M)=(h_{2m}(t,\
M))_{9\times1}\ ,\ |\cfrac{\partial^{2} h_{2}}{\partial
x_{i}\partial x_{j}}\ |=\max_{1\leq m\leq 9}\max_{(t,\ M)\in\
K_{1}^{\prime}}|\cfrac{\partial^{2} h_{2m}(t,\ M)}{\partial
x_{i}\partial x_{j}}\ |\ ,\\&& |\cfrac{\partial\cfrac{\partial^{2}
h_{2}}{\partial x_{i}\partial x_{j}}}{\partial t}\ |=\max_{1\leq
m\leq 9}\max_{(t,\ M)\in\
K_{1}^{\prime}}|\cfrac{\partial\cfrac{\partial^{2} h_{2m}(t,\
M)}{\partial x_{i}\partial x_{j}}}{\partial t}\ |\ ,\
|\cfrac{\partial^{3} h_{2}}{\partial x_{i}\partial x_{j}\partial
x_{k}}\ |=\max_{1\leq m\leq 9}\max_{(t,\ M)\in\
K_{1}^{\prime}}|\cfrac{\partial^{3} h_{2m}(t,\ M)}{\partial
x_{i}\partial x_{j}\partial x_{k}}\ |\ ,\\&&|\cfrac{\partial
\cfrac{\partial^{3} h_{2}}{\partial x_{i}\partial x_{j}\partial
x_{k}}}{\partial t}\ |=\max_{1\leq m\leq 9}\max_{(t,\ M)\in\
K_{1}^{\prime}}|\cfrac{\partial\cfrac{\partial^{3} h_{2m}(t,\
M)}{\partial x_{i}\partial x_{j}\partial x_{k}}}{\partial t}\ |\ ,\
M(K_{1})=\max_{M_{0}\in
K_{1}}\frac{1}{4\pi}\int_{K_{1}}\frac{1}{r_{MM_{0}}}
dx_{1}dx_{2}dx_{3}\ ,\\&&|\cfrac{\partial^{4} h_{2}}{\partial
x_{i}\partial x_{j}\partial x_{k}\partial x_{l}}\ |=\max_{1\leq
m\leq 9}\max_{(t,\ M)\in\ K_{1}^{\prime}}|\cfrac{\partial^{4}
h_{2m}(t,\ M)}{\partial x_{i}\partial x_{j}\partial x_{k}\partial
x_{l}}\ |\ ,\ 1\leq i,\ j,\ k,\ l\leq3\ .\end{eqnarray*}
\end{corollary}
\begin{corollary} \label{corollary1}$\exists\ M_{T,\ 3}>0\ ,\ \mbox{and\ $M_{T,\ 3}$\ is independent with\ $M\ ,\ K_{1}$}\ , \forall\ h_{1}\in\
\Omega_{C}$\ , we can get the following,
$$\max\{\ |W_{jl}(h_{2})|\ ,
\ |\cfrac{\partial W_{jl}(h_{2})}{\partial t}|\ ,\ |\cfrac{\partial
W_{jl}(h_{2})}{\partial x_{i}}|\ ,\ 1\leq j\leq9\ ,\ 1\leq l\leq2\
,\ 1\leq i\leq 3\ .\}\leq\ MM_{T,\ 3}M(K_{1})\ .$$\end{corollary}
From the corollary 3.4 , we can
get$$|W_{jl}(h_{2})(x_{t})-W_{jl}(h_{2})(x_{t}^{\prime})|\leq
4MM_{T,\ 3}M(K_{1})|x_{t}-x_{t}^{\prime}|\leq MM_{T,\
4}M(K_{1})|x_{t}-x_{t}^{\prime}|^{\alpha}\ ,\ \forall\ x_{t}\ ,\
x_{t}^{\prime}\in K_{1}^{\prime}\ ,$$where$$1\leq j\leq9\ ,\ 1\leq
l\leq2\ ,\ M_{T,\ 4}=4M_{T,\ 3}\max_{x_{t},\ x_{t}^{\prime}\in
K_{1}^{\prime}}|x_{t}-x_{t}^{\prime}|^{1-\alpha}\ .$$Now we can
get\begin{eqnarray*}
|g_{j}(h_{1})|&=&|F^{-1}(\alpha_{j1}^{T}Y_{1})F^{-1}(\alpha_{j2}^{T}Y_{1})+F^{-1}(\alpha_{j1}^{T}Y_{1})W_{j2}(h_{2})+
\\&&F^{-1}(\alpha_{j2}^{T}Y_{1})W_{j1}(h_{2})+W_{j1}(h_{2})W_{j2}(h_{2})|
\\&\leq &\theta M+2\theta M^{2}\ M_{T,\ 3}M(K_{1})+[MM_{T,\ 3}M(K_{1})]^{2}
\ ,\ 1\leq j\leq 9\ . \end{eqnarray*}And from
\begin{eqnarray*}&&g_{j}(h_{1})(x_{t})-g_{j}(h_{1})(x_{t}^{\prime})
\\&=&F^{-1}(\alpha_{j1}^{T}Y_{1})(x_{t})F^{-1}(\alpha_{j2}^{T}Y_{1})(x_{t})+F^{-1}(\alpha_{j1}^{T}Y_{1})(x_{t})W_{j2}(h_{2})(x_{t})+
\\&&F^{-1}(\alpha_{j2}^{T}Y_{1})(x_{t})W_{j1}(h_{2})(x_{t})+W_{j1}(h_{2})(x_{t})W_{j2}(h_{2})(x_{t})-
\\&&[F^{-1}(\alpha_{j1}^{T}Y_{1})(x_{t}^{\prime})F^{-1}(\alpha_{j2}^{T}Y_{1})(x_{t}^{\prime})+F^{-1}(\alpha_{j1}^{T}Y_{1})(x_{t}^{\prime})W_{j2}(h_{2})(x_{t}^{\prime})+
\\&&F^{-1}(\alpha_{j2}^{T}Y_{1})(x_{t}^{\prime})W_{j1}(h_{2})(x_{t}^{\prime})+W_{j1}(h_{2})(x_{t}^{\prime})W_{j2}(h_{2})(x_{t}^{\prime})]
\end{eqnarray*}\begin{eqnarray*}&=&F^{-1}(\alpha_{j1}^{T}Y_{1})(x_{t})F^{-1}(\alpha_{j2}^{T}Y_{1})(x_{t})-F^{-1}(\alpha_{j1}^{T}Y_{1})(x_{t}^{\prime})F^{-1}(\alpha_{j2}^{T}Y_{1})(x_{t}^{\prime})
+\\&&F^{-1}(\alpha_{j1}^{T}Y_{1})(x_{t})(W_{j2}(h_{2})(x_{t})-W_{j2}(h_{2})(x_{t}^{\prime}))+(F^{-1}(\alpha_{j1}^{T}Y_{1})(x_{t})-F^{-1}(\alpha_{j1}^{T}Y_{1})(x_{t}^{\prime}))W_{j2}(h_{2})(x_{t}^{\prime})
+\\&&F^{-1}(\alpha_{j2}^{T}Y_{1})(x_{t})(W_{j1}(h_{2})(x_{t})-W_{j1}(h_{2})(x_{t}^{\prime}))+(F^{-1}(\alpha_{j2}^{T}Y_{1})(x_{t})-F^{-1}(\alpha_{j2}^{T}Y_{1})(x_{t}^{\prime}))W_{j1}(h_{2})(x_{t}^{\prime})+
\\&&W_{j1}(h_{2})(x_{t})(W_{j2}(h_{2})(x_{t})-W_{j2}(h_{2})(x_{t}^{\prime}))+(W_{j1}(h_{2})(x_{t})-W_{j1}(h_{2})(x_{t}^{\prime}))W_{j2}(h_{2})(x_{t}^{\prime})\
,\end{eqnarray*} we can
get\begin{eqnarray*}&&|g_{j}(h_{1})(x_{t})-g_{j}(h_{1})(x_{t}^{\prime})|\leq\\&&
\{C_{1}+2[\theta M^{2}\ M_{T,\ 4}M(K_{1}) +C_{1}MM_{T,\
3}M(K_{1})+M^{2}M_{T,\ 3}M_{T,\
4}M(K_{1})^{2}]\}|x_{t}-x_{t}^{\prime}|^{\alpha}\ ,\ 1\leq j\leq9\
.\end{eqnarray*} Hence\ $g(h_{1})=(g_{j}(h_{1}))_{9\times1}\in
\Omega_{C}$\ , if the followings are true.\[\begin{cases}\theta
M+2\theta M^{2}\ M_{T,\ 3}M(K_{1})+[MM_{T,\ 3}M(K_{1})]^{2}\leq M\
,\\C_{1}+2[\theta M^{2}\ M_{T,\ 4}M(K_{1}) +C_{1}MM_{T,\
3}M(K_{1})+M^{2}M_{T,\ 3}M_{T,\ 4}M(K_{1})^{2}]\leq C\
.\end{cases}\]If we let\ $C$\ is big enough, and from the assumption
1.1, if we let$$M=\cfrac{1-\theta}{M_{T,\ 3}M(K_{1})(2\theta+M_{T,\
3}M(K_{1}))}\ ,$$ then (3.8) will be true.
\\At last we need to prove the mapping\
$T :\ h_{1}\rightarrow\ g(h_{1})$\ is continuous about\ $h_{1}$\ in
the\ $
\Omega_{C}$\ .\\
We assume\ $|h_{1}^{(n)}-h_{1}^{*}|\rightarrow 0$\ , when\
$n\rightarrow +\infty$\ , and\ $h_{1}^{(n)}\ ,\ h_{1}^{*}\in
\Omega_{C}\ ,\ n\geq 1$\ , from the lemma 3.1 we know that there
exist\ $h_{2}^{(n)}=h_{2}^{(n)}I_{K_{1}^{\prime}}\in C^{1}[0,\
T]\bigcap C^{\infty}(K_{1})\ ,\
h_{2}^{*}=h_{2}^{*}I_{K_{1}^{\prime}}\in C^{1}[0,\ T]\bigcap
C^{\infty}(K_{1})\ ,$\\such that\ $F(h_{1}^{(n)})=a^{2}b\tau
F(h_{2}^{(n)})\ ,\ F(h_{1}^{*})=a^{2}b\tau F(h_{2}^{*})$\ ,
where\begin{eqnarray*}&&v^{(n)}(t,M_{0})=-\frac{1}{4\pi}\int_{K_{1}}\frac{h_{1}^{(n)}(t,M)}{r_{MM_{0}}}
dx_{1}dx_{2}dx_{3}\ ,\ \mbox{and}\ M=(x_{1},x_{2},x_{3})\ ,\
M_{0}=(x_{10},x_{20},x_{30})\
,\\&&h_{2}^{(n)}(t,M)=(\frac{1}{2\sqrt{\pi\mu}})^{3}\int_{0}^{t}\int_{R^{3}}\frac{v^{(n)}(\tau_{1},
y_{1}, y_{2} , y_{3})}{(\sqrt{t-\tau_{1}})^{3}}\
e^{-\frac{(x_{1}-y_{1})^{2}+(x_{2}-y_{2})^{2}+(x_{3}-y_{3})^{2}}{4\mu(t-\tau_{1})}}
dy_{1}dy_{2}dy_{3}d\tau_{1}\
,\\&&v^{*}(t,M_{0})=-\frac{1}{4\pi}\int_{K_{1}}\frac{h_{1}^{*}(t,M)}{r_{MM_{0}}}
dx_{1}dx_{2}dx_{3}\ ,\ \mbox{and}\ M=(x_{1},x_{2},x_{3})\ ,\
M_{0}=(x_{10},x_{20},x_{30})\
,\\&&h_{2}^{*}(t,M)=(\frac{1}{2\sqrt{\pi\mu}})^{3}\int_{0}^{t}\int_{R^{3}}\frac{v^{*}(\tau_{1},
y_{1}, y_{2} , y_{3})}{(\sqrt{t-\tau_{1}})^{3}}\
e^{-\frac{(x_{1}-y_{1})^{2}+(x_{2}-y_{2})^{2}+(x_{3}-y_{3})^{2}}{4\mu(t-\tau_{1})}}
dy_{1}dy_{2}dy_{3}d\tau_{1}\ .
\end{eqnarray*} From the Lebesgue dominated convergence theorem, we
can get\ $|v^{(n)}-v^{*}|\rightarrow 0$\ , when\ $n\rightarrow
+\infty\ .$ And we know the partial derivation of\ $h_{2}^{(n)}$\
which is no more than the third order, only do the partial
derivation with the variables\ $x_{1}\ ,\ x_{2}\ ,\ x_{3}\ ,$\ can
pass through the integral. Again from the Lebesgue dominated
convergence theorem, we can get when\ $n\rightarrow +\infty$\
,\begin{eqnarray*}&&|h_{2}^{(n)}-h_{2}^{*}|\rightarrow 0\ ,\
|\cfrac{\partial h_{2}^{(n)}}{\partial x_{i}}-\cfrac{\partial
h_{2}^{*}}{\partial x_{i}}\ |\rightarrow 0\ ,\ |\cfrac{\partial^{2}
h_{2}^{(n)}}{\partial x_{i}\partial x_{j}}-\cfrac{\partial^{2}
h_{2}^{*}}{\partial x_{i}\partial x_{j}}\ |\rightarrow 0\ ,\\&&
|\cfrac{\partial^{3} h_{2}^{(n)}}{\partial x_{i}\partial
x_{j}\partial x_{k}}-\cfrac{\partial^{3} h_{2}^{*}}{\partial
x_{i}\partial x_{j}\partial x_{k}}\ |\rightarrow 0\ ,\ 1\leq i\ ,\
j\ ,\ k\leq3\ .\end{eqnarray*} Hence we can get when\ $n\rightarrow
+\infty$\ ,$$ |W_{j1}(h_{2}^{(n)})-W_{j1}(h_{2}^{*})|\rightarrow 0\
,\ |W_{j1}(h_{2}^{(n)})-W_{j1}(h_{2}^{*})|\rightarrow 0\ ,\ 1\leq
j\leq 9\ .$$ Now we can get\begin{eqnarray*}&&
|g_{j}(h_{1}^{(n)})-g_{j}(h_{1}^{*})|\\&=&|F^{-1}(\alpha_{j1}^{T}Y_{1})F^{-1}(\alpha_{j2}^{T}Y_{1})+F^{-1}(\alpha_{j1}^{T}Y_{1})W_{j2}(h_{2}^{(n)})+
F^{-1}(\alpha_{j2}^{T}Y_{1})W_{j1}(h_{2}^{(n)})+W_{j1}(h_{2}^{(n)})W_{j2}(h_{2}^{(n)})-
\\&&[F^{-1}(\alpha_{j1}^{T}Y_{1})F^{-1}(\alpha_{j2}^{T}Y_{1})+F^{-1}(\alpha_{j1}^{T}Y_{1})W_{j2}(h_{2}^{*})+
F^{-1}(\alpha_{j2}^{T}Y_{1})W_{j1}(h_{2}^{*})+W_{j1}(h_{2}^{*})W_{j2}(h_{2}^{*})]|
\\&\leq&|F^{-1}(\alpha_{j1}^{T}Y_{1})||W_{j2}(h_{2}^{(n)})-W_{j2}(h_{2}^{*})|+|F^{-1}(\alpha_{j2}^{T}Y_{1})||W_{j1}(h_{2}^{(n)})-W_{j1}(h_{2}^{*})|
+\\&&|W_{j1}(h_{2}^{(n)})-W_{j1}(h_{2}^{*})||W_{j2}(h_{2}^{(n)})|+|W_{j1}(h_{2}^{*})||W_{j2}(h_{2}^{(n)})-W_{j2}(h_{2}^{*})|
\\&\leq& [\theta M + MM_{T,\ 3}M(K_{1})]\ [\ |W_{j1}(h_{2}^{(n)})-W_{j1}(h_{2}^{*})|+|W_{j2}(h_{2}^{(n)})-W_{j2}(h_{2}^{*})|\ ]\
,\end{eqnarray*}this means that\
$|g_{j}(h_{1}^{(n)})-g_{j}(h_{1}^{*})|\rightarrow 0\ ,\ 1\leq
j\leq9$\ , when\ $n\rightarrow +\infty$\ . \\Hence\ $T$\ is
continuous about\ $h_{1}$\ in the\ $ \Omega_{C}$\ . According to the
Schauder fixed-point theorem we can learn that there exists\
$h_{1}\in \Omega_{C}$\ , such that\ $h_{1}=g(h_{1})$\ . And\
$\forall\ j\ ,\ 1\leq j\leq9$\ , we know that
$$g_{j}(h_{1})=F^{-1}(\alpha_{j1}^{T}Y_{1})F^{-1}(\alpha_{j2}^{T}Y_{1})+F^{-1}(\alpha_{j1}^{T}Y_{1})W_{j2}(h_{2})+
F^{-1}(\alpha_{j2}^{T}Y_{1})W_{j1}(h_{2})+W_{j1}(h_{2})W_{j2}(h_{2})\
,$$where\ $h_{2}=h_{2}I_{K_{1}^{\prime}}\in C^{1}[0,\ T]\bigcap
C^{\infty}(K_{1})\ ,$\ and\ $W_{j1}(h_{2})\ ,\ W_{j2}(h_{2})\ ,\
1\leq j\leq9\ ,$\ are the functions to do the partial derivation
with the components of\ $h_{2}$\ no more than the third order and
their linear combination, moreover only do the partial derivation
with the variables\ $x_{1}\ ,\ x_{2}\ ,\ x_{3}$\ . And the
components of\ $F^{-1}(\alpha_{j1}^{T}Y_{1})\ ,\
F^{-1}(\alpha_{j2}^{T}Y_{1})\ ,\ 1\leq j\leq9$\ , are all in the\ \
$H[F^{-1}(Y_{1})]$\ , moreover
$$H[F^{-1}(Y_{1})]=H[F^{-1}(Y_{1})]I_{K_{1}^{\prime}}\in
C^{1}(K_{1}^{\prime})\ .$$ This means that\ $g_{j}(h_{1})\in
C^{1}(K_{1}^{\prime})\ ,\ 1\leq j\leq 9 $\ , hence the fixed-point
$$h_{1}=g(h_{1})=(g_{j}(h_{1}))_{9\times1}\in C^{1}(K_{1}^{\prime})\
.$$ At last if we let\ $Z_{1}=F(h_{1})$\ , where\ $h_{1}$\ is the
fixed-point of\ $g(h_{1})$\ , then\ $Z_{1}\in \Omega_{2}$\ , and\
$Z_{1}$\ is the fixed point of\ $f(Z_{1})$\ . This is also to say
that the smoothing solution for the Navier-Stokes equations is
globally exist, and we can get the
following.\begin{eqnarray*}u_{1}&=&e_{3}^{T}Z=e_{3}^{T}[F^{-1}(Y_{1}+A_{1}(\eta)Z_{1})]=w_{17}+w_{27}
\\&=&\cfrac{\partial^{2}F_{21}}{\partial
x_{1}\partial x_{2}}+\cfrac{\partial^{2}F_{31}}{\partial
x_{1}\partial x_{3}}-\cfrac{\partial^{2}F_{11}}{\partial
x_{2}^{2}}-\cfrac{\partial^{2}F_{11}}{\partial
x_{3}^{2}}+\cfrac{\partial^{2}h_{31}}{\partial
x_{2}^{2}}+\cfrac{\partial^{2}h_{31}}{\partial
x_{3}^{2}}-\cfrac{\partial^{2}h_{32}}{\partial x_{1}\partial
x_{2}}-\cfrac{\partial^{2}h_{33}}{\partial x_{1}\partial x_{3}}\ ,
\end{eqnarray*}\begin{eqnarray*}u_{2}&=&e_{7}^{T}Z=e_{7}^{T}[F^{-1}(Y_{1}+A_{1}(\eta)Z_{1})]=w_{111}+w_{211}
\\&=&\cfrac{\partial^{2}F_{11}}{\partial x_{1}\partial
x_{2}}+\cfrac{\partial^{2}F_{31}}{\partial x_{2}\partial
x_{3}}-\cfrac{\partial^{2}F_{21}}{\partial
x_{1}^{2}}-\cfrac{\partial^{2}F_{21}}{\partial
x_{3}^{2}}+\cfrac{\partial^{2}h_{32}}{\partial
x_{1}^{2}}+\cfrac{\partial^{2}h_{32}}{\partial
x_{3}^{2}}-\cfrac{\partial^{2}h_{31}}{\partial x_{1}\partial
x_{2}}-\cfrac{\partial^{2}h_{33}}{\partial x_{2}\partial x_{3}}\
,\\u_{3}&=&e_{11}^{T}Z=e_{11}^{T}[F^{-1}(Y_{1}+A_{1}(\eta)Z_{1})]=w_{115}+w_{215}
\\&=&\cfrac{\partial^{2}F_{11}}{\partial
x_{1}\partial x_{3}}+\cfrac{\partial^{2}F_{21}}{\partial
x_{2}\partial x_{3}}-\cfrac{\partial^{2}F_{31}}{\partial
x_{1}^{2}}-\cfrac{\partial^{2}F_{31}}{\partial
x_{2}^{2}}+\cfrac{\partial^{2}h_{33}}{\partial
x_{1}^{2}}+\cfrac{\partial^{2}h_{33}}{\partial
x_{2}^{2}}-\cfrac{\partial^{2}h_{31}}{\partial x_{1}\partial
x_{3}}-\cfrac{\partial^{3}h_{32}}{\partial x_{2}\partial x_{3}}\
,\\p&=&e_{16}^{T}Z=e_{16}^{T}[F^{-1}(Y_{1}+A_{1}(\eta)Z_{1})]=w_{116}+w_{216}
\\&=&\cfrac{1}{\tau}\ [\ \mu
\triangle(\cfrac{\partial F_{11}}{\partial x_{1}}+\cfrac{\partial
F_{21}}{\partial x_{2}}+\cfrac{\partial F_{31}}{\partial
x_{3}})-\cfrac{\partial^{2} F_{11}}{\partial t\partial
x_{1}}-\cfrac{\partial^{2} F_{21}}{\partial t\partial
x_{2}}-\cfrac{\partial^{2} F_{31}}{\partial t\partial x_{3}}\
]+\\&&\cfrac{1}{\tau}\ (\cfrac{\partial^{2}h_{31}}{\partial
t\partial x_{1}}+\cfrac{\partial^{2}h_{32}}{\partial t\partial
x_{2}}+\cfrac{\partial^{2}h_{33}}{\partial t\partial
x_{3}})-\cfrac{\mu}{\tau}\ (\cfrac{\partial\triangle
h_{31}}{\partial x_{1}}+\cfrac{\partial\triangle h_{32}}{\partial
x_{2}}+\cfrac{\partial\triangle h_{33}}{\partial x_{3}})\ ,
\end{eqnarray*}where\begin{eqnarray*}&&v(t,M_{0})=-\frac{1}{4\pi}\int_{K_{1}}\frac{h_{1}(t,M)}{r_{MM_{0}}}
dx_{1}dx_{2}dx_{3}\ ,\ \mbox{and}\ M=(x_{1},x_{2},x_{3})\ ,\
M_{0}=(x_{10},x_{20},x_{30})\
,\\&&h_{2}(t,M)=(\frac{1}{2\sqrt{\pi\mu}})^{3}\int_{0}^{t}\int_{R^{3}}\frac{v(\tau_{1},
y_{1}, y_{2} , y_{3})}{(\sqrt{t-\tau_{1}})^{3}}\
e^{-\frac{(x_{1}-y_{1})^{2}+(x_{2}-y_{2})^{2}+(x_{3}-y_{3})^{2}}{4\mu(t-\tau_{1})}}
dy_{1}dy_{2}dy_{3}d\tau_{1}\ ,\\&&h_{31}=h_{21}+h_{22}+h_{23}\ ,\
h_{32}=h_{24}+h_{25}+h_{26}\ ,\ h_{33}=h_{27}+h_{28}+h_{29}\
,\\&&F_{j1}(t,M)=(\frac{1}{2\sqrt{\pi\mu}})^{3}\int_{0}^{t}\int_{R^{3}}\frac{F_{jv}(\tau_{1},
y_{1}, y_{2} , y_{3})}{(\sqrt{t-\tau_{1}})^{3}}\
e^{-\frac{(x_{1}-y_{1})^{2}+(x_{2}-y_{2})^{2}+(x_{3}-y_{3})^{2}}{4\mu(t-\tau_{1})}}
dy_{1}dy_{2}dy_{3}d\tau_{1}\
,\\&&F_{jv}(t,M_{0})=-\frac{1}{4\pi}\int_{K_{1}}\frac{F_{j}(t,M)}{r_{MM_{0}}}
dx_{1}dx_{2}dx_{3}\ ,\ 1\leq\ j\leq 3\ ,\end{eqnarray*}and\ $h_{1}$\
is the fixed-point of\ $g(h_{1})$\ , where\
$g(h_{1})=(g_{j}(h_{1}))_{9\times1}$\ , and$$
g_{j}(h_{1})=F^{-1}(\alpha_{j1}^{T}Y_{1})F^{-1}(\alpha_{j2}^{T}Y_{1})+F^{-1}(\alpha_{j1}^{T}Y_{1})W_{j2}(h_{2})+
F^{-1}(\alpha_{j2}^{T}Y_{1})W_{j1}(h_{2})+W_{j1}(h_{2})W_{j2}(h_{2})\
,$$ $ 1\leq\ j\leq\ 9$\ . And because\ $Z_{1}=F(h_{1})\in
\Omega_{2}$\ , we can get
$$u_{j}\in C^{2}[0,\ T]\bigcap C^{\infty}(K_{1})\ ,\ 1\leq j\leq 3\
,\ p\in C^{1}[0,\ T]\bigcap C^{\infty}(K_{1})\ .$$ \\Finally we say
our sincerely thanks to Schauder according to at least two points as
following.\\(1)The Schauder fixed-point theorem.\\(2)If\ $h_{1}$\ is
H\"{o}lder continuous,\ $\partial K_{1}$\ satisfies the exterior
ball condition, then the smoothing solution of Poisson's equation
exists. This is also proved by Schauder, and it leads to that the
Newtonian potential of\ $h_{1}$\
,$$\cfrac{-1}{4\pi}\int_{K_{1}}\cfrac{h_{1}}{r_{MM_{0}}}dx_{1}dx_{2}dx_{3}$$is
a solution of Poisson's equation.\\If there is no important help
from Schauder, we can not see so far with the Navier-Stokes
equations. Maybe you will ask why not to do the Laplace
transformation with the variable\ $t$\ just as usual. It is very
difficult to discuss the fixed-point of\ $f(Z_{1})$\ if we do that.
We have tried to do the Laplace transformation with the variable\
$t$\ , but we can not discuss the inverse Laplace transformation in\
$f(Z_{1})$\ . This is the reason why we also do the Fourier
transformation with the variable\ $t$\ .

\end{CJK*}

\end{document}